\documentclass[11pt] {article}
\usepackage[round]{natbib}
\usepackage{amssymb}
\usepackage{amsmath}
\usepackage{fullpage}
\usepackage{booktabs}  
\usepackage{multirow}  
\usepackage{amsthm}
\usepackage{graphicx}
\usepackage{caption}
\usepackage{xcolor}
\usepackage{enumitem}
 \usepackage{algorithm} 
\usepackage{algpseudocode}
\usepackage{subcaption}
\usepackage{authblk}
\usepackage{mathrsfs}
\usepackage{bm}
\usepackage{bbm}
\newcommand{\bvec}[1]{\bm{#1}}
\usepackage{comment}
\algnewcommand\algorithmicinput{\textbf{INPUT:}}
\algnewcommand\INPUT{\item[\algorithmicinput]}  
\algnewcommand\algorithmicoutput{\textbf{OUTPUT:}}
\algnewcommand\OUTPUT{\item[\algorithmicoutput]}
\usepackage{hyperref}[]
\hypersetup{
    colorlinks=true,
    linkcolor=blue,
    filecolor=magenta,      
    urlcolor=blue,
      citecolor=blue,
}
\usepackage{cleveref}
 \usepackage{xcolor}

\newcommand{\mi}{{l}}

\newcommand{\vx}{x}

\newcommand{\dd}{\mathrm{d}}

\DeclareMathOperator*{\rank}{Rank}

\DeclareMathOperator*{\col}{Col} 
\DeclareMathOperator*{\row}{Row} 
 \DeclareMathOperator*{\Span}{Span}
  \DeclareMathOperator*{\svd}{SVD} 
\DeclareMathOperator*{\op}{2}

\DeclareMathOperator*{\bigO}{O}

 \newcommand{\p}{\mathbb P}

\DeclareMathOperator{\var}{Var}

   \newcommand{\lt}{{  { L}_2 } }

    \newcommand{\E}{ { \mathbb{E} } }

 \newcommand{\alphaspace}{\beta_0}

\newtheorem{assumption}{Assumption}
\newtheorem{lemma}{Lemma}

\newtheorem{corollary}{Corollary}
\newtheorem{definition}{Definition}
\newtheorem{remark}{Remark}
\newtheorem{theorem}{Theorem}

\date{\vspace{-5ex}}
\title{Tensor Density Estimator by Convolution-Deconvolution}
\author[1]{Yifan Peng}
\author[2]{Siyao Yang\footnote{Corresponding author}}
\author[3]{Yuehaw Khoo}
\author[4]{Daren Wang}
\affil[1,2]{Committee on Computational and Applied Mathematics, University of Chicago}
\affil[3]{Committee on Computational and Applied Mathematics and  Department of Statistics, University of Chicago}
\affil[4]{Department of Mathematics, University of California, San Diego}

\begin{document}
\maketitle

\begin{abstract}
We propose a linear algebraic framework for performing density estimation. It consists of three simple steps: convolving the empirical distribution with certain smoothing kernels to remove the exponentially large variance; compressing the empirical distribution after convolution as a tensor train, with efficient tensor decomposition algorithms; and finally, applying a deconvolution step to recover the estimated density from such tensor-train representation. Numerical results demonstrate the high accuracy and efficiency of the proposed methods.
\end{abstract}

\section{Introduction}
Density estimation is a fundamental problem in data science with wide-ranging applications in fields such as science (\cite{cranmer2001kernel,chen2017tutorial}), engineering (\cite{chaudhuri2014consistent}), and economics (\cite{zambom2013review}). More precisely, given $N$ data points $\{x^{(i)}\}_{i=1}^N \sim p^*$, where each $x^{(i)}$ is an independent and identically distributed (i.i.d.) sample from an unknown ground truth distribution $p^*:\mathbb{R}^{d}\rightarrow \mathbb{R}$, the goal is to estimate $p^*$ based on empirical distribution 
\begin{equation}
\hat{p}(x)=\frac{1}{N}\sum_{i=1}^N \delta_{x^{(i)}}(x) 
\end{equation}
where $\delta_{x^{(i)}}(x)$ is the Dirac delta function centered at $x^{(i)}$. In general, one cannot use the empirical distribution as an estimator of the density $p^*$, due to the lack of absolutely continuity or exponentially large variance.

\subsection{Our Contributions}
\begin{figure}[ht!]
\centering
\includegraphics[width=0.6\textwidth]{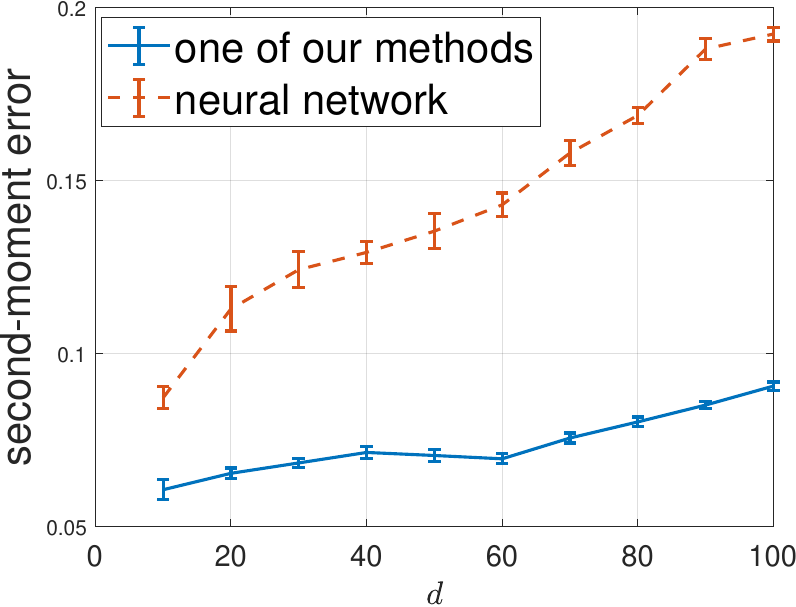} 
\caption{Density estimation for data sampled from a Boltzmann distribution. Figure visualizes the second-moment error for different dimensionality $d$. We compare the result of one of our proposed tensor-train-based method called TT-SVD-kn (see Section \ref{sec:tt-svd-kn}) with that of one neural network approach (Masked Autoregressive Flow).}
\label{fig:intro}
\end{figure}

In this work, we introduce a novel framework for high-dimensional density estimation entirely based on linear algebra. Our goal is to achieve overall computational complexity that is linear in dimensionality. This is accomplished by viewing the density estimation process as applying a low-dimensional orthogonal projector in function space onto the noisy $\hat{p}$. Due to the dimension reduction, the exponentially large variance in $\hat{p}$ is suppressed. Additionally, we leverage tensor-train (TT) algebra to perform this projection with linear computational complexity in terms of both dimensionality and sample size. This approach also yields a final density estimator in the form of a tensor train. Furthermore, we theoretically show that, in addition to achieving linear scaling in computational complexity, the estimation error also enjoys linear scaling under certain mild regularity assumptions. Finally, we validate our method on Boltzmann sampling and real-world generative-modeling tasks. The numerical results demonstrate that our estimator attains high accuracy, comparable to or even surpassing that of deep learning–based methods in several cases. A representative result is shown in Figure~\ref{fig:intro}, with full details provided in Section~\ref{sec:numerical}. 

\subsection{Main idea}\label{sec: main idea}
In this subsection, we outline the main idea of the proposed method. Our method is based on a crucial assumption about the true density: 
\begin{equation}\label{eq:two cluster ground truth} 
p^*(x) = \mu(x)  \left(1+ \sum^d_{j_1=1} \sum_{l_1 =2 }^n b^{j_1}_{l_1} \phi^{j_1}_{l_1}(x_{j_1}) + \sum_{(j_1,j_2)\in \binom{d}{2}} \sum_{l_1,l_2 = 2}^n b^{j_1 j_2}_{l_1 l_2} \phi^{j_1}_{l_1}(x_{j_1}) \phi^{j_2}_{l_2}(x_{j_2})\right), \end{equation}
where $\mu(x) = \mu_1(x_1)\cdots \mu_d(x_d)$ is a mean-field density, i.e., a separable density. $\{\phi^{j}_{l}\}_{l=1}^{n}$ are orthogonal basis functions with respect to $\mu_j$ and $\phi^{j}_{1} \equiv 1$, i.e.,
\begin{equation}\label{eq:orthogonal poly def} 
\int \phi^{j}_{l_1}(x_{j})\phi^{j}_{l_2}(x_{j}) \mu_j(d x_j) = \delta_{l_1,l_2},
\end{equation}
and $b^{j_1}_{l_1},b^{j_1j_2}_{l_1l_2}$ are the corresponding coefficients. The indices $(j_1,j_2)\in \binom{d}{2}$ denote all unordered pairs of distinct dimensions. In other words, we view $p^*$ as a perturbative expansion around a mean-field density $\mu$, where each term in the perturbation series involves only a few variables (up to two variables in \eqref{eq:two cluster ground truth}). Although it is straightforward to generalize \eqref{eq:two cluster ground truth} to an expansion that includes higher-order interactions involving more than two variables with decaying coefficients, for clarity of presentation, we use \eqref{eq:two cluster ground truth} as our default model. Such an expansion is common in perturbation and cluster expansion theory in statistical mechanics.

Based on this setting, we discuss several estimators with increasing levels of sophistication. The first approach to estimate $p^*$ is based on a simple and natural subspace projection:  
\begin{equation}\label{eq:projection based density estimator}
\tilde{p} = \mu \odot \left(\Phi \text{diag}(w_K) \Phi^T \hat{p}\right).
\end{equation}  
For notational simplicity, we treat all functions and operators above as discrete matrices and vectors: we regard the concatenation of basis functions
\[
\Phi := [\Phi_l(x)]_{x,l}
\]  
as a matrix with rows indexed by $x$ and columns indexed by $l$. The operation $\Phi^T \hat{p}$ corresponds to convolution between basis functions and $\hat{p}$, i.e., summing or integrating out the variable $x$. The vector $w_K$, indexed by $l$, contains the weights associated with the basis functions and depends on a chosen positive integer $K$. $\text{diag}(\cdot)$ is the operator producing diagonal matrix with the elements of the input vector on the main diagonal. The mean-field $\mu$ is also treated as a vector, indexed by $x$, and $\odot$ denotes the Hadamard (entrywise) product. Intuitively, the operator $\mu \odot ( \Phi \text{diag}(w_K) \Phi^T \cdot )$ resembles a projector that projects $\hat{p}$ onto a function space that only involves a small number of variables determined by the chosen integer chosen integer $K$.  More precisely, let the complete (up to a certain numerical precision) set of $k$-variable basis functions be defined as: 
\begin{equation}\label{k-cluster basis}
\mathcal{B}_{k,d} = \left\{ \phi^{j_1}_{l_1}(x_{j_1}) \cdots \phi^{j_k}_{l_k}(x_{j_k}) \ |\    1 \leqslant j_1 < \cdots < j_k \leqslant d \quad  \text{and} \quad  2\leqslant l_1,\cdots , l_k \leqslant n   \right\}.
\end{equation}  
which we refer to as the \emph{$k$-cluster basis} throughout the paper. In particular, $\mathcal{B}_{0,d} = \{1\}$ only contains the constant basis. Furthermore, we say a function is \emph{$k$-cluster} if it can be expanded by up to $k$-cluster basis. For example, the perturbed part of $p^*$ in \eqref{eq:two cluster ground truth} is 2-cluster. Based on these definitions, $\Phi$ can be formulated as the concatenation of the $k$-cluster basis across different $k$:
\begin{equation}\label{cluster-basis-matrix}
\Phi(x) := \left[\mathcal{B}_{k,d}(x)\right]_{k=0}^{d}.
\end{equation}  
Furthermore, we define the weight vector $w_K$ as:  
\begin{equation}\label{eq:wK definition}
w_K(j) := \begin{cases} 
1 & j \in \binom{d}{K}, \\ 
0 & \text{otherwise},
\end{cases}
\end{equation}  
which is an indicator vector that selects the $0$- to $K$-cluster basis.  With these definitions, one can verify that if $p^*$ takes the form of \eqref{eq:two cluster ground truth}, and if the univariate basis functions $\{\phi^j_l\}_{l=1}^n$ satisfy the orthogonality condition \eqref{eq:orthogonal poly def}, then we can recover
\begin{equation}\label{2-cluster_recover}
p^* =  \mu \odot \left( \Phi \text{diag}(w_K) \Phi^T p^*\right)
\end{equation}  
as long as $K \geqslant 2$. This fact motivates the proposed density estimator in \eqref{eq:projection based density estimator}.  Moreover, unlike $\hat{p}$, whose variance scales as $\bigO(n^d)$, \Cref{lemma:k-cluster var} shows that the variance of $\tilde{p}$ in \eqref{eq:projection based density estimator} is significantly reduced, scaling as polynomially dependent term $\bigO(d^K)$ due to the projection onto a lower-dimensional function space. This substantial reduction in variance is a key advantage of our method. 

While the estimator in \eqref{eq:projection based density estimator} is intuitive, it is possible to further reduce the variance by imposing additional assumptions on $p^*$. Specifically, we assume that the low-order cluster coefficients $\text{diag}(w_K)\Phi^T p^*$ can be represented as a low-rank tensor train. The details of such a representation are provided in Section \ref{sec:notation}. Under this assumption, our estimation procedure can be modified as:  
\begin{equation}\label{eq:projection based density estimator with SVD}
\tilde{p} =  \mu \odot \left(\Phi \left(\texttt{TT-compress}\left(\text{diag}(w_K) \Phi^T \hat{p}\right)\right)\right),
\end{equation}  
where we apply certain tensor-train compression algorithm, denoted as \texttt{TT-compress}, to the tensor $\text{diag}(w_K) \Phi^T \hat{p}$. The singular value thresholding involved in the compression procedure ($\texttt{TT-compress}$) can further suppress noise in $\text{diag}(w_K) \Phi^T \hat{p}$, leading to a more accurate density estimator.

While the estimator in \eqref{eq:projection based density estimator with SVD} is expected to achieve a lower variance, a major drawback is that the projection $\text{diag}(w_K)\Phi^T \hat{p}$ onto the cluster basis up to order $K$ has a computational complexity of $\bigO(d^K n^K)$, corresponding to the number of cluster basis functions. This scaling quickly becomes prohibitive as $d$ and $K$ increase. To address this computational bottleneck, we propose a final estimator, which can be viewed as a more efficient variant of \eqref{eq:projection based density estimator with SVD}:  
\begin{equation}\label{eq:projection based density estimator with soft-threshold}
\tilde{p} = \mu \odot\left(\Phi  \text{diag}(\tilde{w}_{\alpha})^{-1} \left(\texttt{TT-compress}\left(\text{diag}(\tilde{w}_{\alpha}) \Phi^T \hat{p}\right)\right)\right),
\end{equation}  
where $\tilde{w}_{\alpha} \in \mathbb{R}^{n^d}$ represents a soft down-weighting of the cluster basis. This approach differs from using $w_K$, which applies a hard threshold that entirely discards cluster bases involving more than $K$ variables. 
A key advantage of using such $\tilde{w}_{\alpha}$ is that it can be designed as a rank-1 tensor, allowing the sequence of operations in \eqref{eq:projection based density estimator with soft-threshold} to be executed with a computational complexity of only $\bigO(dn)$. This drastically reduces the cost compared to the $\bigO(d^K n^K)$ complexity in \eqref{eq:projection based density estimator with SVD}. Furthermore, the final density estimator $\tilde{p}$ naturally takes the form of a tensor train, allowing essential tasks such as moment computation, sample generation, and pointwise evaluation of the density to scale linearly in $d$. In Section~\ref{sec:mf-example}, we present a simple example demonstrating that the estimator in \eqref{eq:projection based density estimator with soft-threshold} yields accurate results comparable to those of \eqref{eq:projection based density estimator with SVD}, thereby illustrating their similar compression capabilities. The study of \eqref{eq:projection based density estimator with soft-threshold} will be the focus of this paper.

\subsection{Related work}
In this subsection, we survey several lines of works that are most closely related to the proposed method.

\quad \textbf{Density estimation literature.} 
Classical density estimation methods such as histograms (\cite{scott1979optimal,scott1985averaged}), kernel density estimators (\cite{parzen1962estimation,davis2011remarks})
and nearest neighbor density estimates (\cite{mack1979multivariate}) are known for their statistical stability and numerical robustness in low-dimensional cases. But they often struggle with the curse of dimensionality, even in moderately high-dimensional spaces. More recently, deep generative modeling has emerged as a powerful tool for generating new samples, particularly in the context of image, based on neural network architecture. 
This approach includes energy-based models (\cite{kim2016deep,du2019implicit,song2019generative})
generative adversarial networks (GANs) (\cite{goodfellow2020generative,liu2021density}), normalizing flows (\cite{dinh2014nice,rezende2015variational,dinh2016density}), and autoregressive models (\cite{germain2015made,uria2016neural,papamakarios2017masked}). While the numerical performance for these modern methods is impressive, statistical guarantees and non-convex optimization convergence analysis remain challenging areas of research. Furthermore, for some of these models (e.g. energy based model), generating i.i.d. samples can still be expensive. For methods like GANs, while it can generate samples efficiently, it does not provide a way to evaluate the value of a density at any given point. In short, it is difficult for a method to enjoy all the desirable properties such as easy to train, easy to generate samples, can support evaluation of density values, and can be estimated without the curse of dimensionality. 
\par 

\textbf{Tensor train from partial sensing of a function.} Tensor train, also known as the matrix product state (MPS) (\cite{white1993density,perez2006matrix}) in the physics literature, has recently been used when to uncover a function, when only limited number of observations of this function are given. Methods for this include TT-cross (\cite{oseledets2010tt}), DMRG-cross (\cite{savostyanov2011fast}), TT completion (\cite{liu2012tensor,song2019tensor}), and others (\cite{sun2020low,kressner2023streaming,shi2023parallel}). It has broad applications in high-dimensional scientific computing problems (\cite{kressner2016low,bachmayr2016tensor} and machine learning (\cite{stoudenmire2016supervised,huggins2019towards,sengupta2022tensor}). 
\par
\textbf{Tensor train for density estimation}. In a density estimation setting, we are not given observed entries of a density function, but just an empirical distribution of it.  The key challenge lies in the randomness of the observed object, i.e. the empirical histogram, which contains an exponentially large variance in terms of the dimensionality. Therefore, reducing variance in a high-dimensional setting is crucial to design effective algorithms. Recently, work (\cite{hur2023generative}) and subsequent studies (\cite{tang2022generative,peng2023generative,chen2023combining,tang2023solving,khoo2024nonparametric}) employ sketching techniques in tensor decomposition step for density estimation problems, which originates from randomized linear algebra literature. 
However, most existing algorithms in these works do not have a general way to ensure a linear computational complexity when sketching, in terms of dimensionality. Another approach for density estimation with tensors is through optimization-based method. Works (\cite{bradley2020modeling,han2018unsupervised,novikov2021tensor})  implement gradient descent methods and optimize tensor-network representations according to some loss functions.  However, the performance strongly depends on the initialization and often suffers from the problem of local minima due to the highly non-convex nature of the optimization landscape. Furthermore, these iterative procedures require multiple iterations over $N$ data points, resulting in higher computational complexity compared to the previous tensor decomposition algorithms. \par

\subsection{Organization}
Here we summarize the organization of the paper. In Section \ref{sec:notation}, we introduce the tensor notation and diagram convetions used throughout. In Section \ref{sec:main approach}, we present the core ideas behind our proposed “convolution–deconvolution” framework,  highlighting its three key steps. The central component is the \texttt{TT-compress} algorithm as discussed previously that compresses cluster coefficients into a tensor train. Section \ref{sec:tt-svd} -- \ref{sec:tt-svd-t} introduce various of such  \texttt{TT-compress} algorithms: Section \ref{sec:tt-svd} presents a baseline algorithm based on sequential truncated singular value decomposition (SVD), along with its fast implementation for density estimation problems. Section \ref{sec:tt-svd-kn} includes an algorithm that utilizes the Nystr\"om approximation to further accelerate the algorithm in Section \ref{sec:main approach}. In Section \ref{sec:tt-svd-c}, we study an algorithm that performs SVD efficiently by truncating the elements in the tensor, and then introduce a further acceleration by employing a hierarchical matrix factorization. Section \ref{sec:tt-svd-t} proposes an algorithm utilizing random sketching matrices with low-rank tensor structures. Section~\ref{sec:pre-process} presents a pre-processing procedure based on principal component analysis (PCA), enabling our compression algorithms to handle very high-dimensional problems. Theoretical analysis for the algorithms is provided in Section~\ref{sec:theory}. In Section~\ref{sec:numerical}, we showcase the comprehensive numerical results to demonstrate the high accuracy and efficiency for our proposed algorithms. Finally, in Section \ref{sec:conclusion}, we summarize our findings.

\subsection{Tensor notations and diagrams}\label{sec:notation}
\begin{figure}[!htb]
    \centering
    \subfloat[Tensor diagrams of 3-dimensional $A$ and 2-dimensional $B$.]{{\includegraphics[width=0.32\textwidth]{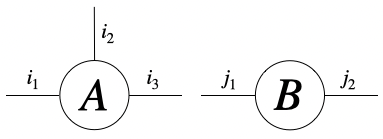}}}
    \qquad 
    \subfloat[Tensor diagram of tensor contraction. ]{{\includegraphics[width=0.48\textwidth]{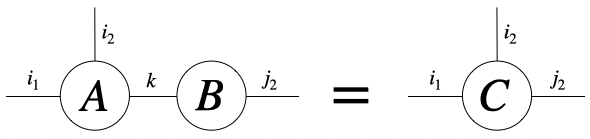}}}
        \caption{Basic tensor diagrams and operations. }
    \label{fig:tt-core-diagram}    
\end{figure}

We begin with some definitions regarding tensors. A $d$-dimensional tensor $A\in \mathbb{R}^{n_1\times n_2\cdots\times n_d}$ is a collection of entries denoted by $A(i_1,i_2,\cdots,i_d)$ with $1\leqslant i_j\leqslant n_j$ for $1\leqslant j\leqslant d$. Diagrammatically, a tensor can be denoted as a node whose dimensionality is indicated by the number of the ``legs" of the node. In Figure \ref{fig:tt-core-diagram}(a), the diagrams for a three-dimensional tensor $A$ and a two-dimensional matrix $B$ are plotted as an example. We use $i_1,i_2,i_3$ to denote the three legs in tensor $A$. We may also use a rectangle to denote a tensor when it has many legs in this paper. 
\par
\begin{figure}[ht!]
\centering
\includegraphics[width=0.4\textwidth]{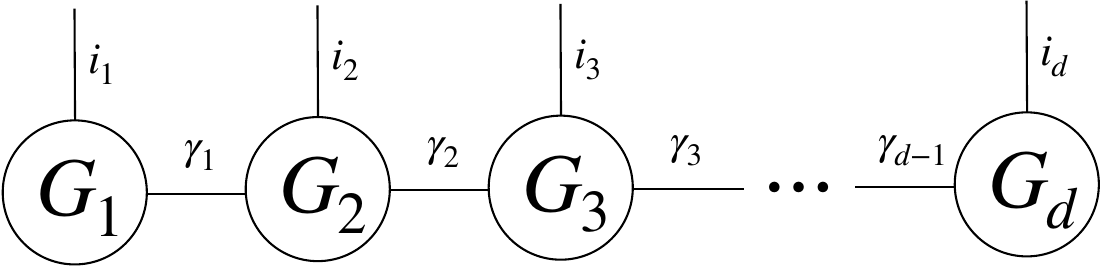} \qquad  \qquad 
\includegraphics[width=0.4\textwidth]{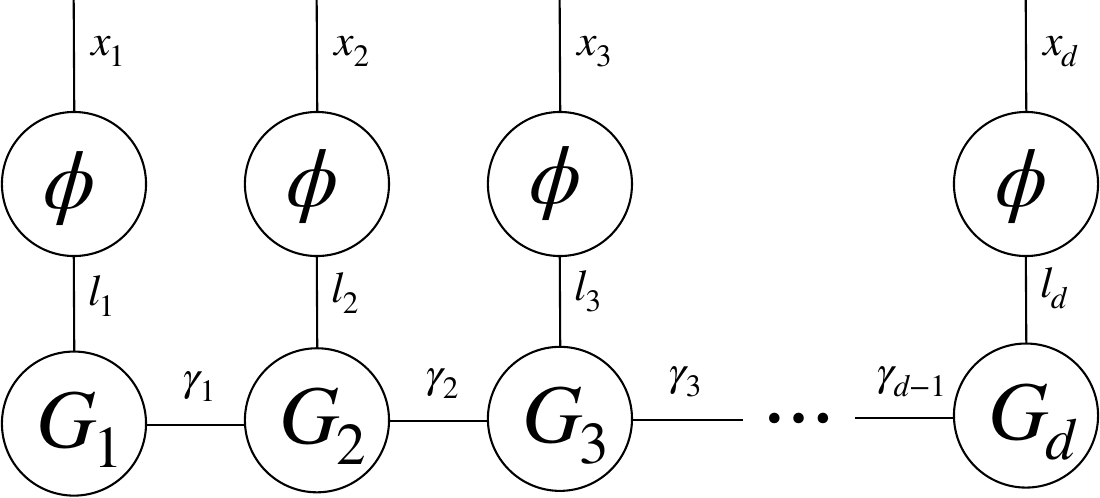}
\caption{Diagram for tensor-train representation. Left: discrete case. Right: continuous case.}
\label{fig:tt-diagram}
\end{figure}

 
Based on the definition of tensor, we list several tensor operations that will be used frequently in the paper:
\begin{enumerate}
    \item \textbf{``$\circ$'' contraction.} We use $\circ$ to represent the tensor cores contraction 
$(A\circ B)(i_1,i_2,j_2) = \sum_{k} A(i_1,i_2,k)B(k,j_2)$. Such operation is denoted as diagram in Figure \ref{fig:tt-core-diagram} (b) which combines the last leg of $A$ and the first leg of $B$. The resulting $C = A\circ B$ is a three-dimensional tensor. 
\item \textbf{Unfolding matrix.} 
Given a $d$-dimensional tensor with entries $A(i_1,i_2,\cdots,i_d)$, if we group the first $j$ indices into a multi-index $i_{1:j}$ 
and group the rest into another multi-index $i_{j+1:d}$, we will obtain a \emph{$j$-th unfolding matrix} $A^{(j)}\in \mathbb{R}^{(n_1\cdots n_j)\times( n_{j+1}\cdots n_d)}$ (also called $j$-th unfolding of $A$) satisfying 
\begin{equation}\label{unfolding matrix}
    A^{(j)}(i_{1:j},i_{j+1:d}) =  A(i_1,i_2,\cdots,i_d). 
\end{equation}
In particular, $A^{(0)} \in \mathbb{R}^{1\times n_1\cdots n_d}$ and $A^{(d)} \in \mathbb{R}^{n_1\cdots n_d \times 1}$ are respectively the row and column vectors reshaped from the tensor.

\end{enumerate}
In general, the memory cost of a high-dimensional tensor suffers from curse of dimensionality. Tensor-train representation is a way to decompose the whole exponential large $d$-dimensional tensor with underlying low-rank structure into $d$ number of components:
\begin{equation}\label{eq:TT representation} 
    A = G_1\circ G_2\circ \cdots \circ G_d \in \mathbb{R}^{n_1\times \cdots \times n_d}
\end{equation}
where $G_1\in \mathbb{R}^{n_1\times r_1}, G_2\in \mathbb{R}^{r_1\times n_2\times r_2},\cdots, G_{d-1}\in \mathbb{R}^{r_{d-2}\times n_{d-1}\times r_{d-1}},G_d\in \mathbb{R}^{r_{d-1}\times n_d}$ are called tensor cores and $\{r_j\}_{j=1}^{d-1}$ are the ranks. The TT representation can be visualized as the left panel of Figure \ref{fig:tt-diagram}. Throughout the paper, we denote the largest rank by $$r = \max_j r_j.$$ Then the memory cost for storing $A$ in TT form is at $\bigO(dnr^2)$ , which has no curse of dimensionality if $r$ is small. In addition, a TT can also be used to represent a continuous function $f(x_1,\cdots,x_d)$ formulated as a tensor product basis function expansion:
\begin{displaymath}
f(x_1,x_2,\cdots,x_d)=\sum_{l_1,l_2,\cdots,l_d=1}^n A(l_1,l_2,\cdots,l_d)\phi_{l_1}(x_1)\phi_{l_2}(x_2)\cdots \phi_{l_d}(x_d)
\end{displaymath}
where the coefficient $A$ takes the TT format as \eqref{eq:TT representation} . Such functional tensor-train representation is illustrated in the right panel of Figure \ref{fig:tt-diagram}. 

\section{Main approach}\label{sec:main approach}
In this section, we detail the framework of our main approach. In order to perform \eqref{eq:projection based density estimator with soft-threshold}, we perform the following three steps: 
\begin{enumerate}
    \item \textbf{Convolution.} Project the empirical distribution $\hat{p}$ to ``coefficients" $\hat{c}$ under certain basis expansion with $\bigO(d N)$ computational cost: 
    \begin{equation}\label{c_hat}
        \hat{c}= \text{diag}(\tilde{w}_{\alpha})\Phi^T \hat{p} = \text{diag}(\tilde{w}_{\alpha})\int \Phi(\vx)^T\hat{p}(\vx)d\vx .
    \end{equation}
    \item \textbf{Tensor compression.} Approximate $\hat c$ as a TT $\tilde c = \texttt{TT-compress}(\hat c)$. This step can done by several proposed algorithms listed in Table~\ref{tab:sketching-table}. Notably, several of these algorithms achieve $\bigO(dN)$ complexity. These complexity estimates will be further validated through numerical experiments in Section~\ref{sec:numerical}.
\begin{table}[H]
    \centering
    \begin{tabular}{|c|c|c|}
    \hline 
      Algorithm   & Corresponding section &Computational complexity \\
         \hline
      TT-SVD   & \S \ref{sec:tt-svd}  & $\bigO(\min(dn^{d+1},dN^2n))$ \\
            \hline 
       TT-SVD-kn   &  \S \ref{sec:tt-svd-kn}& $\bigO(dNn)$ \\
      \hline 
       TT-SVD-c ($K$-cluster)  &  \S \ref{sec:tt-svd-c}&$\bigO(d^{K+1}Nn^{K+1})$ \\
       \hline  
     TT-SVD-c (hierarchical)  & \S \ref{section: hie cluster}&$\bigO(d\log(d)Nn)$ \\
       \hline 
       TT-rSVD-t   & \S \ref{sec:tt-svd-t} & $\bigO(dNn)$ \\
       \hline 
    \end{tabular}
    \caption{Computational complexity of proposed \texttt{TT-compress}.}
    \label{tab:sketching-table}
\end{table}
    \item \textbf{Deconvolution.} Get the desired density estimator based on the TT coefficient $\tilde{c}$ with $\bigO(dN)$ computational cost:
   \begin{equation}
       \tilde{p} = \mu \odot (\Phi\text{diag}(\tilde{w}_{\alpha})^{-1}\tilde{c}) .
   \end{equation} 
   Here $\tilde{p}$ is also a TT, if $\Phi\text{diag}(\tilde{w}_{\alpha})^{-1}$ has an appropriate low-rank tensor structure.
\end{enumerate}
Overall, the computational cost our approach (with proper choices of \texttt{TT-compress}) can achieve computational complexity that scales linearly with both the dimensionality $d$ and sample size $N$, demonstrating the high efficiency. In the rest of this section, we will further elaborate convolution (Step 1) and deconvolution (Step 3). We defer the details of tensor compression (Step 2) to Section~\ref{sec:tt-svd} -- \ref{sec:tt-svd-t}.

\subsection{Convolution step}\label{sec: kernel choice}
\begin{figure}[ht!]
\centering
\includegraphics[width=0.8\textwidth]{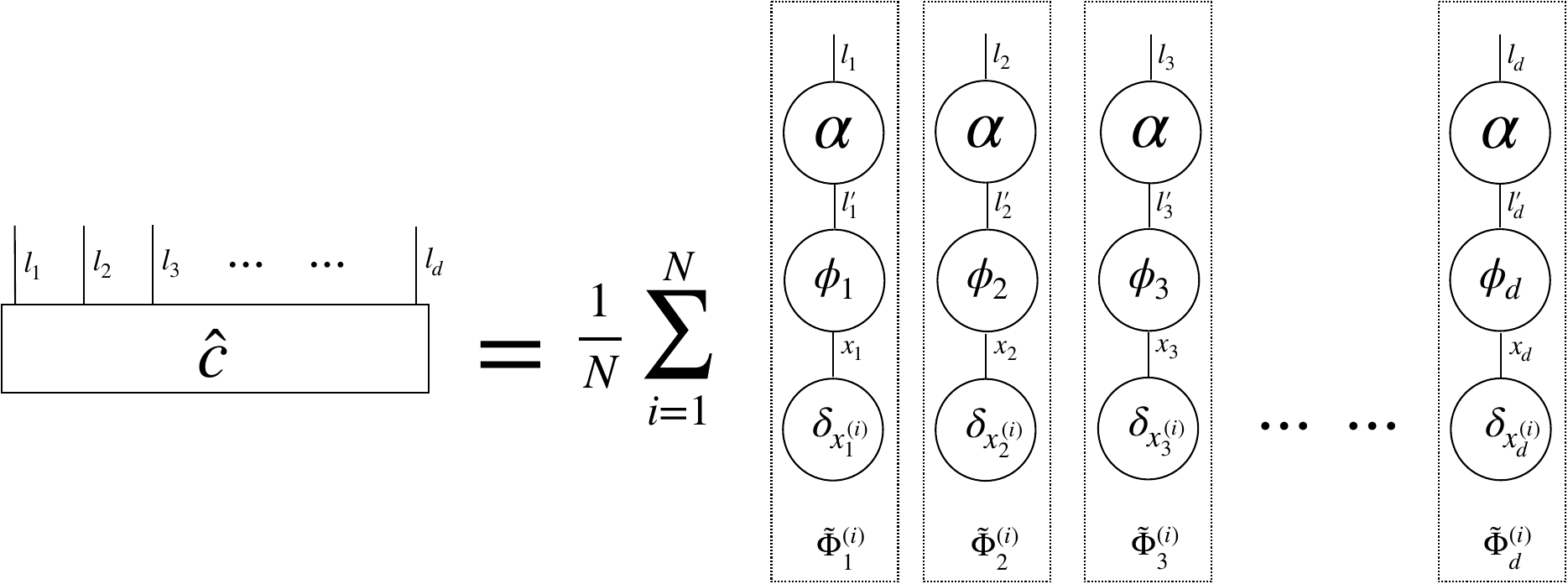} 
\caption{Tensor diagram of $\hat{c}= \text{diag}(\tilde{w}_{\alpha})\Phi^T \hat{p}$. ``$\alpha$" denotes diagonal matrix $\text{diag}(1,\alpha,\cdots,\alpha)\in \mathbb{R}^{n\times n}$; ``$\phi_j$" denotes ``matrix" $[\phi^j_l(x)]_{x,l}$. }
\label{fig:phi_left_right}
\end{figure}

The key part of this step is the choice of $\tilde{w}_{\alpha}\in \mathbb{R}^{n^d}$, the soft down-weighting of the cluster basis \eqref{cluster-basis-matrix}. In this work, we choose the weight of any $k$-cluster basis $\phi_{l_1}^{j_1}(x_{j_1})\cdots\phi_{l_k}^{j_k}(x_{j_k})$ to be 
\begin{equation}\label{eq: alpha choice}
    \alpha^{k}, \text{~where parameter~}  \alpha \in (0,1].
\end{equation}
In this way, the coefficient tensor $\hat{c}= \text{diag}(\tilde{w}_{\alpha})\Phi^T \hat{p}$ has entries
\begin{displaymath}
    \hat{c}(l_1,l_2,\cdots,l_d)  = \frac{1}{N}\sum_{i=1}^N {\alpha}^{\|l-\boldsymbol{1}\|_0} \prod_{j=1}^d \phi^j_{l_j}(x_j^{(i)}) .
\end{displaymath}
where $\|\cdot\|_0$ is the zero norm counting the number of non-zero entries of a vector, and $\boldsymbol{1}$ is the all-one $d$-dimensional vector. $\|l-\boldsymbol{1}\|_0$ is then the number of non-constant univariate basis functions among $\phi_{l_1}^{j_1}(x_{j_1})\cdots\phi_{l_k}^{j_k}(x_{j_k})$. Under such choice, $\text{diag}(\tilde{w}_\alpha)$ can be formulated as a rank-1 tensor operator, allowing the fast convolution (and also fast deconvolution later) with only $\bigO(dN)$ computational cost. In fact, given samples $\{x^{(i)}\}_{i=1}^N$, $\hat{c}$ is formulated as 
\begin{equation}\label{chat_density}
    \hat{c} =  \frac{1}{N}\sum_{i=1}^N 
    \tilde{\Phi}^{(i)}_{1} \otimes \tilde{\Phi}^{(i)}_{2} \otimes \cdots   \otimes \tilde{\Phi}^{(i)}_{d}
\end{equation}
where $\tilde{\Phi}^{(i)}_{j}=[1,\alpha \phi^j_{2}(x^{(i)}_j),\cdots, \alpha \phi^j_{n}(x^{(i)}_j)] \in \mathbb{R}^{ n}$. This expression is plotted as tensor diagrams in Figure \ref{fig:phi_left_right}. For future use, we introduce 
\begin{equation}\label{Phi_tilde_matrix}
    \tilde{\Phi}_{j} \in \mathbb{R}^{N\times n}
\end{equation}
 to denote the matrix whose rows are $\tilde{\Phi}^{(i)}_{j}$ .

The idea of such choice is to assign smaller weights for higher-order cluster basis to suppress the variance. The parameter $\alpha$ describes the denoising strength --- a smaller $\alpha$ produces stronger denoising and thus smaller variance. In particular when $\alpha = 1$, this soft-thresholding method \eqref{eq:projection based density estimator with soft-threshold} reduces to the hard-thresholding \eqref{eq:projection based density estimator with SVD} using all cluster bases. We point out that by choosing sufficiently small $\alpha$, this soft-thresholding strategy has smaller variance compared to the hard-thresholding approach \eqref{eq:projection based density estimator with SVD}. The following lemma considers a toy model to show the variance of the hard-thresholding approach:
\begin{lemma}
\label{lemma:k-cluster var}
Suppose $\{\phi_l\}_{l=1}^n$ are Legendre polynomials and the ground truth density $p^*$ satisfies \eqref{2-cluster_recover} with the mean-field component $\mu \sim \emph{Unif}([0,1]^d)$, the following variance estimation holds
\begin{equation}
\mathbb{E}\|\Phi \emph{diag}(w_K)\Phi^T\hat{p} - \Phi \emph{diag}(w_K)\Phi^Tp^* \|_{\lt}^2 = \bigO(\frac{d^Kn^{2K}}{N}) .
\end{equation}

\end{lemma}
\begin{proof}
Note that the error 
$\|\text{diag}(w_K)\Phi^T\hat{p} - \text{diag}(w_K)\Phi^Tp^* \|_F^2$ has at most $ {d \choose K}n^K $ terms, where each term follows
\begin{equation}
    \left(\int (\prod_{k=1}^K \phi^{j_k}_{l_k}(x_{j_k}))(\hat{p}(x) - p^*(x) ) dx \right)^2
\end{equation}
for indices $(j_1,\cdots,j_K)\in {d\choose K}$ and $l_1,\cdots,l_K \in \{2,\cdots,n\}$. Since Legendre polynomial satisfies $\|\phi^{j_k}_{l_k}\|_\infty\leqslant \sqrt{n}$, then
\begin{equation}
    \mathbb{E}\left(\int (\prod_{k=1}^K \phi^{j_k}_{l_k}(x_{j_k}))(\hat{p}(x) - p^*(x) ) dx \right)^2 = \frac{1}{N}\int (\prod_{k=1}^K \phi^{j_k}_{l_k}(x_{j_k}))^2 p^*(x) dx \leqslant \frac{n^K}{N}\int p^*(x) dx = \frac{n^K}{N}
\end{equation}
Therefore, the total variance is bounded by 
\begin{equation}
\mathbb{E}\|\Phi \text{diag}(w_K)\Phi^T\hat{p} - \Phi \text{diag}(w_K)\Phi^Tp^* \|_{\lt}^2 = \mathbb{E}\|\text{diag}(w_K)\Phi^T\hat{p} - \Phi^Tp^* \|_F^2 \leqslant 
 {d \choose K}n^K \frac{n^K}{N} = \bigO(\frac{d^Kn^{2K}}{N}).
\end{equation}
\end{proof}

On the other hand, the variance of the soft-thresholding approach under the model assumption in \Cref{lemma:k-cluster var} is estimated by
\begin{align}\label{eq: variance bound}
  & \mathbb{E}\left[ \| \Phi \text{diag}(\tilde{w}_\alpha)\Phi^T \hat{p}- \Phi \text{diag}(\tilde{w}_\alpha)\Phi^T p^*\|^2_{\lt} \right] \nonumber\\ 
  = \ &   \mathbb{E}\left[ \|\text{diag}(\tilde{w}_\alpha)\Phi^T \hat{p}-\text{diag}(\tilde{w}_\alpha)\Phi^T p^*\|^2_{F} \right]\nonumber\\
 = \ &  \sum_{l_1,\cdots,l_d=1}^{n} {\alpha}^{2\|l-\boldsymbol{1}\|_0} \mathbb{E}\left[ \int (\hat{p}(\vx)-p^*(\vx))\left( \prod_{j=1}^d \phi^j_{l_j}(x_j)\right) \dd \vx \right]^2 \nonumber\\
= \ & \frac{1}{N}\sum_{l_1,\cdots,l_d=1}^{n}  {\alpha}^{2\|l-\boldsymbol{1}\|_0} \int \prod_{j=1}^d  (\phi^j_{l_j}(x_j))^2 p^*(x) \dd x  \nonumber\\
\leqslant \ & \frac{1}{N}\sum_{l_1,\cdots,l_d=1}^{n}  {\alpha}^{2\|l-\boldsymbol{1}\|_0} n^{\|l-\boldsymbol{1}\|_0} \int p^*(x)\dd x \nonumber\\
= \ &  \frac{1}{N} \sum_{l_1,\cdots,l_d=1}^{n} (\alpha^2n)^{\|l-\boldsymbol{1}\|_0}
 = \frac{1}{N} (1 + (n-1)n{\alpha}^2)^d
\end{align}
where the inequality follows $\|\phi^j_{l_j}\|_\infty \leqslant \sqrt{n}$ for $l_j\geq 2$ and we use the binomial theorem for the last equality. With this estimation, we conclude that if we set $\alpha=1$, the variance scales at least $\bigO(n^d)$, indicating the curse of dimensionality. To avoid this, we choose 
\begin{displaymath}
    \alpha = \sqrt{\frac{C}{(n-1)nd}}.
\end{displaymath}
for some constant $C$ without $d$-dependency.
Under this choice, \eqref{eq: variance bound} can be bounded by 
\begin{equation}
   \mathbb{E}\left[ \|\Phi\text{diag}(\tilde{w}_\alpha)\Phi^T \hat{p}-\Phi\text{diag}(\tilde{w}_\alpha)\Phi^T p^*\|^2_{\lt} \right] \leqslant  \frac{1}{N}(1+\frac{C}{d})^d \leqslant \frac{1}{N}  e^C
\end{equation}
The resulting error has no curse of dimensionality. However, we also point out that $\alpha$ should not be chosen to be too small, which might lead to instability in our later operations. This will be further elaborated in Section \ref{sec: kernel inverse}.

\subsection{Deconvolution step}\label{sec: kernel inverse} 
\begin{figure}[ht!]
\centering
\includegraphics[width=0.4\textwidth]{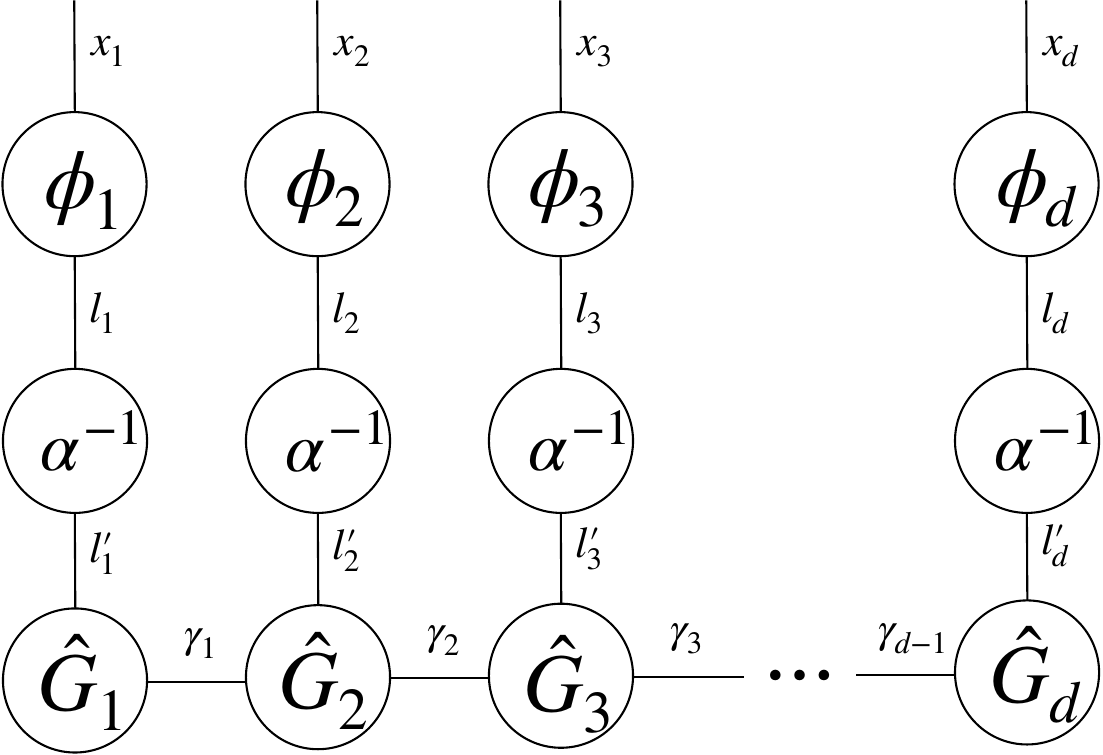}
\caption{Tensor diagram for the final density estimator $\tilde{p}$. ``$\alpha^{-1}$" denotes diagonal matrix $\text{diag}(1,1/\alpha,\cdots,1/\alpha)\in \mathbb{R}^{n\times n}$; ``$\phi_j$" denotes ``matrix" $[\phi^j_l(x)]_{x,l}$. }
\label{fig:solution_diagram}
\end{figure}

After we have received the approximated coefficient $\tilde{c}$ in TT format:
\begin{displaymath}
    \tilde{c} = \hat{G}_1 \circ   \hat{G}_2 \circ \cdots   \circ  \hat{G}_d \in \mathbb{R}^{n\times n\times \cdots \times n}. 
\end{displaymath}
in tensor compression step, the last step then obtains the final density estimator in TT format by implementing the deconvolution 
\begin{equation}\label{eq: kernel inverse}
    \tilde{p} = \mu \odot (\Phi\text{diag}(\tilde{w}_\alpha)^{-1}\tilde{c}) .
\end{equation}
The structure of $\Phi\text{diag}(\tilde{w}_\alpha)^{-1}\tilde{c}$ is illustrated in Figure \ref{fig:solution_diagram}. 
Since $\tilde{c}$ is a low-rank tensor train, the complexity of the evaluation on a single point $\vx$ is only $\bigO(dn)$.

\begin{remark}
    In practice, as dimensionality $d$ grows, some elements in $\text{diag}(\tilde{w}_\alpha)$ become exponentially small, which potentially lead to instability when computing the inverse. One solution is to implement inverse with given threshold: $\tilde{p}= \Phi(\emph{diag}(\alpha)+\lambda I)^{-1}\tilde{c}$. Here $\lambda$ is a small thresholding constant to avoid instability issues and $I$ is identity matrix. In practice, we find that even without this post-procedure, the tensor-train estimator does not suffer from such an instability.
    

\end{remark}

\subsection{Motivation: mean-field example}
\label{sec:mf-example}
In this subsection, we walk through the main steps of the proposed approach using a mean-field example of the ground truth density, i.e., $p^* = \mu$ in \eqref{eq:two cluster ground truth}. This toy example is used to illustrate why the soft-thresholding strategy in \eqref{eq:projection based density estimator with soft-threshold} effectively projects an empirical density onto the space spanned by low-order cluster bases, providing a similar effect as \eqref{eq:projection based density estimator with SVD}. The analysis can be extended to general examples. 

For convenience, we first introduce some notations for moments: we define the first-order moment $\hat{q}_j\in \mathbb{R}^{(n-1)\times 1}$ as 
\begin{equation}
    \hat{q}_j(l_j) = \int \phi^j_{l_j+1}(x_j)\hat{p}(x) \dd x, \qquad l_j=1,\cdots,n-1,
\end{equation}
and the second-order moment $\hat{q}_{j,j'}\in \mathbb{R}^{(n-1)\times (n-1)}$
\begin{equation}
    \hat{q}_{j,j'}(l_j,l'_{j'}) = \int \phi^j_{l_j+1}(x_j)\phi^{j'}_{l'_{j'}+1}(x_{j'})\hat{p}(x)\dd x, \qquad l_j,l'_{j'}=1,\cdots,n-1.
\end{equation}
Higher-order moments can be defined following the same logic.\par 
We consider to compress the tensor $\text{diag}(\tilde{w}_\alpha)\Phi^T\hat{p}$ into a tensor train using one specific approach called CUR approximation (\cite{mahoney2009cur}). Specifically, a given matrix $A \in \mathbb{R}^{m\times n}$ can be approximated by the following rank-$r$ factorization 
\begin{displaymath}
    A = C U^{+} R
\end{displaymath}
where $C \in \mathbb{R}^{m\times r}$ and $R \in \mathbb{R}^{r\times n}$ respectively consist of selected $r$ columns and rows. $U$ is the intersection (cross) submatrix formed by the selected rows and columns, and $U^+$ is the Moore-Penrose generalized inverse of $U$. To ensure a good approximation, the columns and rows are typically chosen to maximize the determinant of $U$. Repeating the procedure for all dimensions, we can receive a tensor-train decomposition eventually. Further details will be discussed in Section \ref{sec:tt-svd-kn}. 
Based on the notations above, we first unfold the weighted projection tensor $\text{diag}(\tilde{w}_\alpha)\Phi^T\hat{p}$ as 
\begin{equation}
     (\text{diag}(\tilde{w}_\alpha)\Phi^T\hat{p})^{(1)} = 
     \begin{bmatrix}
         1 & \alpha [ \hat{q}_j^{(0)} ]_{j=2}^d & \alpha^2 [ \hat{q}_{j,j'}^{(0)} ]_{j,j'=2}^d & \cdots \\
         \alpha \hat{q}_1 & \alpha^2 [ \hat{q}_{1,j}^{(1)} ]_{j=2}^d & \alpha^3 [ \hat{q}_{1,j,j'}^{(1)}]_{j,j'=2}^d & \cdots 
     \end{bmatrix} =: B_1
\end{equation}
in terms of the order of $\alpha$, where we use $[\cdot]$ to represent the concatenation of those inside terms through the column. Since the ground truth density $p^*$ is mean-field (rank-1), we choose only $r=1$ column and row for CUR approximation of $(\text{diag}(\tilde{w}_\alpha)\Phi^T\hat{p})^{(1)}$. With $\alpha$ chosen to be sufficiently small, the $(1,1)$-entry is the largest in $B_1$. Therefore, first row and first column are chosen for the CUR approximation: 
\begin{equation}
B_1
     \approx 
     \begin{bmatrix}
     1 \\ \alpha \hat{q}_1
     \end{bmatrix}
     \underbrace{\begin{bmatrix}
      1 &\alpha [ \hat{q}_j^{(0)}]_{j=2}^d & \alpha^2[  \hat{q}_{j,j'}^{(0)}]_{j,j'=2}^d & \cdots 
     \end{bmatrix}}_{=:B_2}  
\end{equation}
The left vector $     \begin{bmatrix}
     1 \\ \alpha \hat{q}_1
     \end{bmatrix}$ 
 is the resulting first TT core (before normalization). Next, we continue to implement the cross approximation on the unfolding of remaining tensor $B_2$. First row and first column are again chosen  following the same logic:
 \begin{equation}
     B_2 =      \begin{bmatrix}
         1 & \alpha[  \hat{q}_j^{(0)}]_{j=3}^d & \alpha^2 [ \hat{q}_{j,j'}^{(0)} ]_{j,j'=3}^d & \cdots \\
         \alpha \hat{q}_2 & \alpha^2[ \hat{q}_{2,j}^{(1)}]_{j=3}^d & \alpha^3 [ \hat{q}_{2,j,j'}^{(1)}]_{j,j'=3}^d & \cdots 
     \end{bmatrix} \approx 
      \begin{bmatrix}
     1 \\ \alpha \hat{q}_2
     \end{bmatrix} 
     \underbrace{\begin{bmatrix}
      1 & \alpha[ \hat{q}_j^{(0)}]_{j=3}^d & \alpha^2[ \hat{q}_{j,j'}^{(0)} ]_{j,j'=3}^d & \cdots 
     \end{bmatrix}}_{=:B_3}.
 \end{equation}

We repeat this procedure until decomposing the tensor into the product of $d$ rank-$1$ cores:
\begin{equation}\label{eq:mean-field CUR}
   \text{diag}(\tilde{w}_\alpha)\Phi^T\hat{p} \approx 
   \begin{bmatrix}
     1 \\ \alpha \hat{q}_1
     \end{bmatrix} \otimes  
        \begin{bmatrix}
     1 \\ \alpha \hat{q}_2
     \end{bmatrix} \otimes  \cdots \otimes     \begin{bmatrix}
     1 \\ \alpha \hat{q}_d
     \end{bmatrix} =: \texttt{TT-compress}\left(\text{diag}(\tilde{w}_{\alpha}) \Phi^T \hat{p}\right).
\end{equation}
Furthermore, we try to estimate the difference between the representation in \eqref{eq:mean-field CUR} and hard truncation 
$\texttt{TT-compress}\left(\text{diag}(w_K) \Phi^T \hat{p}\right)$ in $\eqref{eq:projection based density estimator with SVD}$, which are the summation of all $k$-cluster terms for $k>K$. Since ground truth density $p^*$ is mean-field, the expectation of empirical one-marginal is $\mathbb{E}[\hat{q}_j]=0$. Besides, variance $\text{Var}[\hat{q}_j]=\bigO(\frac{n-1}{N})$ for each $j$ due to the orthogonality of $\{\phi^j_l\}_{l=1}^n$ with respect to $\mu$. Therefore,
we can have an upper bound of $k$-cluster, 
$\alpha^k \|\hat{q}_{j_1}\otimes \cdots \otimes \hat{q}_{j_k}\|_F \leqslant C\alpha^k (\frac{n}{N})^{k/2}$    
for some constant $C>0$. Then compared to the hard-thresholding $K$-cluster estimator in $\eqref{eq:projection based density estimator with SVD}$, the difference towards the soft-thresholding estimator~\eqref{eq:projection based density estimator with soft-threshold} is upper bounded by 
\begin{align}
  & \|\left(\texttt{TT-compress}\left(\text{diag}(w_K) \Phi^T \hat{p}\right)\right)- \left(\texttt{TT-compress}\left(\text{diag}(\tilde{w}_{\alpha}) \Phi^T \hat{p}\right)\right)\|_F \leqslant  C\sum_{k=K+1}^d {d\choose k} \alpha^k (\frac{n}{N})^{k/2} \nonumber\\
  & \leqslant C \sum_{k=K+1}^d (\frac{\alpha de}{K}\sqrt{\frac{n}{N}})^k \leqslant Cd (\frac{\alpha de}{K}\sqrt{\frac{n}{N}})^{K+1}
\end{align}
when $\frac{\alpha de}{K}\sqrt{\frac{n}{N}}<1$ for sufficient large $N$. The difference between two estimators is polynomially small with respect to sample size $N$. This turns out to be a clear example to demonstrate the reliability of soft-thresholding approach~\eqref{eq:projection based density estimator with soft-threshold}. \par

\section{Tensor-train SVD (TT-SVD)}\label{sec:tt-svd}
Starting from this section, we focus on Step 2 of the main approach --- tensor compression. To this end, we first propose to use the tensor-train SVD (TT-SVD) algorithm. The TT-SVD algorithm is originally proposed in \cite{oseledets2011tensor}, which compresses a high-dimensional tensor into a low-rank tensor train by performing truncated SVD sequentially. For our problem, we use TT-SVD to find the approximation $\tilde{c}$ in TT format of $\hat{c}$ depicted in Figure \ref{fig:phi_left_right}, which can be obtained following the procedures in Algorithm \ref{algorithm: TT-SVD}.
\begin{algorithm}[H]
\begin{algorithmic}[1]
	\INPUT Tensor $\hat{c} \in \mathbb{R}^{n\times n\times \cdots \times n}$. Target rank $\{r_1,\cdots,r_{d-1}\}$.	 
 
     \For {$ j  \in \{ 1,  \cdots,   d - 1 \}$} 
  \State Perform SVD to $B_j \in \mathbb{R}^{nr_{j-1}\times n^{d-j}}$, which is formed by the following tensor contraction
\begin{center}
    \includegraphics[width=0.5\textwidth]{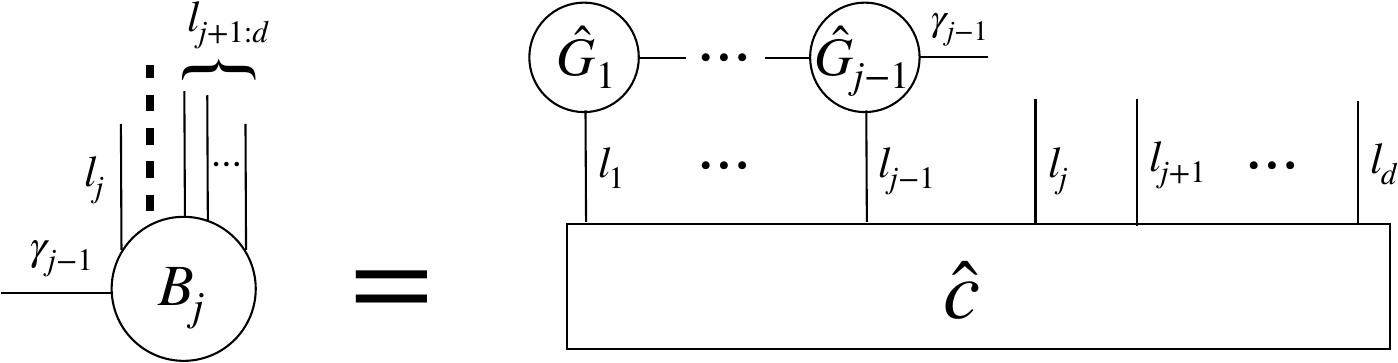} 
\end{center}
In particular, $B_1 = \hat{c}^{(1)}$ which is the first unfolding matrix of $\hat{c}$ defined in \eqref{unfolding matrix}. With SVD $B_j = U_j \Sigma_j V^T_j$, we obtain the $j$-th core $\hat{G}_j \in \mathbb{R}^{r_{j-1}\times n \times r_j}$, which is reshaped from the first $r_j$ principle left singular vectors in $U_j$. 
  \EndFor
    \State Obtain the last core $\hat{G}_d = B_d \in \mathbb{R}^{r_{d-1}\times n }$ . 
  \OUTPUT Tensor train $\tilde{c}=\hat{G}_1\circ \hat{G}_2 \circ \cdots \circ \hat{G}_d \in \mathbb{R}^{n\times n\times \cdots \times n}$. 
 \caption{TT-SVD} 
\label{algorithm: TT-SVD}
\end{algorithmic}
\end{algorithm}
Below, we briefly outline the idea behind this algorithm. When $j=1$, we want to obtain the first core $\hat{G}_1$ by performing SVD to $\hat{c}^{(1)}$ (which is $B_1$ in Algorithm~\ref{algorithm: TT-SVD}). Then we let $\hat{G}_1\in \mathbb{R}^{n\times r_1}$ to be first $r_1$ principle left singular vectors. The rationale of the first step is to approximate $\hat c^{(1)}$ as 
$\hat c^{(1)} \approx \hat G_1 (\hat G_1^T \hat c^{(1)} )$. In the $j=2$ step, we want to obtain $\hat{G}_2$, which gives a further approximation to $\hat G_1 (\hat G_1^T \hat c^{(1)} )$. The way we do this is to factorize $\hat G_1^T \hat c^{(1)} \in \mathbb{R}^{r_1\times n^{d-1}}$, by applying SVD to  $B_2$, which is the reshaping of $\hat G_1^T \hat c^{(1)}$. This completes the $j=2$ step. Then we do the above procedure recursively till $j=d-1$. After implementing sequential SVD for $d-1$ times, the remaining $B_d\in \mathbb{R}^{r_{d-1}\times n}$ is a matrix and we take it as the final core in the tensor-train representation since no more decomposition is required.

Naively, the computational cost from SVD in Algorithm \ref{algorithm: TT-SVD} scales exponentially large as $\bigO(dn^{d+1})$ in terms of dimensionality. However, for the density estimation problem, we can take advantage of the structure of $\hat{c}$ as depicted in Figure \ref{fig:phi_left_right} to achieve a fast implementation of Algorithm \ref{algorithm: TT-SVD}. We will use the next subsection to elaborate such fast implementation.

\subsection{Fast implementation of TT-SVD}\label{sec:TT-SVD fast algorithm}
In this subsection, we present a fast implementation of Algorithm \ref{algorithm: TT-SVD} with a complexity of $\bigO(dN^2nr)$, which scales linearly with $d$ and is therefore well-suited for high-dimensional cases. Notably, this efficient implementation introduces no additional approximations and produces the same output as the vanilla implementation in Algorithm \ref{algorithm: TT-SVD}.

\begin{figure}[!htb]
    \centering
    \includegraphics[width=0.8\textwidth]{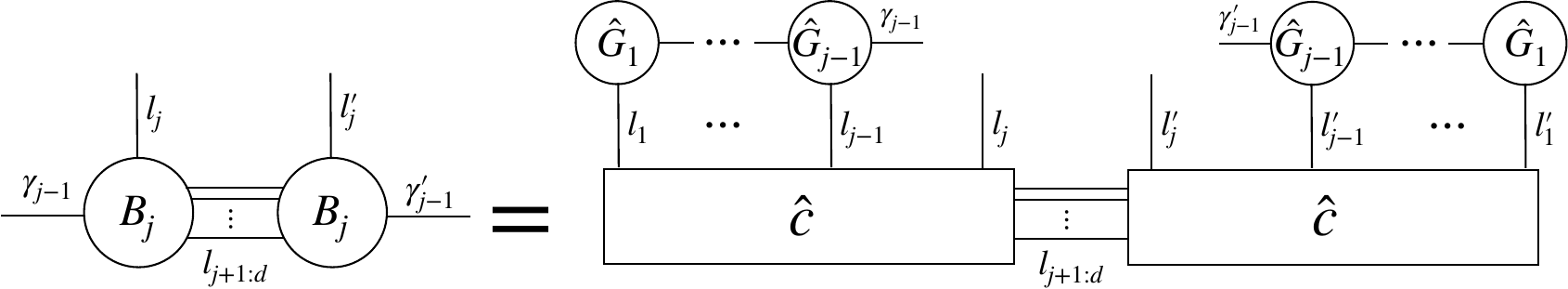}
     \caption{
    Tensor components for forming $B_jB^T_j$.}\label{fig:Bj} 
\end{figure}

The key of the fast implementation lies in fast calculations of SVD of $B_j$ in Algorithm \ref{algorithm: TT-SVD}. 
Instead of performing SVD to $B_j$ directly to get its left singular vectors as in Line 2 of Algorithm~\ref{algorithm: TT-SVD}, we consider computing the eigenvectors of $B_jB_j^T$, which is a relatively small matrix of size $nr_{j-1}\times nr_{j-1}$.  However, one concern is that the cost of forming $B_jB_j^T$ can in principle scale exponentially with dimensionality, since the column size of $B_j$ is $n^{d-j}$ as indicated by multiple legs between the two nodes in Figure \ref{fig:Bj}. Fortunately, by using the structure of $\hat{c}$ as shown in Figure \ref{fig:phi_left_right}, one can form $B_jB_j^T$ much more efficiently. This is shown in Line 2 of Algorithm~\ref{algorithm: TT-SVD-fast}, which summarizes the TT-SVD procedure that is implemented in a computationally feasible manner. 

\begin{figure}[ht!]
\centering
    \includegraphics[width=0.4\textwidth]{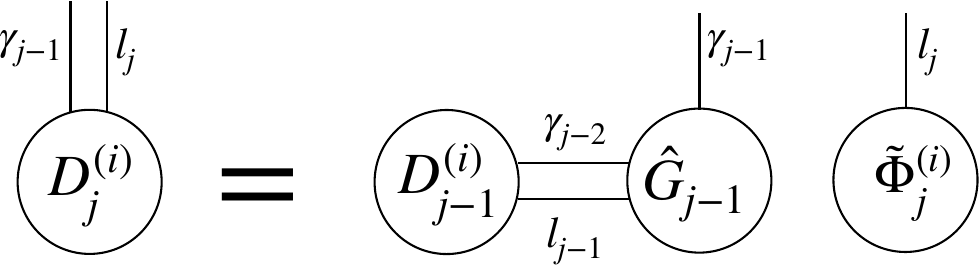} 
\caption{Tensor diagram illustrating the recursive calculation of $D_j^{(i)}$, which is the $i$-th column of the matrix $D_j$ in Algorithm \ref{algorithm: TT-SVD-fast}.}
\label{fig:Dj_iter}
\end{figure}

In Algorithm \ref{algorithm: TT-SVD-fast}, the main computation lies in forming $D_j$ and $E_j$ in Line 2. A naive way is to compute them independently, leading to $\bigO(dN^2)$ complexity. However, $D_j$ and $E_j$ can be computed in a recursive way efficiently, where computations can be reused. More specifically, one can precompute $E_j$ by
\begin{equation}\label{tt-svd-E-iter}
     E_j = E_{j+1} \odot  ( \tilde{\Phi}_{j} \tilde{\Phi}_{j}^T )  
 \end{equation}
 starting from $E_d = I_{N\times N}$ which is $N\times N$ the identity matrix. $D_j$ is computed recursively based on core $\hat{G}_{j-1}$ starting from $D_1 = \tilde{\Phi}^T_{1}$ according to Figure \ref{fig:Dj_iter}. The last core $\hat{G}_d$ is assigned as the matrix $B_d$ in algorithm~\ref{algorithm: TT-SVD}, which is equal to the average of $D_d^{(i)}$ over sample index $i$. Therefore, we compute $\hat{G}_d$ as the average of $D_d^{(i)}$ in Algorithm~\ref{algorithm: TT-SVD-fast}.
 In this way, the computational cost of $E_j$, $D_j$ and $B_jB_j^T$ is $\bigO(N^2nr)$ which has no dependency on dimensionality. Afterwards, one can solve the core by eigen-decomposition of $B_j B^T_j \in \mathbb{R}^{nr_{j-1}\times nr_{j-1}}$ with a minor cost as $\bigO((nr)^3)$. Overall, the computational cost of this fast implementation is $\bigO(dN^2nr)$.

While this fast implementation significantly reduces the computational complexity from exponential to linear in $d$, we point out that such implementation 
still requires $\bigO(N^2)$ computations to form each $E_j$, which can make it challenging with a large number of samples even in low dimensional cases. In Section \ref{sec:tt-svd-kn}, we will further improve Algorithm \ref{algorithm: TT-SVD-fast} using kernel Nystr\"om approximation, which will bring computational cost to $\bigO(dN)$ only.

\begin{algorithm}[H]
\begin{algorithmic}[1]
	\INPUT Samples $\{x^{(i)}\}_{i=1}^N$. Target rank $\{r_1,\cdots,r_{d-1}\}$.	
          \For {$ j  \in \{ 1,  \cdots,   d - 1 \}$} 
\State Perform eigen-decomposition to $A_j := B_jB_j^T \in \mathbb{R}^{nr_{j-1}\times nr_{j-1}}$, which is computed by $A_j = \frac{1}{N^2}D_jE_j D^T_j$ where $D_j \in \mathbb{R}^{nr_{j-1 \times N }}$ and $ E_j \in  \mathbb{R}^{N\times N}$ are the resulting matrices by contraction in the boxes of the below figure:
\begin{center}
    \includegraphics[width=0.7\textwidth]{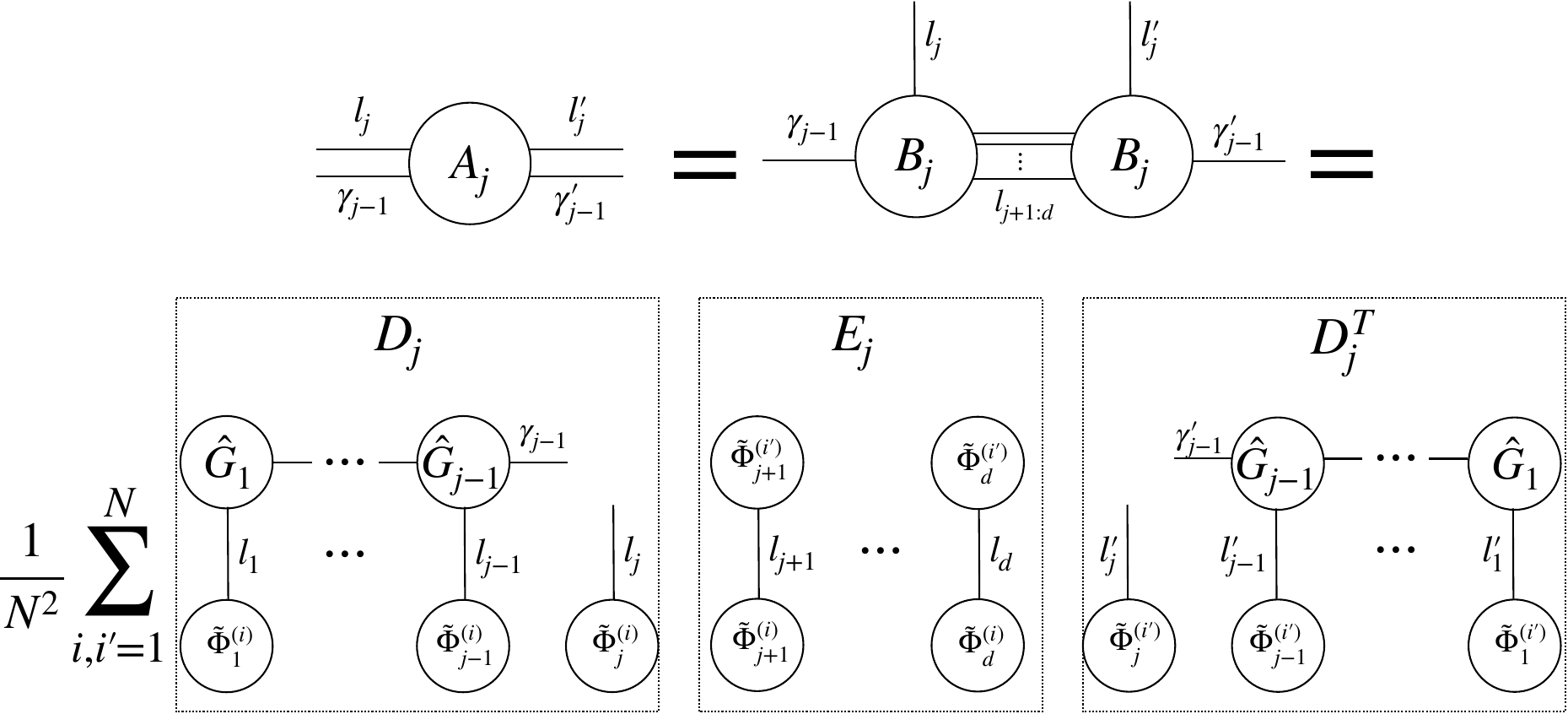} 
\end{center}
In particular, $D_1 = \tilde{\Phi}_1^T$ defined in \eqref{Phi_tilde_matrix}. Here $D_j$ and $E_j$ should be computed recursively to reuse calculations. With eigen-decomposition $A_j = U_j \Lambda_j U_j^T$, we obtain the $j$-th core $\hat{G}_j \in \mathbb{R}^{r_{j-1}\times n \times r_j}$, which is reshaped from the first $r_j$ principle eigenvectors in $U_j$.
  \EndFor
   \State Obtain the last core:
\begin{displaymath}
    \hat{G}_d = \frac{1}{N}\sum_{i=1}^N D^{(i)}_d \in \mathbb{R}^{r_{d-1}\times n }
\end{displaymath}
where $D^{(i)}_d$ is the reshaping of $i$-th column of $D_d$.
  \OUTPUT Tensor train $\tilde{c}=\hat{G}_1\circ \hat{G}_2 \circ \cdots \circ \hat{G}_d \in \mathbb{R}^{n\times n\times \cdots \times n}$. 
 \caption{Fast implementation of TT-SVD} 
\label{algorithm: TT-SVD-fast}
\end{algorithmic}
\end{algorithm}

\section{TT-SVD with kernel Nystr\"om approximation (TT-SVD-kn)}
 \label{sec:tt-svd-kn}
In this section, we propose a further modification based on the fast implementation of TT-SVD developed in Section \ref{sec:TT-SVD fast algorithm}. This modification accelerates Algorithm \ref{algorithm: TT-SVD-fast} using kernel Nystr\"om approximation to reduce the complexity from $\bigO(dN^2)$ to $\bigO(dN)$. As we have pointed out at the end of Section \ref{sec:TT-SVD fast algorithm}, the $\bigO(N^2)$ complexity of Algorithm \ref{algorithm: TT-SVD-fast} is the consequence of exact evaluation of the $N\times N$ matrix $E_j$ in the step of eigen-decomposition for $B_jB_j^T$. The proposed modification in this section will utilize kernel Nystr\"om approximation (\cite{williams2000using}) to form $E_j$ efficiently.

\begin{figure}[ht!]
\centering
\includegraphics[width=0.5\textwidth]{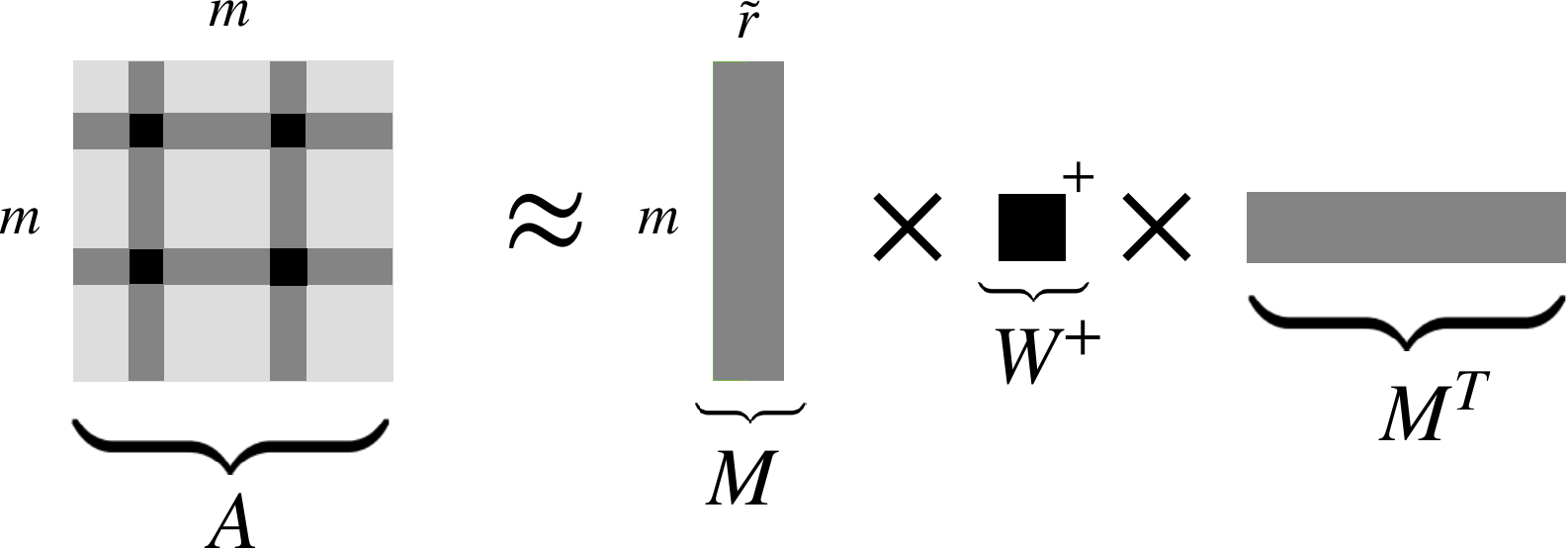}
\caption{Cross approximation of a symmetric matrix.}
\label{fig:cur}
\end{figure}

To begin with, we briefly introduce the idea of kernel Nystr\"om method to approximate a symmetric matrix as shown in Figure \ref{fig:cur}. Given a symmetric matrix $A \in \mathbb{R}^{m \times m}$, the kernel Nystr\"om method chooses $\tilde{r}$ columns of $A$ uniformly at random without replacement and approximates $A$ by the following low-rank cross approximation 
\begin{equation}\label{nystrom}
    A \approx M W^+ M^T
\end{equation}
where $M \in \mathbb{R}^{m\times \tilde{r}}$ consists of the $\tilde{r}$ chosen columns, $W \in \mathbb{R}^{\tilde{r}\times \tilde{r}}$ is the intersection of the chosen $\tilde{r}$ columns and the corresponding $\tilde{r}$ rows and $W^+$ is the Moore-Penrose generalized inverse of $W$.

\begin{figure}[ht!]
\centering
\includegraphics[width=1.0\textwidth]{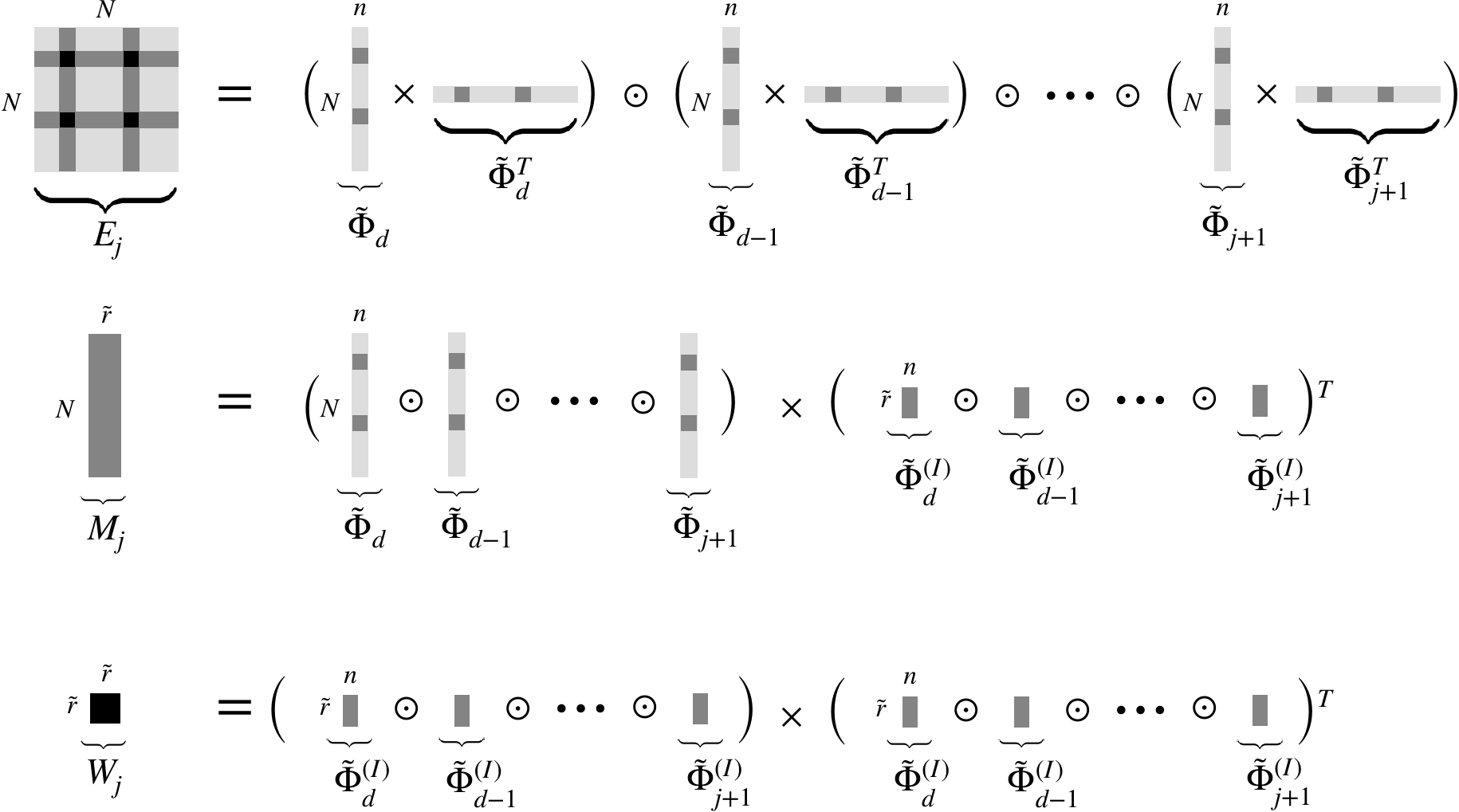}
\caption{Cross approximation of $E_j$, the most expensive component in Algorithm \ref{algorithm: TT-SVD-fast}.}
\label{fig:kn-tt}
\end{figure}

Such approximation can be applied to the symmetric matrix $E_j$ in Algorithm \ref{algorithm: TT-SVD-fast}, resulting in the more efficient TT-SVD-kn algorithm summarized in Algorithm \ref{sec:tt-svd-kn}. By the recurrence relation \eqref{tt-svd-E-iter}, $E_j$ is formulated as 
\begin{equation} \label{Ej_full}
    E_j = (\tilde{\Phi}_{d}\tilde{\Phi}_{d}^T) \odot  (\tilde{\Phi}_{d-1}\tilde{\Phi}_{d-1}^T) \odot \cdots \odot  (\tilde{\Phi}_{j+1}\tilde{\Phi}_{j+1}^T),
\end{equation}
which is depicted by the first equation in Figure \ref{fig:kn-tt}. To apply kernel Nystr\"om method, we randomly sample $\tilde{r}$ rows in $E_j$ with index $\mathcal{I} \subset \{1,2,\cdots,N\}$, and we define $\tilde{\Phi}^{(\mathcal{I})}_{m}\in\mathbb{R}^{\tilde{r}\times n}$, which is denoted by the blocks in each $\tilde{\Phi}_{m}$ in Figure \ref{fig:kn-tt}, as the resulting submatrices formed by the sampled rows. 

At this point, we perform the cross approximation for $E_j$:
\begin{displaymath}
    E_j \approx  M_j W_j^+ M_j^T
\end{displaymath}
where $M_j \in \mathbb{R}^{N\times \tilde{r}}$ is the submatrix formed by the chosen columns in $E_j$ and $W_j \in \mathbb{R}^{\tilde{r} \times \tilde{r}}$ is the corresponding intersection. Using the expression \eqref{Ej_full}, $M_j$ and $W_j$ are respectively formulated as the second and third equation in Figure \ref{fig:kn-tt}. The structures of $M_j$ and $W_j$ allow us to compute them in a recursive manner as in Line 2 of Algorithm \ref{algorithm: TT-SVD-kn}: suppose we have already obtained $\tilde{\Phi}_{j+2:d}$ and $\tilde{\Phi}^{(\mathcal{I})}_{j+2:d}$, $M_j$ and $W_j$ can be evaluated with complexity $\bigO(Nn)$ (direct computations based on their definitions in Figure \ref{fig:kn-tt} is $\bigO(dNn)$). Consequently, $B_jB_j^{T}$ in the previous Algorithm \ref{algorithm: TT-SVD-fast} can now be approximated by 
\begin{displaymath}
    B_jB_j^{T} =  \frac{1}{N^2} D_j E_j D_j^T \approx \frac{1}{N^2} D_j  M_j W_j^+ M_j^T D_j^T
\end{displaymath}
as in Line 4 of Algorithm \ref{sec:tt-svd-kn}. The computational cost for the above matrix multiplications is $\bigO(Nnr\max(nr,\tilde{r}))$, which is significantly smaller compared to the original $\bigO(N^2nr)$ complexity in Algorithm \ref{algorithm: TT-SVD-fast}. By considering the calculations for all $d$ tensor cores, the total computational cost is $\bigO(dNnr\max(nr,\tilde{r}))$. Here we remark that we have employed the same indices $\mathcal{I}$ for computing all $E_j$. Such practice is the key to reuse the calculations in this recursive implementation so that we can achieve a complexity that is only linear in sample size.

\begin{algorithm}[H]
\begin{algorithmic}[1]
	\INPUT Samples $\{x^{(i)}\}_{i=1}^N$. Target rank $\{r_1,\cdots,r_{d-1}\}$. Number of sampled columns $\tilde{r}$. 
 \State Sample $\tilde{r}$ integers from $\{1,\cdots,N\}$ and let $\mathcal{I}$ be the set of sampled integers.
 \State Pre-compute $\{W_j\}_{j=1}^{d-1}$ and $\{M_j\}_{j=1}^{d-1}$ based on sampled indices $\mathcal{I}$ recursively by Figure \ref{fig:kn-tt}.  
  \For {$ j  \in \{ 1,  \cdots,   d-1  \}$} 
\State Perform eigen-decomposition to $B_jB_j^T \in \mathbb{R}^{nr_{j-1}\times nr_{j-1}}$, which is approximated by
\begin{displaymath}
       B_jB_j^T \approx A_j := \frac{1}{N^2} D_j  M_j W_j^+ M_j^T D_j^T 
\end{displaymath}
where $D_j$ is computed recursively as Figure \ref{fig:Dj_iter}. With eigen-decomposition $A_j = U_j \Lambda_j U_j^T$, we obtain the $j$-th core $\hat{G}_j \in \mathbb{R}^{r_{j-1}\times n \times r_j}$, which is reshaped from the first $r_j$ principle eigenvectors in $U_j$.
  \EndFor
  \State Obtain the last core:
\begin{displaymath}
    \hat{G}_d = \frac{1}{N}\sum_{i=1}^N D^{(i)}_d \in \mathbb{R}^{r_{d-1}\times n }
\end{displaymath}
where $D^{(i)}_d$ is the $i$-th column of $D_d$.
  \OUTPUT Tensor train $\tilde{c}=\hat{G}_1\circ \hat{G}_2 \circ \cdots \circ \hat{G}_d \in \mathbb{R}^{n\times n\times \cdots \times n}$.
 \caption{TT-SVD-kn} 
\label{algorithm: TT-SVD-kn}
\end{algorithmic}
\end{algorithm}

\section{TT-SVD with truncated cluster basis (TT-SVD-c)}\label{sec:tt-svd-c}
In this section, we introduce the TT-SVD-c algorithm (Algorithm \ref{algorithm: cluster basis sketching}), which is a modified TT-SVD algorithm based on truncated cluster basis sketching. The central idea is to utilize sketching technique to approximate the full SVDs of $B_j\in \mathbb{R}^{nr_{j-1}\times n^{d-j}}$ in Line 2 of Algorithm~\ref{algorithm: TT-SVD}. 
More specifically, each step when we compute the left singular vectors $U_j$ from $B_j$ in Algorithm~\ref{algorithm: TT-SVD}, we instead compute the left singular vectors from a sketched $B_j$: 
\begin{equation}\label{Bj_sketching}
    B'_j=B_jS_j,
\end{equation}
where $S_j$ is some appropriately chosen sketching matrix with size $n^{d-j}\times \tilde{r}_j$, with sketching size $\tilde{r}_j$ much smaller than $n^{d-j}$ (Line 2 of Algorithm \ref{algorithm: cluster basis sketching}). In this way, we can reduce the exponential complexity of full SVD for $B_j$ to $\bigO(\max(nr\tilde{r}_j^2,n^2r^2\tilde{r}_j))$ complexity of SVD for $B_j'$.

\begin{algorithm}[H]
\begin{algorithmic}[1]
	\INPUT Samples $\{x^{(i)}\}_{i=1}^N$. Target rank $\{r_1,\cdots,r_{d-1}\}$. Cluster order $K$.
  \For {$ j  \in \{ 1,  \cdots,   d-1\}$} 

  \State Perform SVD to $B'_j  \in \mathbb{R}^{nr_{j-1}\times \tilde{r}_j}$, which is formed by the following tensor contraction 
  \begin{center}
    \includegraphics[width=0.5\textwidth]{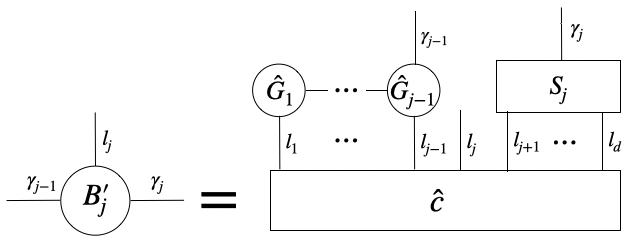} 
   \end{center} 
where $S_j \in \mathbb{R}^{n^{d-j}\times \tilde{r}_j}$ is the sketching matrix depending on cluster order $K$. The contraction of the righ-hand side is done efficiently as shown in Figure~\ref{fig: SVD-c-diagram}. With SVD $B'_j = U_j \Sigma_j V_j^T$, we obtain the $j$-th core $\hat{G}_j \in \mathbb{R}^{r_{j-1}\times n \times r_j}$, which is reshaped from the first $r_j$ principle left singular vectors in $U_j$ . 
  \EndFor
    \State Obtain the last core $\hat{G}_d = B'_d \in \mathbb{R}^{r_{d-1}\times n }$.
    \OUTPUT Tensor train $\tilde{c} =\hat{G}_1\circ \hat{G}_2 \circ \cdots \circ \hat{G}_d \in \mathbb{R}^{n\times n\times \cdots \times n}$. 
 \caption{TT-SVD-c} 
\label{algorithm: cluster basis sketching}
\end{algorithmic}
\end{algorithm}

The key of this approach is to choose good sketching matrix $S_j$ such that $B'_j$ retains the same range, and thus the same left singular space, as 
$B_j$. Here, we follow a specific way to construct sketching matrix $S_j$ (\cite{chen2023combining,peng2023generative}), where we only include the low-order cluster basis. We first introduce the set of \emph{$k$-cluster indices}:  
\begin{displaymath}
    \mathcal{J}_{k,d} =\bigl\{(l_1,\cdots,l_d) \ |\ l_1,\cdots,l_d \in \{1,\cdots,n\} \quad \text{and} \quad \|l-\boldsymbol{1}\|_0 = k \bigr\},
\end{displaymath}
which are the basis indices of those cluster bases in $\mathcal{B}_{k,d}$ defined in \eqref{k-cluster basis} since $\phi_1^j \equiv 1$.

Our choice of sketching matrix is based on the low-order cluster indices. To sketch $B_j$, we construct $S_j \in \mathbb{R}^{n^{d-j}\times \tilde{r}_j}$ as follows: each column of $S_j$ is defined based on a multivariate index $l \in \bigcup_{k=0}^{K} \mathcal{J}_{k,d-j}$ where $K$ is the chosen cluster order:
\begin{displaymath}
     \bvec{1}_{l_1} \otimes \bvec{1}_{l_2} \otimes \cdots \otimes \bvec{1}_{l_{d-j}}
\end{displaymath}
where $\bvec{1}_{l}$ is a $n\times 1$ column vector evaluating as 1 on the $l$th component and $0$ elsewhere. Under this construction, $S_j$ is a sparse matrix with column size $$\tilde{r}_j=\sum_{k=0}^{K} |\mathcal{J}_{k,d-j}|= \bigO(d^Kn^K).$$

\begin{figure}[!htb]
    \centering
    \subfloat[]{{\includegraphics[width=0.4\textwidth]{fig_D_iter.pdf} }} \qquad \quad 
    \subfloat[ ]{{\includegraphics[width=0.4\textwidth]{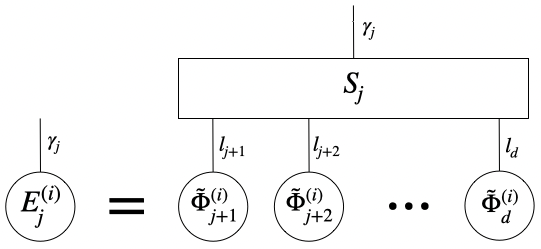} }}  \par 
    \subfloat[ ]{{\includegraphics[width=0.6\textwidth]{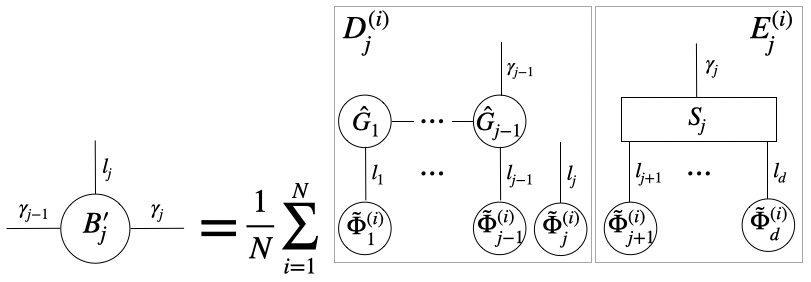} }} 
    \caption{Tensor diagrams illustrating recursive computation of $B'_j$ in Line 2 of Algorithm \ref{algorithm: cluster basis sketching}. (a) recursive calculation of $D_j^{(i)}$, which is the $i$-th column of the matrix $D_j$; (b)  recursive calculation of $E_j^{(i)}$, which is the $i$-th column of the matrix $E_j$; (c) calculation of $B_j'$ based on $D_j^{(i)}$ and $E_j^{(i)}$. }\label{fig: SVD-c-diagram}
\end{figure}

At this point, the only concern is that forming $B'_j = B_jS_j$ has exponential complexity by direct matrix multiplication. Similar to Algorithm \ref{algorithm: TT-SVD-fast}, we utilize the structure of $\hat{c}$ \eqref{chat_density}:
\begin{equation}\label{tt-svd-c-form-B'}
    B'_j=B_jS_j =\frac{1}{N}D_j (E_{j})^T = \frac{1}{N}\sum_{i=1}^N  D_j^{(i)}(E_{j}^{(i)})^T \in \mathbb{R}^{nr_{j-1}\times \tilde{r}_j}
\end{equation}
where $D_j,E_j$ are defined in Figure~\ref{fig: SVD-c-diagram} (c) and
$D^{(i)}_j,E^{(i)}_j$ represent the $i$-th column of matrix $D_j,E_j$, respectively. According to the structure of $D_j,E_j$, their computation can be done in a recursive way following the diagram in Figure~\ref{fig: SVD-c-diagram} (a) and (b).\par

In terms of the computational complexity of Algorithm~\ref{algorithm: cluster basis sketching}, the cost for $\{E_j\}_{j=1}^{d-1}$ is $\bigO(d^{K+1}Nn^K)$ for $k$-cluster truncation. Besides, the cost for $\{D_j\}_{j=1}^{d-1}$ follows the analysis in Figure~\ref{fig:Dj_iter} and Algorithm~\ref{algorithm: TT-SVD-fast}, which is $\bigO(dNnr^2)$. Furthermore, the computational complexity of multiplication in $B_j'$ and the following SVD step requires $\bigO(d^{K+1}Nn^{K+1})$. To summarize, the dominant term in computational complexity is $\bigO(d^{K+1}Nn^{K+1})$. The reduction is from standard TT-SVD algorithm $\bigO(\min(dn^{d+1},dN^2n))$ to $\bigO(d^{K+1}Nn^{K+1})$, which is significant especially for a large sample size $N$ and low-order $K$. In the following subsection, we further propose a modification of TT-SVD-c which constructs sketching matrix $S_j$ based on hierarchical matrix decomposition. This acceleration can achieve a lower computational complexity as $\bigO(d\log(d))$.

\subsection{Acceleration of TT-SVD-c via hierarchical matrix decomposition}\label{section: hie cluster}
In this subsection, we propose a modification that further accelerates TT-SVD-c method with cluster order $K=1$. The key insight of this acceleration lies in the assumption that the ground truth density $p^*$ has high correlations between nearby variables (such as $x_1$ and $x_2$) and low correlation between distant variables (such as $x_1$ and $x_d$), which is satisfied by many real cases. In Algorithm~\ref{algorithm: cluster basis sketching}, the sketching matrix $S_j$ considers all correlations equally without taking the distance into account, which includes redundant correlations with distant variables and leads to a large column size of $S_j$. In fact, even if cluster order $K=1$, the column size is $\bigO(d)$, which results in the overall computational cost of $\bigO(d^{2})$. This motivates us to further sketch $S_j$, which is done by hierarchical matrix decomposition (\cite{hackbusch1999sparse,hackbusch2000sparse,bebendorf2008hierarchical}) of covariance matrix $M$, whose $(j_1,j_2)$-block is formulated as 
\begin{equation}
    M_{j_1,j_2} = \big[\mathbb{E}_{\hat{p}}[\phi^{j_1}_{l_1}(x_{j_1})\phi^{j_2}_{l_2}(x_{j_2})] \big]_{l_1,l_2} \in \mathbb{R}^{n\times n}
\end{equation}
for all $1\leqslant j_1,j_2 \leqslant d$. 

\begin{figure}[ht!]
\centering
\includegraphics[width=0.4\textwidth]{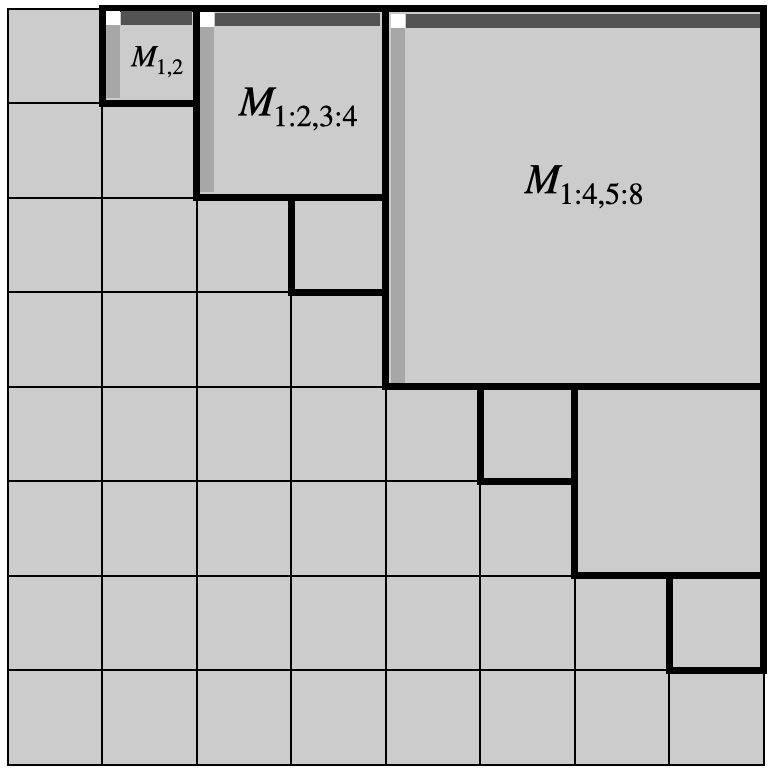}
\caption{Hierarchical matrix decomposition of covariance matrix for $d=8$. The figure depicts the three block matrices $M_{1:4,5:8},M_{1:2,3:4},M_{1,2}$ whose top $\tilde{r}$ right singular vectors (visualized as dark gray matrix in each block matrix).}
\label{fig:hie_visual}
\end{figure}

To illustrate the idea, we consider an example with $d=8$. Covariance matrix $M$ is depicted in Figure~\ref{fig:hie_visual}, where we use $M_{j_1:j_2,j_3:j_4}$ to represent the submatrix formed by blocks with row index from $j_1$ to $j_2$ and column index from $j_3$ to $j_4$. 
To implement hierarchical matrix decomposition, we perform SVD on upper triangular off-diagonal blocks $M_{1:4,5:8},M_{1:2,3:4},M_{1,2}$ and get top $\tilde{r}$ right singular vectors $V'_{1},V'_{2},V'_{3}$, respectively. Here the sketching size $\tilde{r}$ (different from that in Algorithm~\ref{algorithm: cluster basis sketching}) is chosen as a constant. $V'_{1},V'_{2},V'_{3}$ respectively capture correlations with variables $x_{5:8}$, $x_{3:4}$ and $x_2$.

Once the singular vectors $V'_{1},V'_{2},V'_{3}$ are obtained, we can use them to construct sketching matrices with smaller column sizes. When computing the first core $\hat{G}_1$, we form $Q_1=\text{diag}[V'_{1},V'_{2},V'_{3}]\in \mathbb{R}^{7n\times 3\tilde{r}}$ which captures the weighted correlation between $x_1$ and $x_{2:8}$, and we use $S_1$ and $Q_1$ to sketch $B_1$ as $B_1'=B_1S_1Q_1\in \mathbb{R}^{n\times 3\tilde{r}}$. Once this is done, we perform SVD to $B'_1$ to obtain $\hat{G}_1$. For the second core, we form $Q_2=\text{diag}[V'_{1},V'_{2}]\in \mathbb{R}^{6n\times 2\tilde{r}}$ to capture the weighted correlation between variables $x_2$ and variables $x_{3:8}$. Then we sketch $B_2'=B_2S_2Q_2$, and perform SVD to $B_2'$ to obtain $\hat{G}_2$. Note that we can reuse $V'_{1},V'_{2}$ in the second step since both of them were already computed in the first step. By repeating the procedure and generating new sketch $B_j'=B_jS_jQ_j$ for the remaining $j$, this hierarchical approach enables the efficient computation of all cores $\{\hat{G}_j\}_{j=1}^{d}$. \par

In terms of computational complexity, there are two major sources of computations. The first is the SVD of off-diagonal block matrices ($M_{1:4,5:8},M_{1:2,3:4},M_{1,2}$ for $d=8$). For each block, we first apply kernel Nystr\"om approximation \eqref{nystrom} with $\tilde{r}$ sampled rows and columns, followed by SVD. The cost for this part of computations is $\bigO(d\log(d)Nn\tilde{r})$. The second major computation is forming $B_j'=B_jS_jQ_j$. The number of columns in the weight matrix $Q_j$ is at most $\log_2(d)\tilde{r}$. Since $S_j$ is a sparse matrix, the computation of $B_j'=B_jS_jQ_j$ for all $j$ needs $\bigO(d\log(d)Nn\tilde{r})$ computational cost upon leveraging a similar strategy for computing \eqref{tt-svd-c-form-B'}. The computational costs discussed above are the dominant cost in Algorithm~\ref{algorithm: cluster basis sketching}. We can replace $S_j$ with $S_jQ_j$ in Algorithm~\ref{algorithm: cluster basis sketching} and follow the operations in Figure~\ref{fig: SVD-c-diagram}, the method achieves a reduced computational complexity $\bigO(d\log(d)Nn)$ compared to the original complexity $\bigO(d^{2}Nn^{2})$ for TT-SVD-c with cluster order $K=1$. \par

\section{TT-SVD with random tensorized sketching (TT-rSVD-t)}\label{sec:tt-svd-t}
In this section, we propose another modification for TT-SVD algorithm via sketching. This approach follows the randomized sketching technique (\cite{huber2017randomized}), originated from randomized linear algebra literature (\cite{halko2011finding}), to accelerate SVD.
The key idea is to choose sketching matrix $S_j \in \mathbb{R}^{n^{d-j}\times \tilde{r}}$ as some random matrix, where $\tilde{r}$ is a given sketching size, which is usually chosen to be reasonably larger than target rank $r$ of the output tensor-train estimator. Afterwards, similar as TT-SVD-c, we perform SVD to $B'_j=B_jS_j$. In general, the matrix multiplication $B_jS_j$ has exponentially large computational cost. Our solution to this curse of dimensionality is to form $S_j$ by a set of random tensor trains, leading to the TT-rSVD-t algorithm summarized in Algorithm \ref{algorithm: efficient tensor train sketching}. \par

\begin{algorithm}[H]
\begin{algorithmic}[1]
	\INPUT Samples $\{x^{(i)}\}_{i=1}^N$. Target rank $\{r_1,\cdots,r_{d-1}\}$. Sketching size $\tilde{r}$. 
 \State Generate random tensor train with cores $\{H_j\}_{j=1}^d$, where $H_{1}\in \mathbb{R}^{n\times \tilde{r}}, H_{2},\cdots ,H_{d-1} \in \mathbb{R}^{\tilde{r}\times n\times \tilde{r}} , H_{d} \in \mathbb{R}^{\tilde{r} \times n} $.
  \For {$ j  \in \{ 1,  \cdots,   d-1  \}$}  

\State Perform SVD to $B'_j \in \mathbb{R}^{nr_{j-1}\times \tilde{r}_j}$, which is formed by the following tensor contraction
    \begin{center}
    \includegraphics[width=0.6\textwidth]{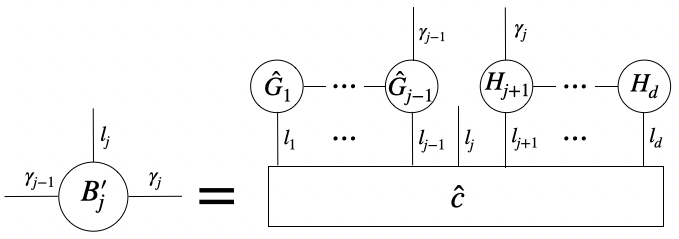} 
    \end{center}
The contraction of the righ-hand side is done efficiently as shown in Figure~\ref{fig:SVD-t-diagram}. With SVD $B'_j = U_j \Sigma_j V_j^T$, we obtain the $j$-th core $\hat{G}_j \in \mathbb{R}^{r_{j-1}\times n \times r_j}$, which is reshaped from the first $r_j$ principle left singular vectors in $U_j$ .   
  \EndFor
    \State Obtain the last core 
   $ \hat{G}_d = B'_d \in \mathbb{R}^{r_{d-1}\times n }$. 
    \OUTPUT Tensor train $\tilde{c} =\hat{G}_1\circ \hat{G}_2 \circ \cdots \circ \hat{G}_d \in \mathbb{R}^{n\times n\times \cdots \times n}$. 
 \caption{TT-rSVD-t} 
\label{algorithm: efficient tensor train sketching}
\end{algorithmic}
\end{algorithm}

To implement Algorithm \ref{algorithm: efficient tensor train sketching}, we first generate a $d$-dimensional random tensor train with fixed rank $\tilde{r}$ whose cores are $H_{1}, H_{2},\cdots, H_{d}$ (Line 1). In practice, we can set each core as i.i.d. random Gaussian matrix/tensor or random uniform matrix/tensor. When sketching $B_j$, we construct sketching matrix $S_j\in \mathbb{R}^{n^{d-j}\times \tilde{r}}$ as the reshaping of $H_{j+1}\circ H_{j+2}\circ \cdots \circ H_{d}$ (Line 3). Afterwards, we contract $B_j$ with $S_j$ and perform SVD to solve the $j$-th core $\hat{G}_j$.

\begin{figure}[!htb]
    \centering
    \subfloat[]{{\includegraphics[width=0.4\textwidth]{fig_D_iter.pdf} }} \qquad \quad 
    \subfloat[ ]{{\includegraphics[width=0.4\textwidth]{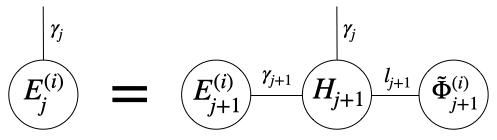} }}  \par 
    \subfloat[ ]{{\includegraphics[width=0.6\textwidth]{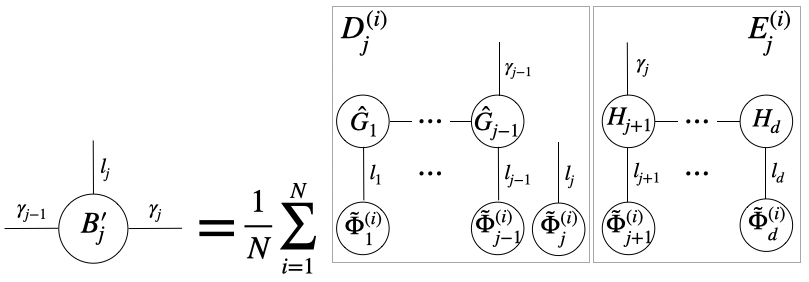} }} 
    \caption{Tensor diagrams illustrating recursive computation of $B'_j$ in Line 2 of Algorithm \ref{algorithm: efficient tensor train sketching}. (a) recursive calculation of $D_j^{(i)}$, which is the $i$-th column of the matrix $D_j$; (b)  recursive calculation of $E_j^{(i)}$, which is the $i$-th column of the matrix $E_j$; (c) calculation of $B_j'$ based on $D_j^{(i)}$ and $E_j^{(i)}$.  }\label{fig:SVD-t-diagram}
\end{figure}

Similar to \eqref{tt-svd-c-form-B'} for TT-SVD-c, \( B'_j \) is computed recursively as illustrated in Figure~\ref{fig:SVD-t-diagram}(c). The key difference lies in the new formulation of \( E_j \), which arises from the new choice of sketching matrix \( S_j \). Leveraging the tensor-train structure, \( E_j \) can also be computed recursively as shown in Figure~\ref{fig:SVD-t-diagram}(b). The total computational cost for \(\{E_j\}_{j=1}^{d-1}\) is \( \bigO(dNn\tilde{r}^2) \). The computation of \(\{D_j\}_{j=1}^{d-1}\), depicted in Figure~\ref{fig:SVD-t-diagram}(a), follows the same approach used in Algorithm~\ref{algorithm: TT-SVD-fast} and incurs a computational cost of \( \bigO(dNnr^2) \). Overall, the SVD step required to compute all tensor cores has a total complexity of \( \bigO(dNn\tilde{r}^2) \), which is linear with respect to both dimensionality and sample size.

\section{Pre-processing steps for general distributions}
\label{sec:pre-process}
In this section, we propose several pre-processing steps to extend the discussed framework to general distributions. The pre-processing procedure consists of a dimension-reduction step using principal component analysis (PCA), followed by a step to determine the mean-field and the univariate bases. The full algorithm is summarized in Algorithm~\ref{algorithm: general distribution alg}.  

\begin{algorithm}[H]
	\begin{algorithmic}[1]
		\INPUT Samples $\{x^{(i)}\}_{i=1}^N$. Embedding dimension $d'$, $d' \leqslant d$. Target rank $\{r_1,\cdots,r_{d-1}\}$.  
		\State Apply PCA to the samples $\{x^{(i)}\}_{i=1}^N$. Denote the top $d'$ orthonormal principal components as $Q \in \mathbb{R}^{d\times d'}$. Define the new variables as $z^{(i)} = Q^T x^{(i)}, \, i=1,\cdots,N$, with the associated empirical distribution denoted by $\hat{p}_z$.  
		\State Estimate the 1-marginals $\mu'_1(z_1)\cdots \mu'_{d'}(z_{d'})$ using one-dimensional kernel density estimation. The product of these 1-marginals forms the mean-field $\mu' = \prod_{j=1}^{d'} \mu'_j$ in the reduced coordinate system.  
		\State Compute the orthogonal basis $\{ \phi^j_{l}(z_j)\}_{l=1}^{n}$ for each $j = 1,\cdots ,d'$, where the basis functions are orthogonalized with respect to $\mu'_j(z_j)$.  
		\State Obtain the $d'$-dimensional tensor density estimator in the reduced coordinate system:
		\[
		\tilde{p}_z = \mu'\odot \text{diag}(\tilde{w}_\alpha)^{-1}\Phi\left(\texttt{TT-compress}\left(\text{diag}(\tilde{w}_\alpha) \Phi^T \hat{p}_z\right)\right),
		\]
		where $\texttt{TT-compress}$ is any tensor-train compression algorithm detailed in Section~\ref{sec:tt-svd} -- \ref{sec:tt-svd-t}. The final density estimator in the original space is obtained as:  
	\begin{displaymath}
		\tilde{p}(x) = \frac{1}{Z}\tilde{p}_z(Q^T x),
	\end{displaymath}
       where $Z$ is the normalization factor.
		\caption{Tensor density estimator after variable re-ordering.}  
		\label{algorithm: general distribution alg}  
    \OUTPUT Tensor density estimator $\tilde{p}$.  
	\end{algorithmic}
\end{algorithm}
When dealing with a general distribution, the optimal ordering of variables is generally unknown. A poor variable ordering may lead to a high-rank TT representation, significantly increasing computational complexity. To address this, Step 1 of Algorithm~\ref{algorithm: general distribution alg} applies PCA to simultaneously decorrelate the variables and reduce the dimensionality. The resulting reduced-coordinate variables $\{z^{(i)}\}_{i=1}^{N}$ often exhibit simpler dependence structures, allowing a more compact TT representation. The eigen-decomposition in PCA can be performed through a ``streaming approach'', as in \cite{henriksen2019adaoja,hoskins2024subspace,khoo2024nonparametric} with a complexity of $\bigO(Ndd')$.

In Step 2 of Algorithm~\ref{algorithm: general distribution alg}, we aim to perform the projection \eqref{eq:projection based density estimator with soft-threshold} on the empirical distribution $\hat{p}_z$. However, to do so, we require knowledge of the mean-field $\mu'(z) = \prod_{j=1}^{d'} \mu'_j(z_j)$. To approximate this, we apply a kernel density estimator to each of the one-dimensional marginals of $\hat{p}_z$. This step ensures that $\mu'$ captures the product structure of the new coordinate system.  

In Step 3, we construct orthonormal basis functions with respect to each one-dimensional marginal. This can be achieved using a Gram-Schmidt procedure or a stable three-term recurrence method (see \cite{gautschi2004orthogonal}) for polynomial basis functions.

Finally, in Step 4, we estimate the density $\hat{p}_z$ in the reduced-coordinate system using the tensor-train density estimator. This involves applying any tensor-train compression algorithms \texttt{TT-compress} detailed in Section~\ref{sec:tt-svd} -- \ref{sec:tt-svd-t} to $\text{diag}(\tilde{w}_\alpha)\Phi^T \hat{p}_z$. The resulting density estimator $\tilde{p}_z$ is then mapped back to the original space via the transformation $\tilde{p}(x) = \tilde{p}_z(Q^T x)$ followed by normalization.

\section{Theoretical results} 
\label{sec:theory}
In this section, we present an error analysis of our proposed methods. For simplicity, we focus primarily on the hard-thresholding approach described in \eqref{eq:projection based density estimator with SVD}. Based on Section~\ref{sec:mf-example}, we know the estimator from hard-thresholding approach in \eqref{eq:projection based density estimator with SVD} is close to the primarily discussed soft-thresholding approach in \eqref{eq:projection based density estimator with soft-threshold}. The detailed analysis of the soft-thresholding approach in \eqref{eq:projection based density estimator with soft-threshold} will be deferred to future work.

For simplicity, we assume the mean-field component $\mu \sim \text{Unif}([0,1]^d)$. In this way, the univariate basis functions $\{\phi_l\}_{l=1}^n$ are orthonormal, and we take Legendre polynomials as an example in our analysis. In addition, we propose the following assumptions on the ground truth density: 
\begin{assumption} \label{assumption:additive model}
    Suppose the ground truth density function $p^*:[0,1]^d \rightarrow \mathbb{R}^+$ satisfies is an additive model (1-cluster):
\begin{equation} \label{eq:additive model}
   p^* (x_1,\cdots, x_d) = 1 + q_1^*(x_1) + \cdots + q_d^*(x_d),
\end{equation}   
where $q_j^* \in H^{\alphaspace}([0,1])$ for some integer $\alphaspace \geqslant 1$   satisfies the regularity
\begin{displaymath}
    \min_{j} \| q^*_j\|_{\lt([0,1])} \geqslant b^* \quad \text{and} \quad \max_{j}\| q_j^*\|_{H^{\alphaspace}([0,1])}\leqslant B^*
\end{displaymath}
for two $\bigO(1)$ positive constants $b ^*$ and $B ^*$ where we set $B^* \geqslant 1$. Under this assumption, we can immediately obtained that $\|p^*\|_\infty \leqslant 1+ B^*d$ according to Sobolev embedding theorem.
\end{assumption}
Assumption~\ref{assumption:additive model} considers only the 1-cluster model and include the regularity constraint for each one-dimensional component. This motivation stems from the fact that the left/right subspaces are computationally straightforward, allowing for an efficient tensor-train representation of $p^*$. This specific structure, in turn, yields a favorable estimation error rate with respect to the dimensionality.

We state our main theorem under these assumptions:
\begin{theorem}\label{thm: main theorem}
Let $c^*=\text{diag}(w_K)\Phi^Tp^*$ and $ \hat {c} =\text{diag}(w_K)\Phi^T {\hat p}$ be the $\mathbb{R}^{ n ^d}$ coefficient tensor associated with the ground truth and empirical density function, respectively. Suppose the ground truth density function $p^*$ satisfies Assumption \ref{assumption:additive model}, 
and the sample size is sufficiently large such that 
\begin{equation}\label{eq:snr}
N\geqslant C_0d^{2K-1}n^{4K-2}\log(dn),
\end{equation}
where $C_0$ is some sufficiently large constant independent of $N$, $n$ and $d$. When we use TT-SVD (Algorithm~\ref{algorithm: TT-SVD}) as $\texttt{TT-compress}$ to compress $\hat{c}$ to $\tilde{c}$ with target TT rank $r=2$, the following bound holds with high probability: 
\begin{equation}
    \frac{\|\tilde{c} - c^*\|_F}{\|c^*\|_F} = \bigO\left( \sqrt{ \frac{dn B^* \log(dn) }{N} \sum_{j=1}^{d-1}\frac{1}{\sigma_2^2(c^{*(j)})} } \right)
\end{equation}
where $c^{*(j)}$ is the $j$th-unfolding matrix of $c^{*}$ as defined in \eqref{unfolding matrix}, and $\sigma_2(\cdot)$ denotes the second largest singular value. 
\end{theorem}
The proof of this theorem is lengthy and therefore deferred to \Cref{app:main proof}. The theorem states that the relative error of the coefficient estimates decreases at the Monte Carlo rate of $\bigO(1/\sqrt{N})$ with respect to the sample size $N$. Furthermore, under the assumption that $p^*$ is a 1-cluster, the relative error grows only linearly with the dimensionality $d$ (up to a logarithmic factor), provided that the second-largest singular values (also the smallest since $r = 2$) of all reshaped ground truth tensors $c^*$ are uniformly bounded below.  As a direct consequence of this theorem, one can derive the final error estimate for the density estimator:
\begin{corollary}
\label{thm: main corollary}
Suppose all the conditions in \Cref{thm: main theorem} hold, the following bound for density estimator $\tilde{p}$ \eqref{eq: kernel inverse} by TT-SVD (\Cref{algorithm: TT-SVD}) with target rank $r=2$ holds with high probability:
\begin{equation}
     \frac{ \|\tilde{p} - p^*\|_{\lt([0,1]^d) }  }{\|p^*\|_{\lt([0,1]^d) } } = \bigO  \bigg( \frac{B^*}{b^*} \Big(    n^{-\alphaspace} + \frac{B^*}{b^{*2}}\sqrt{\frac{dn  \log^2(dn)  }{N }    } \Big) \bigg),
\end{equation}
where $\alphaspace$ is the regularity parameter specified in \Cref{assumption:additive model}.   
\end{corollary}
The proof of \Cref{thm: main corollary} can be found in \Cref{app:main proof-2}.

\section{Numerical results}\label{sec:numerical}
In this section, we carry out numerical experiments to examine the performance of our proposed density estimation algorithms. In specific, we use TT-SVD-kn (Algorithm \ref{algorithm: TT-SVD-kn}), TT-SVD-c (Algorithm \ref{algorithm: cluster basis sketching}) with hierarchical matrix decomposition and TT-rSVD-t (Algorithm \ref{algorithm: efficient tensor train sketching}). Unless otherwise noted, we use the following parameter settings throughout: we set the number of sampled columns in Nystr\"om approximation $\tilde{r} = 100$ in TT-SVD-kn, sketching size $\tilde{r}=10$ in TT-SVD-c with hierarchical matrix decomposition, sketching size $\tilde{r} = 30$ in TT-rSVD-t and the kernel parameter $\alpha = 0.01$. These parameter choices are based on empirical observations.

For simplicity, we choose the mean-field density $\mu$ as a uniform distribution, so that $\{\phi^j_l\}_{l=1}^{n}$ can be chosen as any orthonormal basis functions. Here we use the Fourier basis 
\begin{equation} \label{fourier_basis}
    \phi_{l}(x) = \begin{cases}
        \frac{1}{\sqrt{2L}}, & \text{~for~} l = 1, \\ 
        \frac{1}{\sqrt{L}}\cos(\frac{l\pi x}{2L})  & \text{~for even~} l > 1 \\ 
        \frac{1}{\sqrt{L}}\sin(\frac{(l+1)\pi x}{2L}) & \text{~for odd~} l > 1 
    \end{cases}
\end{equation}
as the univariate basis function for all $j=1,\cdots,d$. 

\subsection{Gaussian mixture model}
We first consider the $d$-dimensional Gaussian mixture model which is defined as :
\begin{displaymath}
      \frac{1}{6}\mathcal{M}(0.18) +  \frac{1}{3}\mathcal{M}(0.2) + \frac{1}{2}\mathcal{M}(0.22)
\end{displaymath}
where 
\begin{displaymath}
   \mathcal{M}(\sigma) =  \frac{2}{3}\mathcal{N}(-0.5,\sigma^2 I_{d\times d}) + \frac{1}{3}\mathcal{N}(0.5,\sigma^2 I_{d\times d}).
\end{displaymath}
We set the domain as the hypercube $[-L,L]^d$ where $L=1.5$ with mesh size $0.1$ in each direction.  

\begin{figure}[ht!]
\centering
\includegraphics[width=0.49\textwidth]{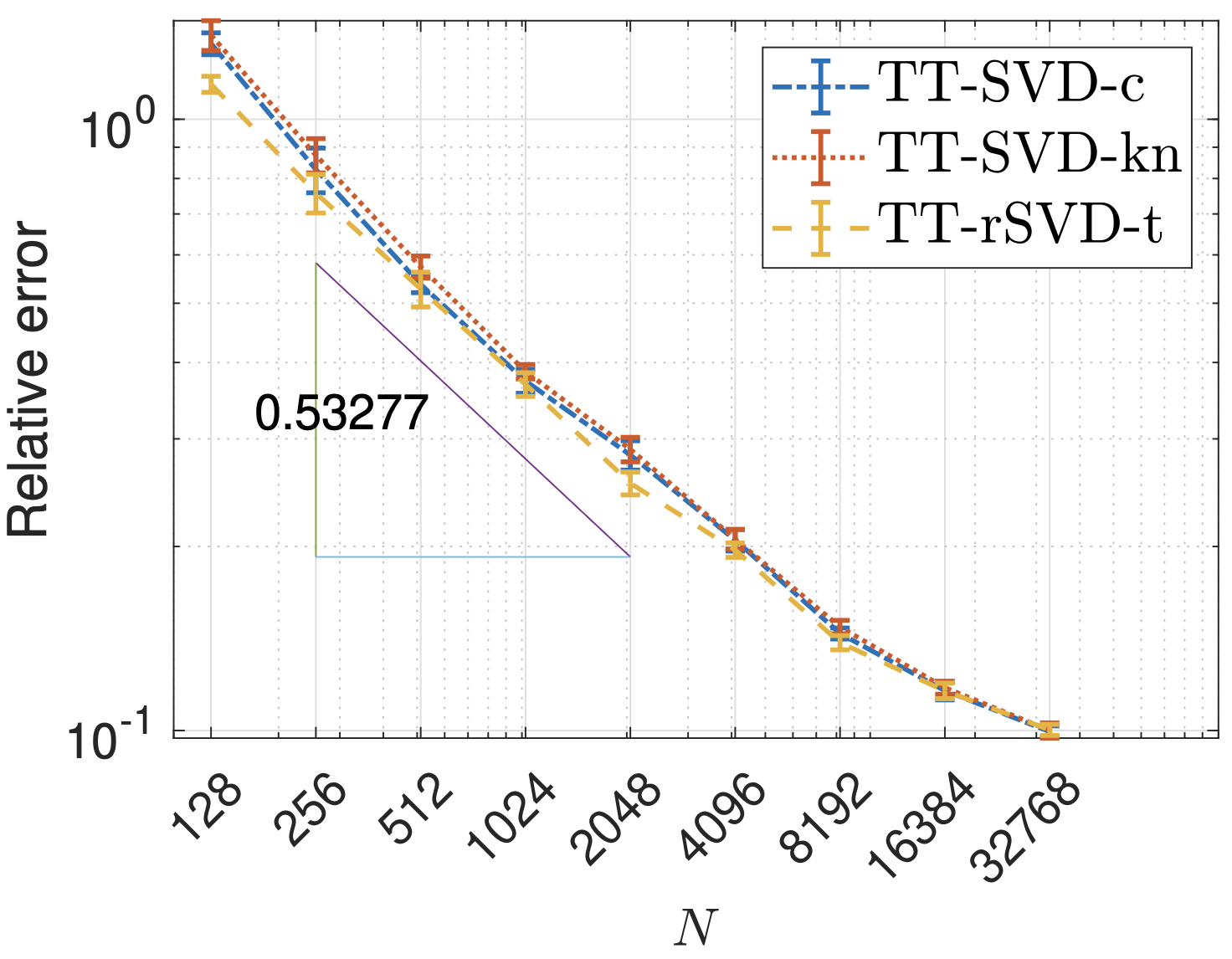}
\includegraphics[width=0.49\textwidth]{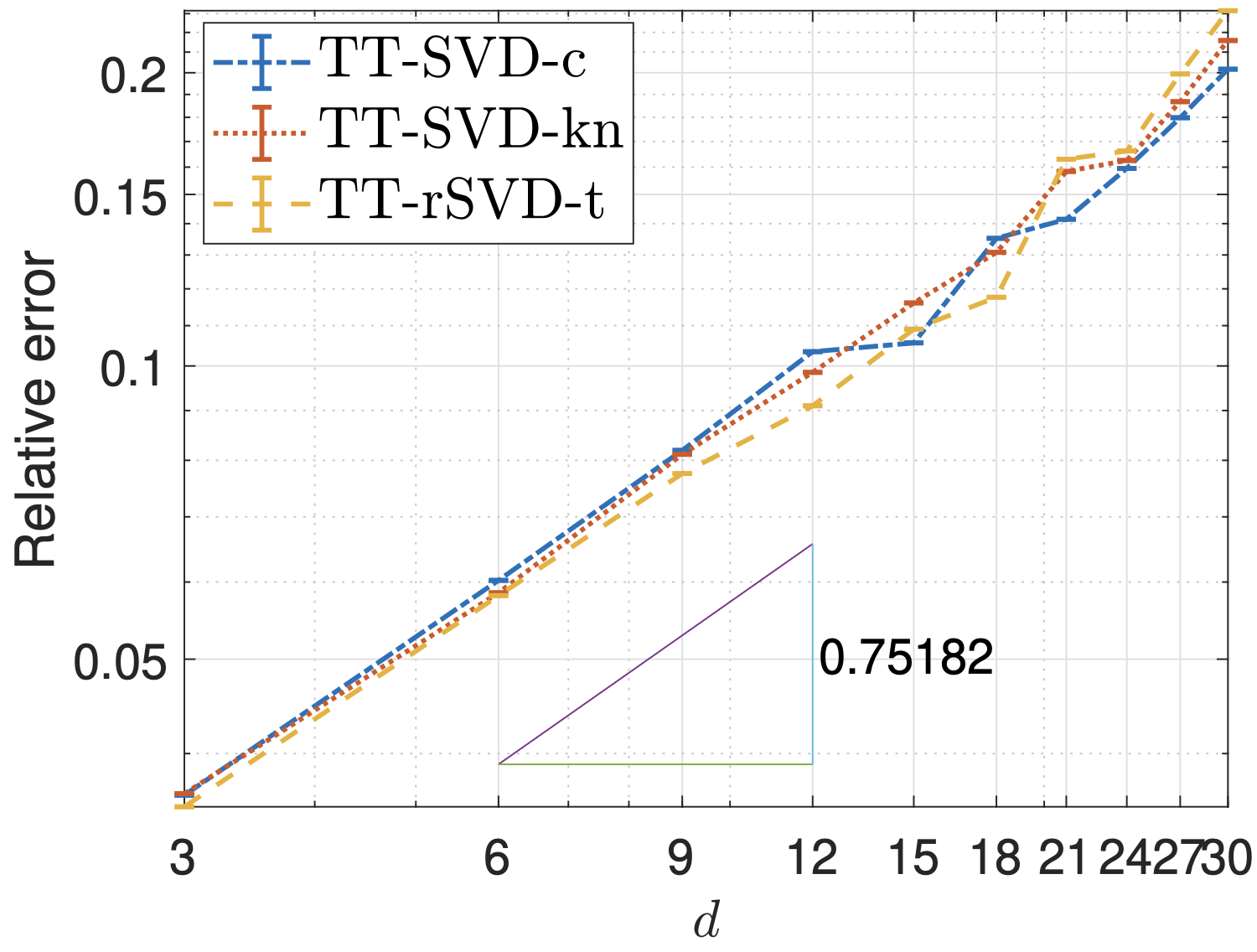}
\caption{{Gaussian mixture model: relative  $L_2$-error for proposed algorithms. Left: $d =10$ and sample number $N=2^7,2^8,\cdots,2^{15}$. Right: $N =10^6$ and dimension $d=3,6,9,\cdots,30$.}}
\label{fig:gm_error}
\end{figure}

We first fix dimension $d = 10$. We draw $N$ i.i.d samples from the mixture model, and use proposed three algorithms to fit a TT approximation for the empirical distribution of these samples. For the algorithms, we set fixed target rank $r = 3$ and use $n=17$ Fourier bases. In the left panel of Figure \ref{fig:gm_error}, we test different sample sizes $N$ and plot the relative  $L_2$-error defined by $$\|\tilde{p} - p^*\|_{\lt}/\|p^*\|_{\lt}$$ in logarithm scale. One can observe that the relative error decays with Monte Carlo rate $\bigO(1/\sqrt{N})$. Here we remark that the ground truth $p^*$ for this Gaussian mixture model can be formulated as a TT without approximation error, and the numerical $L_2$-norm based on quadrature of TTs can be easily computed with $\bigO(d)$ complexity \cite{oseledets2010tt}. 

Next, we fix the sample size $N = 10^6$ and test our algorithms with various dimensions $d$. The  relative  $L_2$-error is plotted in right panel of Figure \ref{fig:gm_error}, which exhibits sublinear growth with increasing dimensionality according to the curves. In particular, with $N=10^6$ samples, our algorithms can achieve two-digit relative $L_2$-error up to $d = 12$. For dimensionality higher than that, we can expect that the error can also be controlled below $0.1$ for a sufficiently large sample size, as suggested by the Monte Carlo convergence observed in the previous experiment.

\begin{figure}[ht!]
\centering
\includegraphics[width=0.49\textwidth]{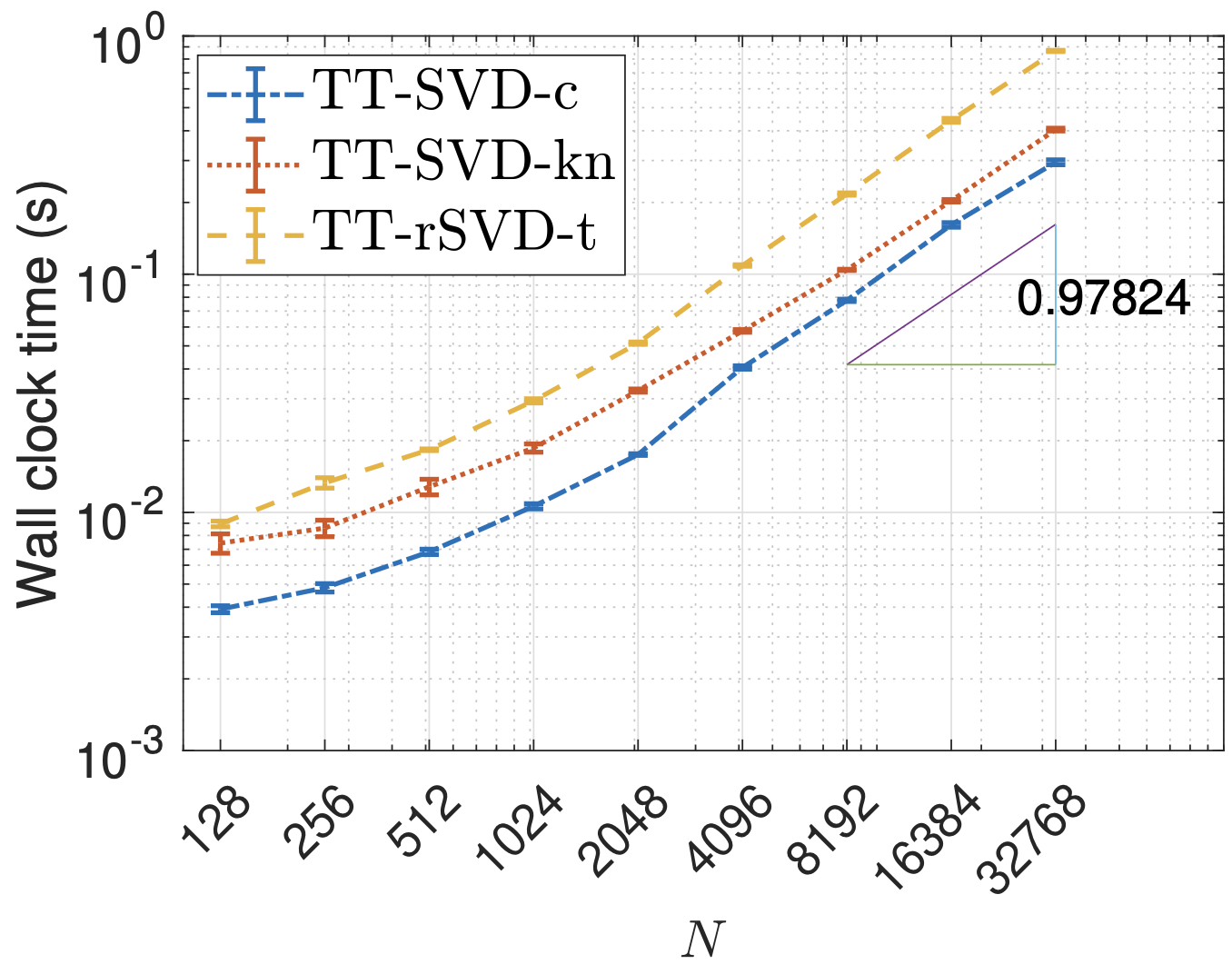}
\includegraphics[width=0.49\textwidth]{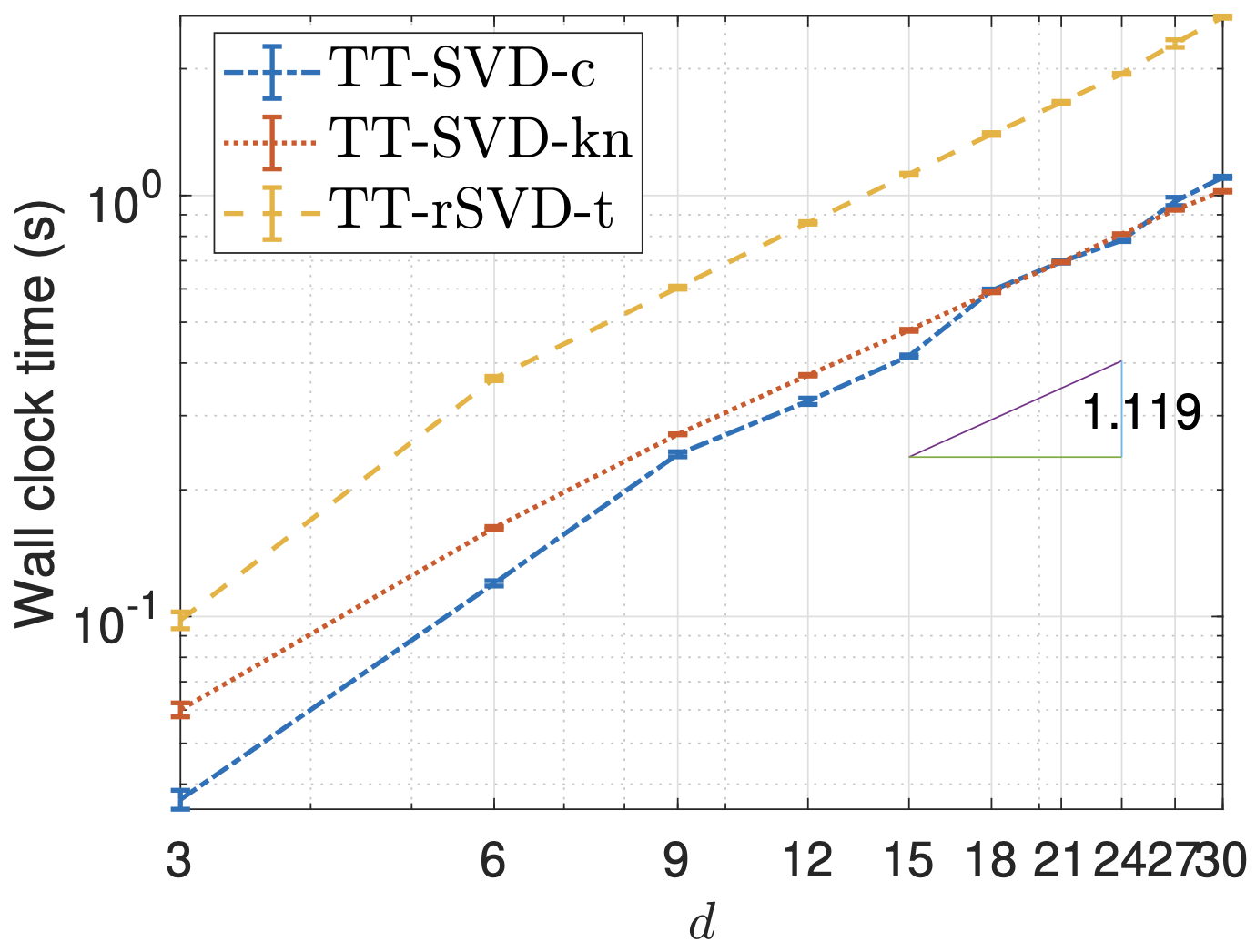}
\caption{Gaussian mixture model: average wall clock time (seconds) for proposed algorithms. Left: $d =5$ and sample number $N=2^7,2^8,\cdots,2^{15}$. Right: $N =10000$ and dimension $d=3,6,9,\cdots,30$.}
\label{fig:gm_time}
\end{figure}

Finally, we verify the time complexity of the algorithms.
In Figure \ref{fig:gm_time}, we plot the wall clock time of three algorithms for fixed $d =5$ with various $N$ (left panel) and fixed $N = 10000$ with various $d$ (right panel). The rates of wall clock time verify the theoretical computational complexity at listed in Table \ref{tab:sketching-table}.

\subsection{Ginzburg-Landau model}
Now we consider a more practical problem in statistical physics as the Ginzburg-Landau model \cite{hoffmann2012GL}. In particular, we consider the density estimation for the Boltzmann distribution 
$$p^*(\vx) = \frac{1}{Z}\mathrm{e}^{-\beta V(\vx)}$$ where $\beta$ is the inverse temperature, $Z$ is the normalization factor such that $\int p^*(\vx) d \vx = 1$, and $V$ is the Ginzburg-Landau potential energy.

\subsubsection{1D Ginzburg-Landau model}\label{sec:num_exp_gl1d}
 We first consider a 1D Ginzburg-Landau model where the  potential energy is formulated as 
\begin{equation}\label{eq:GLpotential}
    V(x_1,\cdots,x_d) = \sum^{d+1}_{i=1} \frac{\lambda}{2} \left( \frac{x_i - x_{i-1}}{h} \right)^2 + \frac{1}{4\lambda}(1 - x_i^2)^2  \text{~for~} \vx \in [-L,L]^d
\end{equation}
where $x_0 = x_{d+1} = 0$, $h = 1/(1+d)$, $\lambda= 0.03$ and $\beta = 1/8$. We consider hypercube domain $[-L,L]^d$ with $L=2.5$ and mesh size $0.05$ in each direction. For the ground truth $p^*$, we use DMRG-cross algorithm (see \cite{savostyanov2011fast}) to get $p^*$ in TT form, which is known to provide highly accurate approximations with arbitrarily prescribed precision. The samples of $p^*$ are generated by running Langevin dynamics for sufficiently long time.

\begin{figure}[ht!]
\centering
\includegraphics[width=0.32\textwidth]{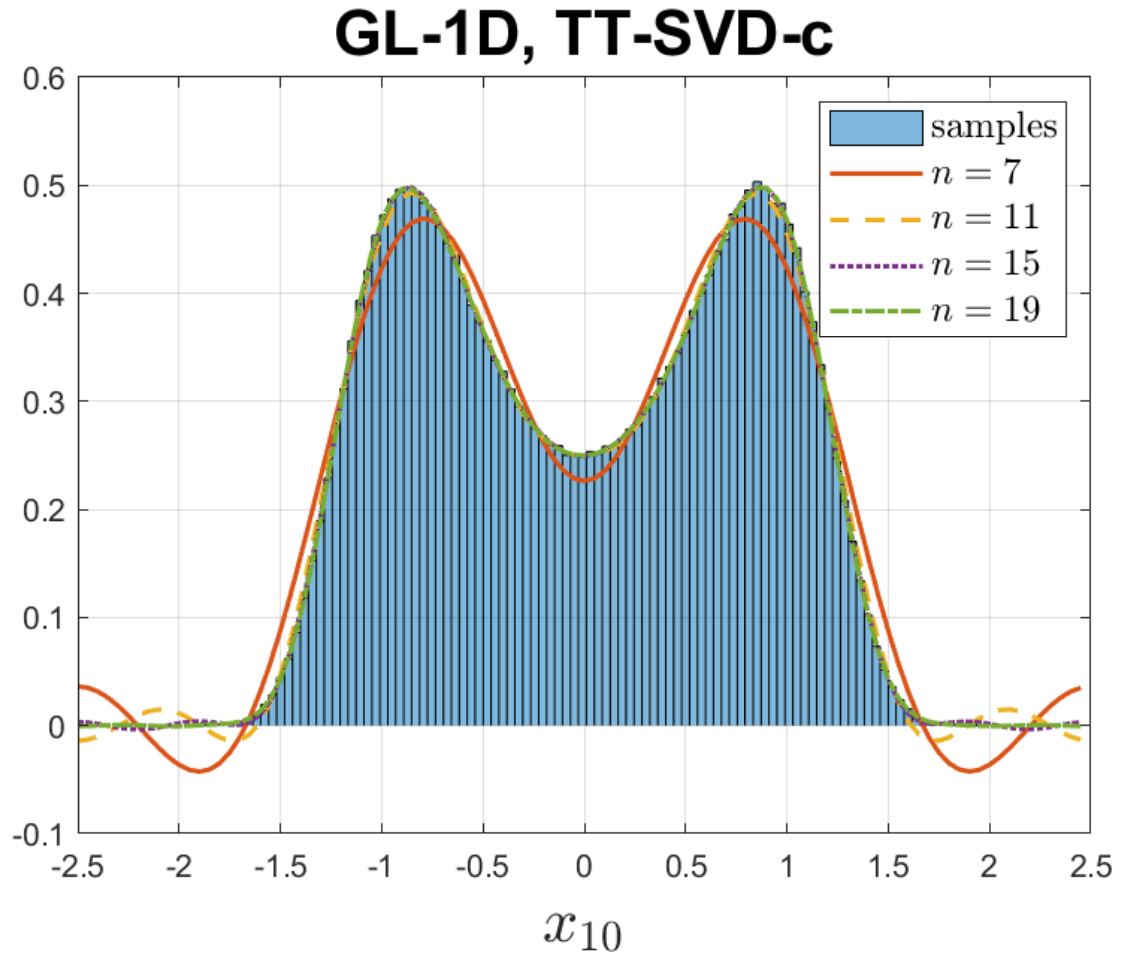}
\includegraphics[width=0.32\textwidth]{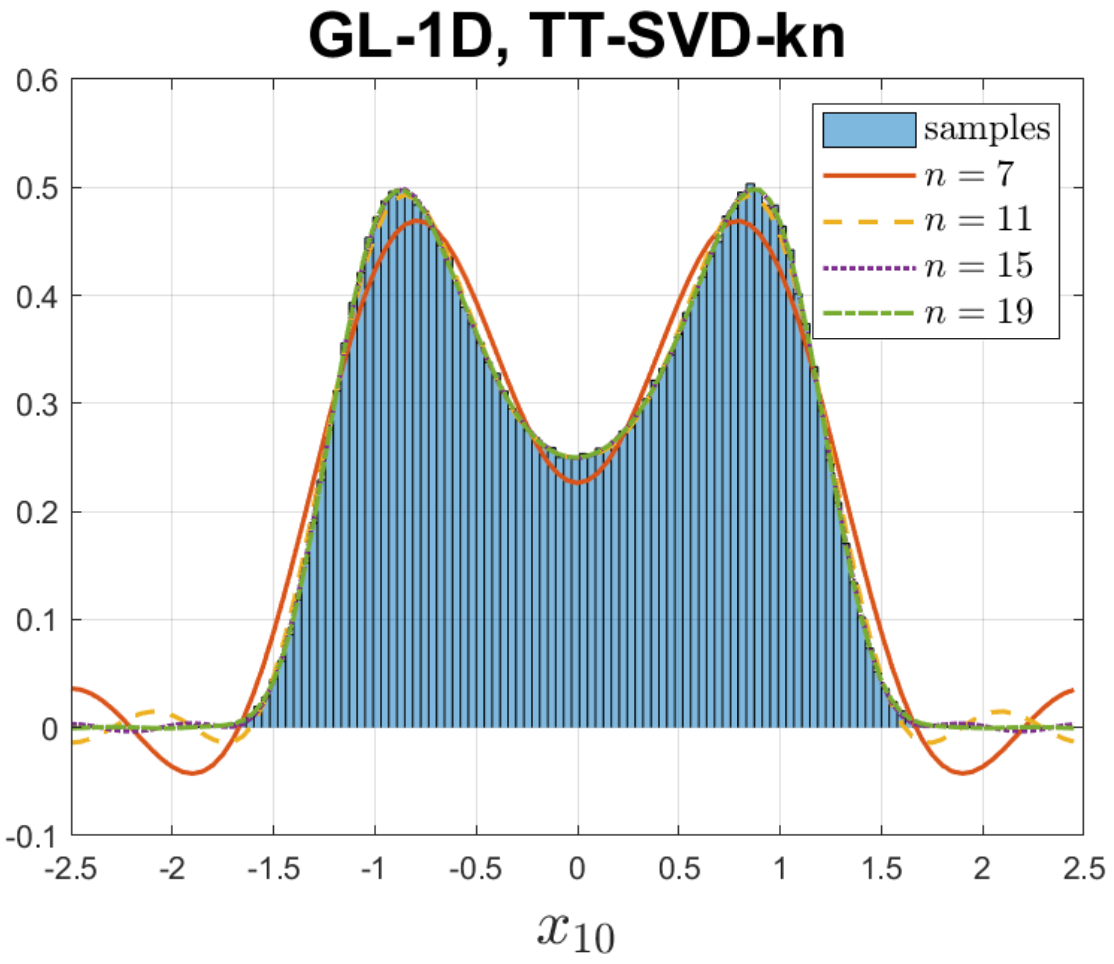}
\includegraphics[width=0.32\textwidth]{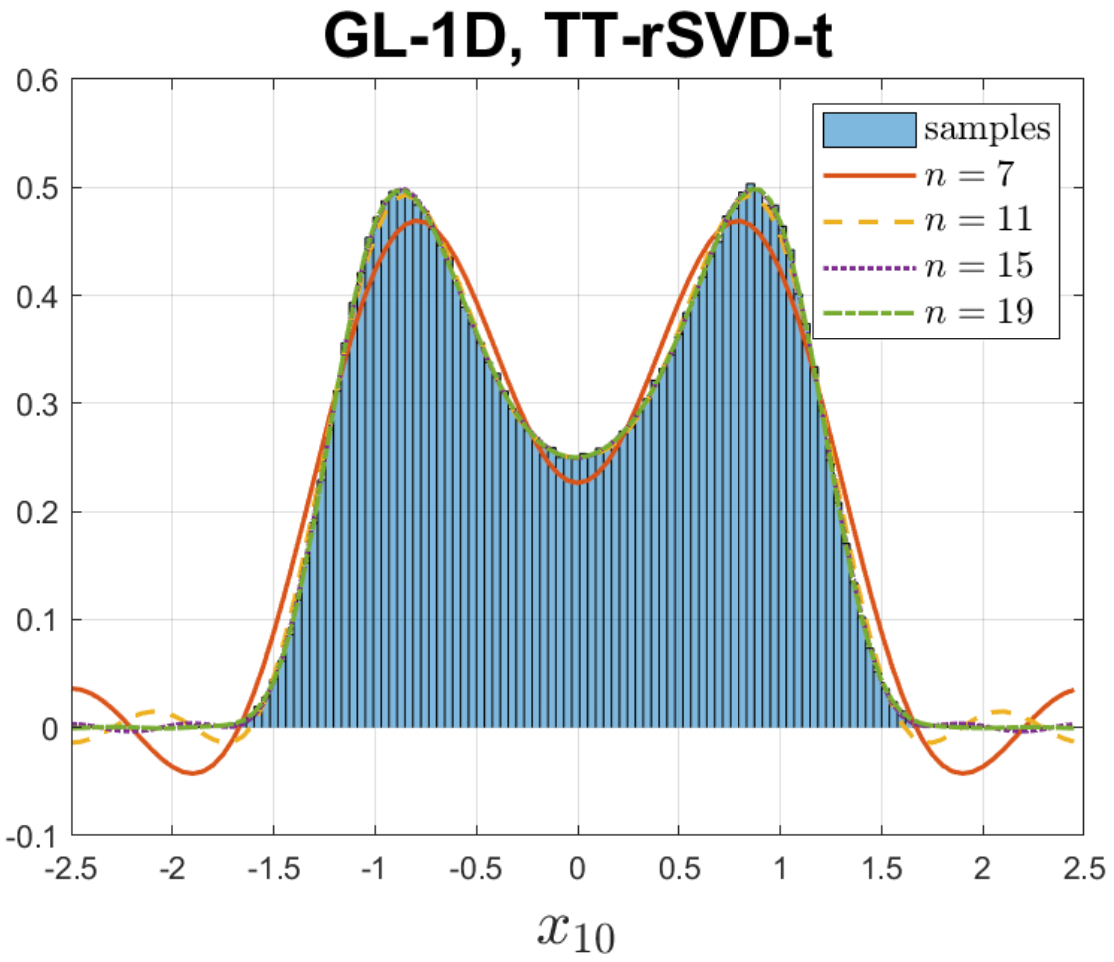}
\caption{1D G-L model: $x_{10}$-marginal of approximated solutions for $d=20$ and $N=10^6$ under different numbers of univariate basis functions used.}
\label{fig:gl1d_mar1}
\end{figure}

In Table \ref{tab:gl1d_err}, we list the relative  $L_2$-error of the approximated densities estimated by $N = 10^6$ samples with different number of univariate basis functions. For this example, we set fixed target rank $r= 4$. One can observe that the approximation improves as we use more basis functions. Such improvement can be further visualized in Figure \ref{fig:gl1d_mar1}, where we plot the marginal distribution of the approximated solutions solved by three algorithms when $d = 20$. All of them give accurate match to the sample histograms when a sufficient number of basis functions are used.

\begin{table}[ht]
\centering
\begin{tabular}{c|ccc|ccc}
\toprule
\multirow{2}{*}{$n$} & \multicolumn{3}{c|}{$d = 5$} & \multicolumn{3}{c}{$d = 10$} \\
\cmidrule(lr){2-4} \cmidrule(lr){5-7}
 & TT-SVD-c & TT-SVD-kn & TT-rSVD-t & TT-SVD-c & TT-SVD-kn & TT-rSVD-t \\
\midrule
7  & 0.2810 & 0.2809  & 0.2810 & 0.4304 & 0.4304 & 0.4305 \\
11 & 0.0782 & 0.0781 & 0.0781 & 0.1365 & 0.1365  & 0.1378 \\
15 & 0.0305 & 0.0304 & 0.0301 & 0.0488 & 0.0483 & 0.0533 \\
19  & 0.0275 & 0.0273 & 0.0289 & 0.0449 & 0.0404 & 0.0526 \\
\bottomrule

\toprule
\multirow{2}{*}{$n$} & \multicolumn{3}{c|}{$d = 15$} & \multicolumn{3}{c}{$d = 20$} \\
\cmidrule(lr){2-4} \cmidrule(lr){5-7}
 & TT-SVD-c & TT-SVD-kn & TT-rSVD-t & TT-SVD-c & TT-SVD-kn & TT-rSVD-t \\
\midrule
7  & 0.6342 & 0.6342 & 0.6343 & 0.7897 & 0.7897 & 0.7898 \\
11 & 0.2169 & 0.2167 &  0.2186 & 0.3017 & 0.3016 & 0.3037 \\
15 &  0.0733 &  0.0740 & 0.0836& 0.1126 & 0.1125 & 0.1158 \\
19  & 0.0613 & 0.0602 & 0.0716 & 0.0884 & 0.0874 & 0.1099 \\
\bottomrule
\end{tabular}
\caption{1D G-L model: comparison of $L_2$ relative errors with $N = 10^6$ and $d=5,10,15,20$ under different numbers of univariate basis functions used.}
\label{tab:gl1d_err}
\end{table}

\subsubsection{2D Ginzburg-Landau model}
Next, we consider the 2D Ginzburg-Landau model which is defined on a $m\times m$ lattice. The potential energy is formulated as  
\begin{equation*}
V(x_{1,1},\cdots,x_{m,m}) =  \frac{\lambda}{2}\sum^{m+1}_{i=1}\sum^{m+1}_{j=1} \left[ \left(\frac{x_{i,j} - x_{i-1,j}}{h}\right)^2 +  \left(\frac{x_{i,j} - x_{i,j-1}}{h}   \right)^2  \right] + \frac{1}{4\lambda} \sum_{i=1}^m \sum_{j=1}^m  (1-x_{i,j}^2)^2 
\end{equation*}
with the boundary condition 
\begin{displaymath}
    x_{0,j} = x_{m+1,j} = 1 , \quad x_{i,0} = x_{i,m+1} = -1 \quad \forall \  i,j =1,\cdots,m
\end{displaymath}
where we choose parameter $h = 1/(1+m)$, $\lambda =0.1$ and inverse temperature $\beta = 1$. We set $L=2$ with mesh size $0.05$ for the domain. We consider a test example where $m = 16$ and thus the total number of dimensionality is $m^2 = 256$. For this 2D model, we apply the PCA pre-processing steps described in Section \ref{sec:pre-process} to order the dimensions and also reduce the dimensionality, where we choose embedding dimension $d'= 100$. In terms of the algorithms, we set target rank $r = 10$, $\alpha = 0.001$ and we use $N=10^4$ samples and $n=21$ Fourier bases.

We evaluate the performance of our fitted tensor-train density estimators as generative models. Specifically, we first fit tensor-train density estimators $\tilde{p}$ from the given $N=10^4$ samples $\{x_i\}_{i=1}^N$ using proposed algorithms, then generate new samples $\{\tilde{x}_i\}_{i=1}^N$ from the learned tensor trains. Sampling each $\tilde{x}_i \sim \tilde{p}$ can be done via tensor-train conditional sampling 
\begin{equation}\label{tt_sampling}
    \tilde{p}(x_1,x_2,\cdots,x_d) =  \tilde{p}(x_1)   \tilde{p}(x_2 | x_1)  \tilde{p}(x_3 | x_2, x_1) \cdots  \tilde{p}(x_d |x_{d-1}, \cdots, x_1) 
\end{equation}
where each 
\begin{displaymath}
    \tilde{p}(x_i |x_{i-1}, \cdots, x_1) = \frac{ \tilde{p}(x_1,\cdots,x_i)}{\tilde{p}(x_1,\cdots,x_{i-1})}
\end{displaymath}
is a conditional distribution of $\tilde{p}$, and $\tilde{p}(x_1),\cdots,\tilde{p}(x_1,\cdots,x_{d-1})$ are the marginal distributions. We point out each conditional sampling can be done in $\bigO(d)$ complexity using the TT structure of $\tilde{p}$ \cite{dolgov2020approximation}.

\begin{figure}[ht!]
\centering
\includegraphics[width=\textwidth]{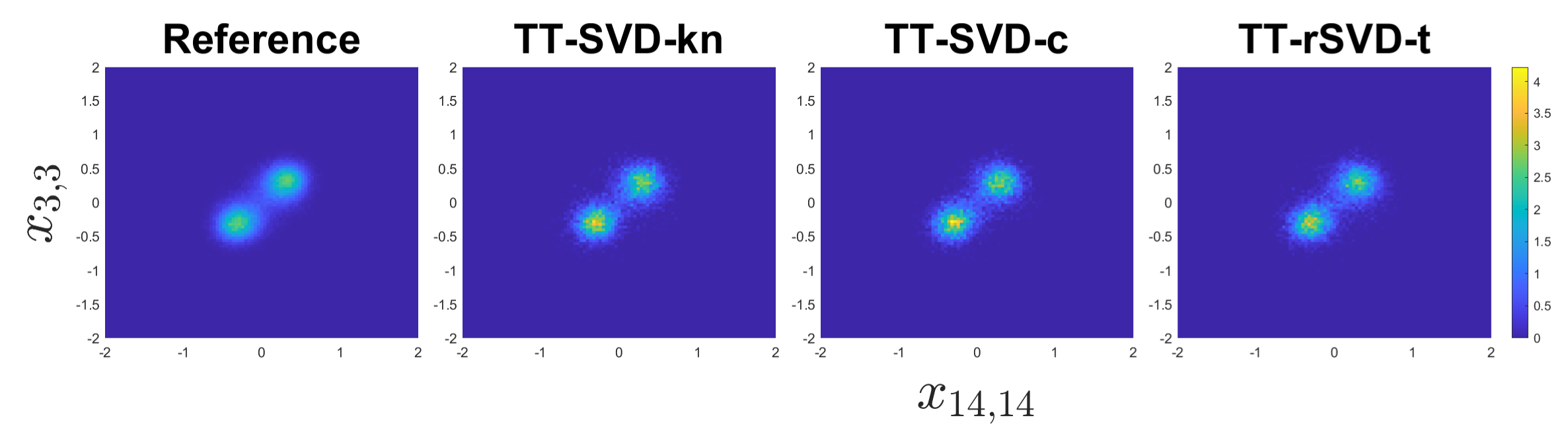}
\includegraphics[width=\textwidth]{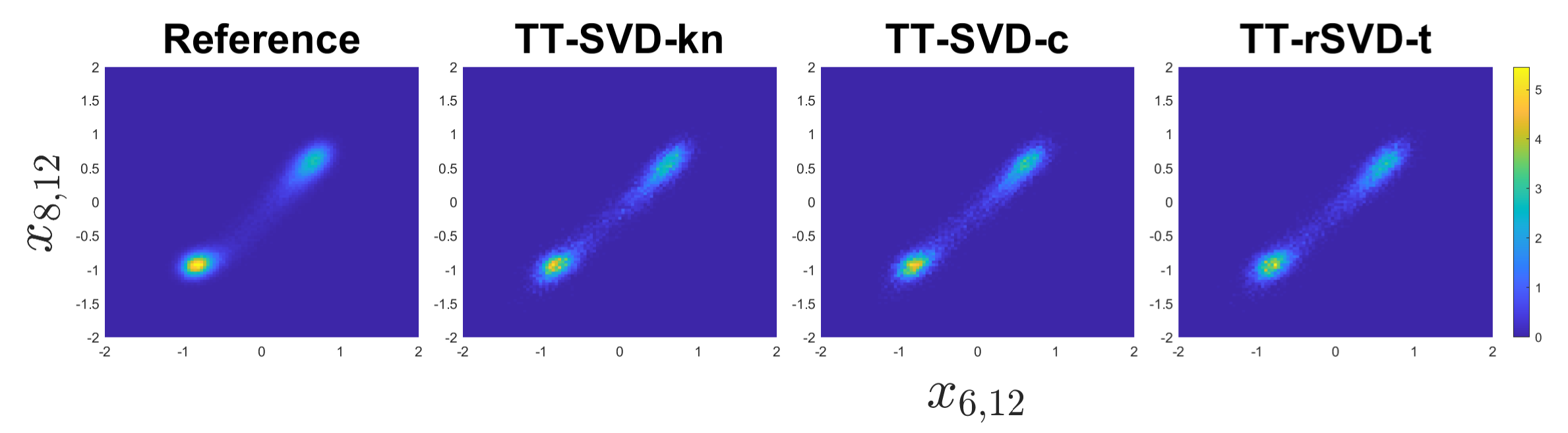}
\caption{2D G-L model on a $16\times 16$ lattice: visualization of $(x_{3,3},x_{14,14})$-marginal (first row) and $(x_{6,12},x_{8,12})$-marginal (second row). Left to right: reference sample histograms, samples from tensor trains fitted by proposed algorithms.}
\label{fig:gl2d_mar2}
\end{figure}

Figure \ref{fig:gl2d_mar2} illustrates the $(x_{3,3},x_{14,14})$-marginal and $(x_{6,12},x_{8,12})$-marginal distributions of the generated samples. For comparison, we also show results for $N^* = 10^5$ reference samples. The visualizations reveal that samples generated from tensor trains learned by all three proposed algorithms can capture local minima accurately. To quantify the accuracy, we compute the second-moment error 
\begin{equation}\label{2ed_mnt_err}
    \|\tilde{X}^T\tilde{X} - {X^*}^T X^* \|_F/\|{X^*}^T X^*\|_F
\end{equation}
where $\tilde{X}\in \mathbb{R}^{N\times d}$ is the matrix formed by the generated samples, and $X^*\in \mathbb{R}^{N^*\times d}$ are reference samples. The second-moment errors for TT-SVD-kn, TT-SVD-c and TT-rSVD-t algorithms are respectively $0.0278, 0.0252, 0.0264$.

\subsection{Comparison with neural network approaches}
 In this section, we compare the performance of generative models fitted by our proposed  tensor-train algorithms with 
one example of neural network-based generative methods (\cite{papamakarios2017masked}). For neural network details, we use 5 transforms and each transform is a 3-hidden layer neural network with width 128. We use Adam optimizer to train the neural network. The experiments are based on selected examples from previous subsections as well as real data from MNIST.

\subsubsection{1D Ginzburg-Landau model}
We first consider the Boltzmann distribution of 1D Ginzburg-Landau potential defined in \eqref{sec:num_exp_gl1d} where we choose model parameters $\lambda =0.005$ and $\beta = 1/20$. We use $N = 10^4$ samples as training data to fit tensor trains, and then generate $N=10^4$ new samples from the learned tensor trains by conditional sampling \eqref{tt_sampling}. 
\begin{figure}[ht!]
\centering
\includegraphics[width=0.6\textwidth]{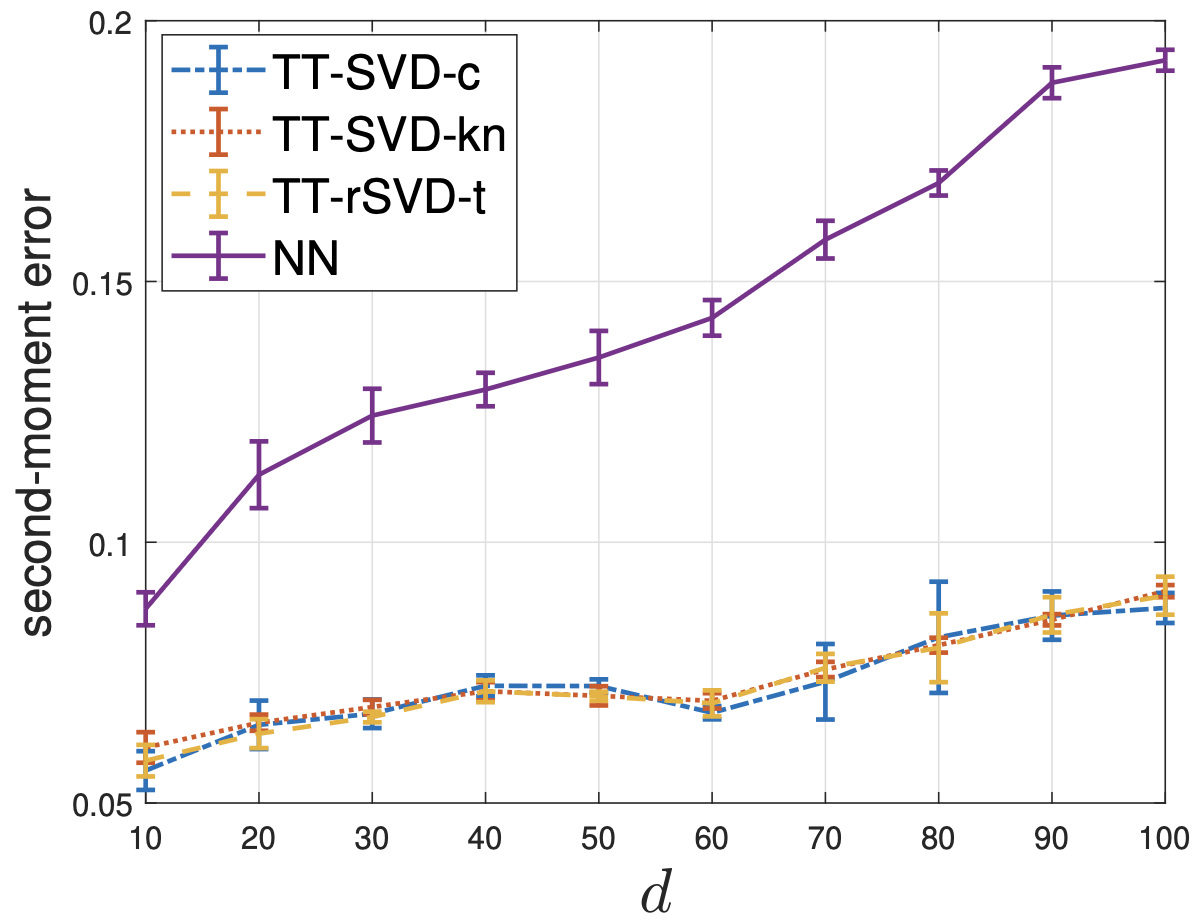} \qquad \qquad 
\caption{Comparison of second-moment errors of generated samples.}
\label{fig:2ed_mnt_err_gl1d}
\end{figure}

\begin{figure}[ht!]
\centering
\includegraphics[width=\textwidth]{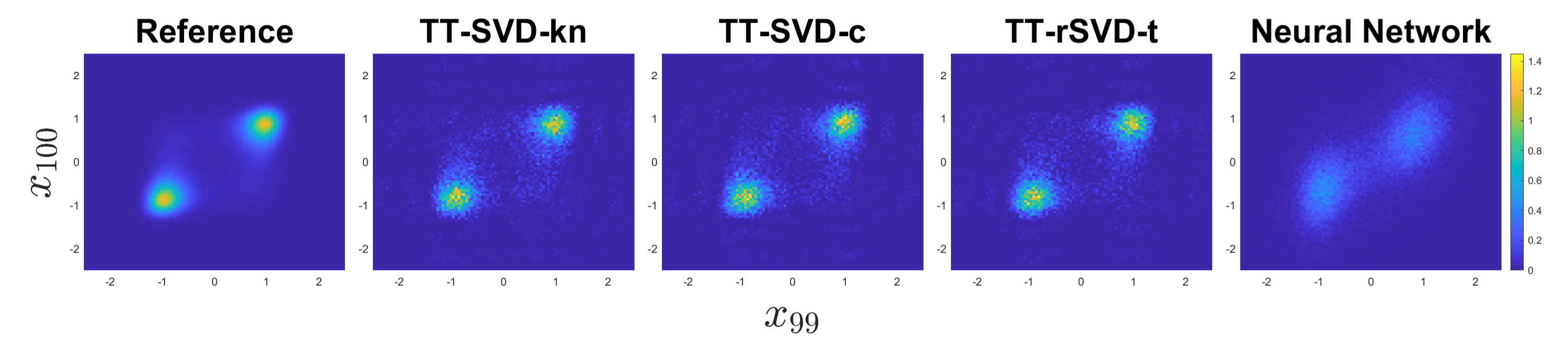}
\includegraphics[width=\textwidth]{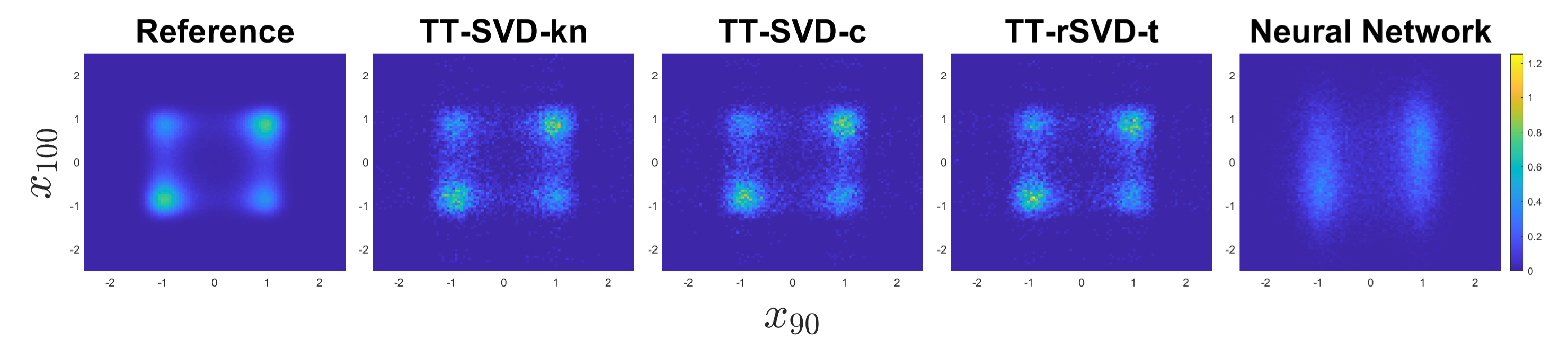}
\caption{1D G-L model for $d=100$: visualization of $(x_{99},x_{100})$-marginal (first row) and $(x_{90},x_{100})$-marginal (second row) of samples from generative models. Left to right: reference samples, samples from tensor trains under three proposed algorithms, and samples from neural network.}
\label{fig:gl1d_mar2}
\end{figure}

In Figure \ref{fig:2ed_mnt_err_gl1d}, we plot the second-moment errors of the samples generated by tensor-train algorithms and a neural network, compared against a reference set of $N^* = 10^6$ samples across various dimensions $d=10,20,\cdots,100$. Here we choose target rank $r=4$ and use $n=17$ Fourier bases for all simulations. The three proposed tensor-train algorithms exhibit similar performance and consistently achieve two-digit accuracy across all dimensions. Moreover, they are considerably more accurate than the neural network. To visualize this accuracy gap, Figure \ref{fig:gl1d_mar2} further shows histograms of two selected marginal distributions from the generated samples for $d=100$. It can be observed that samples generated by the tensor-train estimators more accurately capture local minima compared to those produced by the neural network.

\subsubsection{Real data: MNIST}
\begin{figure}[!htb]
    \centering
    \subfloat[Real images]{{\includegraphics[width=0.4\textwidth]{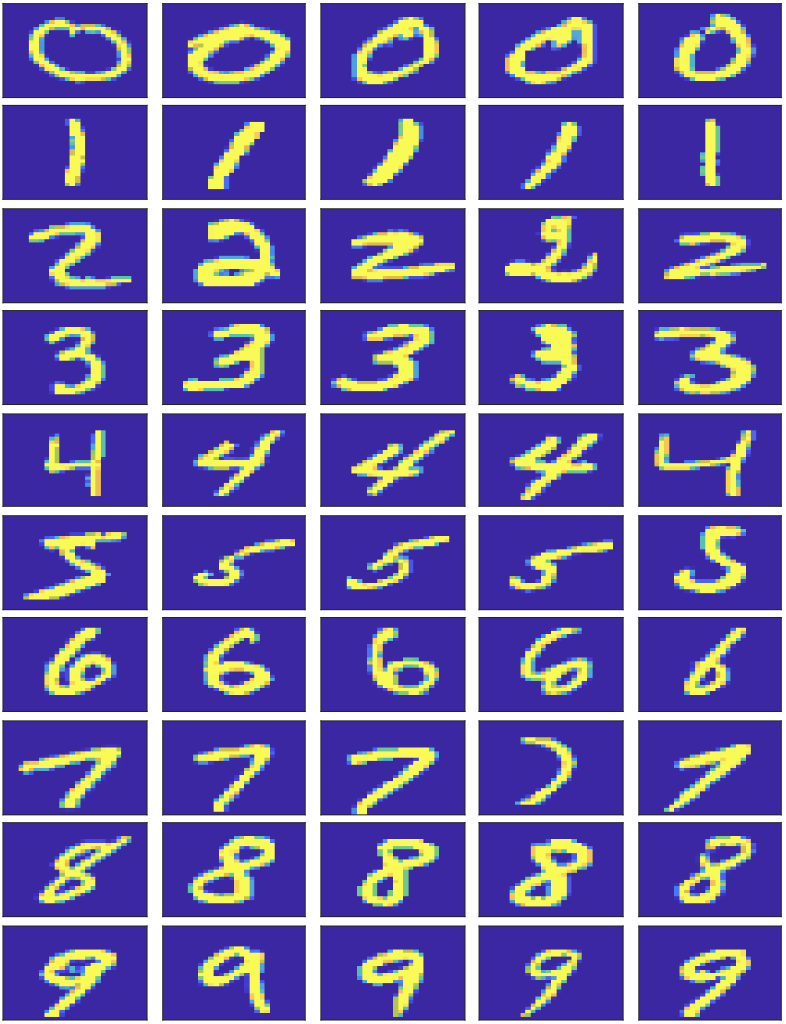}}}\quad 
    \subfloat[Generated images]{{\includegraphics[width=0.4\textwidth]{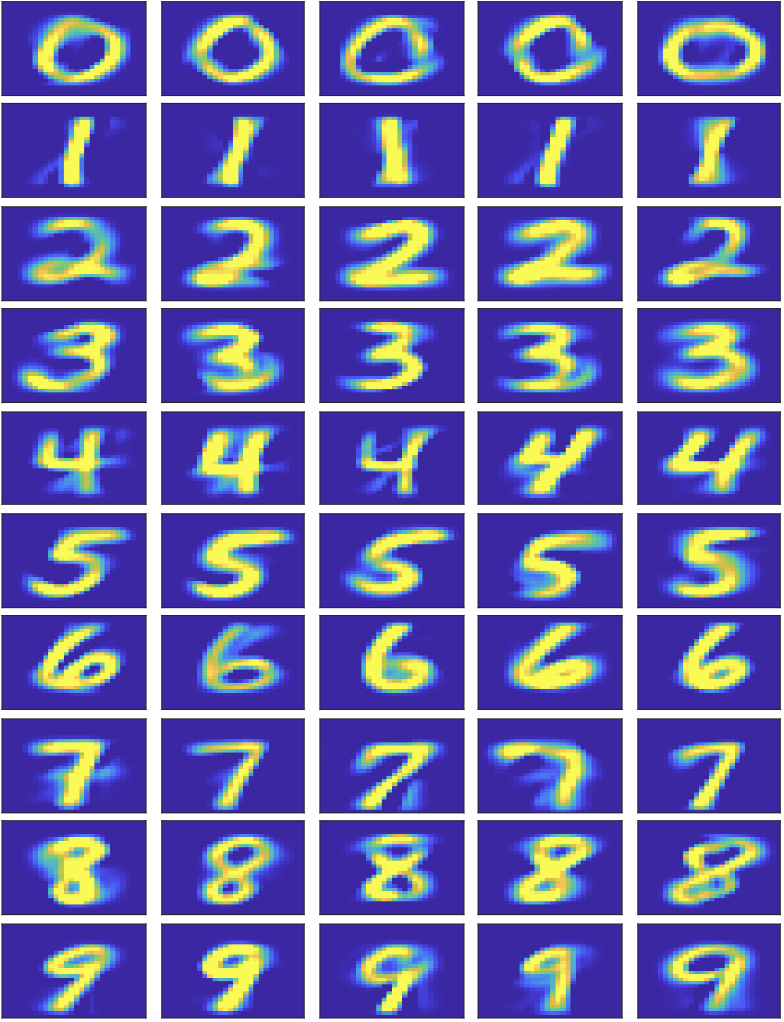}}}
        \caption{Generated and real images from MNIST. For each digit class, five images 
        are randomly selected for illustration. 
        }
    \label{fig:mnist}    
\end{figure}

We also evaluate our algorithms on the MNIST dataset, which consists of 70,000 images of handwritten digits (0 through 9), each labeled with its corresponding digit. Each image is represented as a $28 \times 28$ pixel matrix. To learn the generative models, we treat the dataset as a collection of $70000\times 28^2$ samples. For each digit class, we apply the TT-rSVD-t algorithm with parameters $\alpha = 0.05$, sketching rank $\tilde{r} = 40$, target rank $r = 10$ and $n=15$ Fourier bases, to learn tensor-train representations from the corresponding samples of dimension $d = 28^2 = 784$. Again, the PCA pre-processing steps are applied to reduce dimensionality, where embedding dimension is set to $d'=8$. After training, we independently sample five images for each digit from the learned TT, as shown in Figure~\ref{fig:mnist}(b). For comparison, Figure~\ref{fig:mnist}(a) displays five randomly selected images per digit from the original MNIST dataset. The generated images from the tensor train exhibit well-preserved digit structures and visually plausible variations across classes.


\section{Conclusion}\label{sec:conclusion}
 In this paper, we propose a natural approach to high-dimensional density estimation entirely based on linear algebra, which allows the construction of a density estimator as a low-rank functional tensor train. To this end, we provide four tensor decomposition algorithms. The first one is based on direct TT-SVD algorithm with fast implementation in a density estimation setting. The second algorithm follows a Nyst\"orm approximation. The third approach modifies the first approach by thresholding certain elements in the tensor and employing a hierarchical matrix factorization to accelerate computational efficiency. The last algorithm employs a tensorized sketching technique. Several of these proposed algorithms can achieve linear computational complexity both in dimension $d$ and in sample size $N$. Theoretical analysis demonstrates that tensor-train estimators can also achieve linear scaling of the estimation error with respect to the dimensionality $d$ under some specific conditions. Comprehensive numerical experiments validate the accuracy of all proposed algorithms, in comparison to other density estimators based on neural networks.

\section{Declaration}
\paragraph{Acknowledgment.} Y. Khoo acknowledges partial supports of DMS-2339439, DE-SC0022232 and Sloan foundation. Y. Peng acknowledges partial support of DE-SC0022232.

\bibliographystyle{plainnat}
\bibliography{citations}  
\appendix

\section{Definitions and notations for the proof of \Cref{thm: main theorem}}

Before the proof, we introduce some definitions and notations that will be used in the main proof in \Cref{app:main proof}.

To begin with, we use $c^*  = \Phi^T p^* \in \mathbb R^{ n ^d}$ to present the coefficient tensor of ground truth under a set of orthonormal basis functions $\{\phi_l^j\}_{l=1}^n$, whose entries are given by
\begin{displaymath}
    c^*(l_1,\cdots,l_d) = \langle  p^* , \phi^1_{l_1}\cdots \phi^d_{l_d} \rangle.
\end{displaymath}
We point out that for $p^*$ with 1-cluster assumption as stated in \Cref{assumption:additive model}, it holds that $c^* = \text{diag}(w_K)\Phi^T p^* $ for any $K \geqslant 1$. In the rest of this appendix, we omit the superscripts on the basis functions unless they are needed to indicate the corresponding coordinate of dimensions. Based on the additive formulation of ground truth $p^*$ in \eqref{eq:additive model} and orthogonality of $\{\phi_l^j\}_{l=1}^n$,  $c^*$ can be expressed as 
\begin{equation}\label{eq:additive model for coefficient tensor} 
c^* = \bvec{1}^d +    c_1^* \otimes \bvec{1}^{d-1} + \bvec{1}  \otimes c_2^* \otimes \bvec{1}^{d-2} + \cdots + \bvec{1}^{d-1}\otimes c_d^*,
\end{equation}
where $\bvec{1}^k \in \mathbb R^{n^k}$ represents the $k$ times Kronecker product of $\bvec{1}=  [1,0,\cdots,0]^T \in \mathbb{R}^n$. $c^*_j \in  \mathbb R^{n}$ is a vector whose $\mi$-th coordinate is given by 
  \begin{equation}\label{eq:explicit  expression of coefficients 1} 
 c^*_{j}(l)= \int_0^1 \cdots \int_0^1 p^*(x_1,\cdots,x_d) \phi_\mi(x_j) \dd x_1 \cdots \dd x_d. 
 \end{equation} 
We can further simplify the expression of $c^*_j$ as follows. By the normalization condition of $p^*$, $q_j^*$ in the expression of $p^*$ satisfies $\sum_{j=1}^d \int_0^1  q_j^*(x_j) dx_j   =0$. Without loss of generality, we assume that 
\begin{align}
    \label{eq:normalize the marginal} \int_0^1  q_j^*(x_j) \dd x_j =0.
\end{align}  
This further gives us 
 \begin{align}\label{eq:explicit  expression of coefficients}  c^*_{j}(\mi) = \begin{cases} 
  0,  &\text{ if } \mi =1;
 \\
 \int_0^1   q_j^*(x_j ) \phi_\mi( x_j)   \dd x_j  , & \text{ if } 2\leqslant \mi \leqslant n .
\end{cases} 
\end{align}

\begin{definition} \label{definition:cluster coefficient tensor estimator}
Let  $\mathcal B_{s,d} $ be the \(s\)-cluster basis functions defined in \eqref{k-cluster basis}. Define the $s$-cluster estimator     
\begin{align}\label{eq:cluster basis function estimator}
    \hat p_s (x):= 1  + \sum_{\varphi_1\in \mathcal B_{1,d}}\langle \hat p, \varphi_1\rangle \varphi_1(x) + \cdots +   \sum_{\varphi_s\in \mathcal B_{s,d}}\langle  \hat p , \varphi_s\rangle \varphi_s(x )  
\end{align}   
The corresponding coefficient tensor $\hat c[s]  \in \mathbb R^{n^d}$ of $ \hat p_s$  is   
  \begin{align} \label{eq:coefficient tensor estimator} \hat c[s](\mi_1,\cdots,\mi_d)   = \langle \hat p_s  ,  \varphi_{\mi_1}    \cdots   \varphi_{\mi_d}   \rangle  . 
  \end{align}
 In this case,   $\hat c[s]$ will be  referred  to as 
  the $s$-cluster coefficient   estimator of $c^*$.
\end{definition}

\begin{definition}
    Given matrix $A \in \mathbb{R}^{m\times n}$ and positive integer $r \leqslant \min\{m,n\}$, we define function $\svd: \mathbb{R}^{m\times n}\times \mathbb{Z}^+ \rightarrow \mathbb{O}^{m\times r} $ satisfying 
    \begin{equation}\label{svd_notation}
        \svd(A,r) = U 
    \end{equation}
  where $A = U\Sigma V^T$ is the top-$r$ SVD of the matrix $A$.  
\end{definition}

\section{Proof of the main theorems}
\subsection{Proof of \Cref{thm: main theorem}}
\label{app:main proof}
In this proof, we simply denote $\hat{c}[K]$ as defined in \eqref{eq:coefficient tensor estimator} by $\hat{c}$.
\begin{proof}
Let $\hat G_{1:d }$ be the TT cores of $\tilde{c}$, to bound the error $\|c^* - \hat G_{1:d }  \|_F$, we first use \Cref{lemma:tensor core orthogonal decomposition 2} to formulate it as 
\begin{align} \nonumber 
       & \|c^* - \hat G_{1:d }  \|_F^2 
       \\ \label{eq:main variance term 2}
       =&   \| \hat G_{1:d-1 } \hat G_{1:d -1}^T   c ^{* (d-1)} -  c^{* (d-1)} \|_F^2 
       \\ 
       + &\label{eq:main variance term 1} \| \hat G_{1:d-1 } \hat G_{1:d -1}^T \hat c^{(d-1)}  -\hat G_{1:d-1 } \hat G_{1:d -1}^T   c^{* (d-1)}  \|_F^2      .
    \end{align}
Later we will bound the above two terms by \textbf{Step 1} and \textbf{Step 2} respectively. Before elaborating these two steps, let us introduce two error estimates that will be used frequently in the proof of this theorem. With high probability, we have
\begin{align} \label{eq:sup projected norm}
    & \max_{k\in\{1,\cdots, d\} } \|  ( U_{k-1}^{*T}   \otimes I_n )  (c^{*(k)} -\hat c^{(k)}) V_k ^* \|_F   \leqslant    C_1 \sqrt{\frac{n\log(dn)\|p^*\|_\infty}{N} }    ,
    \\ \label{eq:s cluster frobenius norm}
      & \|      c^{* } -\hat c     \|_F    \leqslant    C_2  \bigg(\sqrt { \frac{     d^K n^{2K}    } { N} }    \bigg). 
\end{align}
Their proofs are provided respectively by \Cref{lemma:maximum projection norm} and \Cref{lemma:Frobenius norm s-cluster}. In the above inequality, $C_1, C_2$ are two sufficiently large universal constants. $U_k^*$ and $V_k^*$ are orthogonal matrices such that $U_k^*U_k^{*T}   \in \mathbb R^{n^k \times n^k} $ is the projection matrix onto the space $\col(c^{* (k)})$, and $V_k^*V_k^{*T}  \in \mathbb R^{n^{d-k} \times n^{d-k}} $ is the projection matrix onto the space $\row(c^{* (k)})$. Due to the additive formulation \eqref{eq:additive model for coefficient tensor} of $c^*$, $U^*_k$ and $V_k^*$ have analytical expressions, shown in \Cref{lemma:additive row and column spaces}. 

\eqref{eq:s cluster frobenius norm} provides an upper bound of the gap between empirical and ground truth coefficient tensor, which scales as $\bigO(d^{K/2})$ upon using $K$-cluster basis for the density estimator. Notably, with the projection $U_k^*$ and $V_k^*$, this gap grows only at $\bigO( \sqrt { d \log(dn)})$ under the \Cref{assumption:additive model} for $\|p^*\|_\infty$ as stated in \eqref{eq:sup projected norm}. This is due to  the 1-cluster structure of $c^*$, which is the key to make our error analysis sharp. 

{\bf Step 1.} To bound  \eqref{eq:main variance term 2}, we show by induction that for $k\in\{1,\cdots, d-1\}$,
\begin{equation}\label{eq:main variance induction general k}
    \|    \hat G_{1:k} \hat G_{1:k}^T   c^{* (k)}  -c^{* (k)}\|_F^2 
    \leqslant    C_3  \frac{n\log(dn)\|p^*\|_\infty \|c^{*} \|_F^2}{N} \sum_{j=1}^k\frac{1}{\sigma_2^2(c^{*(j)})}
\end{equation}
for a universal sufficiently large constant $C_3$ with high probability, where $\sigma_{2}(c^{*(j)})$ is the second largest singular value of the unfolding matrix $c^{*(j)}$. 

We prove $k=1$ first:
\begin{align}
   \|    \hat G_{1} \hat G_{1}^T   c^{* (1)}  -c^{* (1)}\|_F^2  = \|    \hat G_{1} \hat G_{1}^T   c^{* (1)}  - U_1^*{U_1^*}^T c^{* (1)}\|_F^2 \leqslant \|    \hat G_{1} \hat G_{1}^T    - U_1^*{U_1^*}^T\|_{\op}^2 \|c^{*}\|_F^2
\end{align}
We use the matrix perturbation theory in \Cref{corollary:Davis Khan 2} to bound $\|\hat G_{1} \hat G_{1}^T    - U_1^*{U_1^*}^T  \|_{\op}$. The spectrum gap assumption in \Cref{corollary:Davis Khan 2} holds since $\|       \hat c^{(1)} -c^{*(1)}     \|_F    \leqslant    O  \bigg(\sqrt { \frac{     d^K n^{2K}    } { N} }    \bigg)$ and $O  \bigg(\sqrt { \frac{     d^K n^{2K}    } { N} }    \bigg) \leqslant \sigma_2(c^{*(1)})$ by \eqref{eq:s cluster frobenius norm} and using a sufficiently large $N$ such that \eqref{eq:snr} holds. Therefore, we have
\begin{equation}
  \|\hat G_{1} \hat G_{1}^T - U_1^*{U^*_1}^T \|_{\op} \leqslant C_4(\frac{\|(\hat c^{(1)} -c^{*(1)}){V_1^*}^T \|_{\op}}{\sigma_2(c^{*(1)})} + \frac{\|\hat c^{(1)} -c^{*(1)} \|^2_{\op}}{\sigma^2_2(c^{*(1)})}).
\end{equation}
for some universal constant $C_4$. According to \eqref{eq:sup projected norm}, the first term
\begin{equation}
\frac{\|(\hat c^{(1)} -c^{*(1)}){V_1^*}^T \|_{\op}}{\sigma_2(c^{*(1)})}\leqslant  C_1 \sqrt{\frac{n\log(dn)\|p^*\|_\infty}{N \sigma_2^2(c^{*(1)})} } 
\end{equation} 
and the second term can be bounded by the first term upon using sufficiently large $N$ such that \eqref{eq:snr} holds. Therefore, 
\begin{equation}\label{core1_est}
    \|    \hat G_{1} \hat G_{1}^T   c^{* (1)}  -c^{* (1)}\|_F^2 \leqslant  C_3 \frac{n\log(dn)\|p^*\|_\infty \|c^{*} \|_F^2}{N \sigma_2^2(c^{*(1)})}   
\end{equation}
for some universal constant $C_3$.

Now suppose \eqref{eq:main variance induction general k} holds for $k$.  Then for $k+1$, 
\begin{align} \nonumber 
    & \| \hat G_{1:k +1} \hat G_{1:k+1}^T   c ^{*(k+1) }  -   c^{* (k+1)}  \|_F^2   \\ \nonumber 
    =  &  \|\hat G_{1:k+1}\hat G_{1:k+1} ^T c^{*(k+1)}   - ( \hat G_{1:k }\otimes I_n )  (  \hat G_{1:k } \otimes I_n  ) ^T c^{* (k+1 ) }  \|_F^2  +\|\hat G_{1:k }\hat G_{1:k } ^T c^{*(k )} -c^{*(k)}  \|_F ^2 
    \\ \label{eq:variance induction term 1}
    \leqslant &  \|\hat G_{1:k+1}\hat G_{1:k+1} ^T c^{*(k+1)}   -  (  \hat G_{1:k }\otimes I_n )  ( \hat G_{1:k } \otimes I_n  )  ^T c^{* (k+1 ) }  \|_F^2  
       \\ \label{eq:variance induction term 2}  + & 
C_3  \frac{n\log(dn)\|p^*\|_\infty \|c^{*} \|_F^2}{N} \sum_{j=1}^k\frac{1}{\sigma_2^2(c^{*(j)})}  .
\end{align}
where the equality follows from 
  \Cref{lemma:tensor core orthogonal decomposition}, and the inequality follows from induction \eqref{eq:main variance induction general k}. 
 It suffices to  bound \eqref{eq:variance induction term 1}.
Note that 
\begin{align} \nonumber 
     \eqref{eq:variance induction term 1} 
    = &\| (\hat G_{1:k } \otimes I_n)  \hat G_{k+1} \hat G_{k+1} ^T (\hat G_{1:k}  \otimes I_n)^T  c^{*(k +1)} - (  \hat G_{1:k }\otimes I_n  )    ( \hat G_{1:k } \otimes I_n  ) ^T c^{* (k+1 ) }  \|_F^2 
    \\ \nonumber 
    = & \| (\hat G_{1:k } \otimes I_n) \big [   \hat G_{k+1} \hat G_{k+1} ^T (\hat G_{1:k}  \otimes I_n)^T  c^{*(k +1)} -   ( \hat G_{1:k } \otimes I_n )   ^T c^{* (k+1 ) }  \big]  \|_F^2 
    \\ \nonumber  
    \leqslant &  \| \hat G_{1:k } \otimes I_n   \|_{\op } ^2 \|   \hat G_{k+1} \hat G_{k+1} ^T (\hat G_{1:k}  \otimes I_n)^T  c^{*(k +1)} -    (  \hat G_{1:k } \otimes I_n )   ^T c^{* (k+1 ) }    \|_F^2 
    \\ \label{eq:induction step 1 term 1}
    =&\|   \hat G_{k+1} \hat G_{k+1} ^T (\hat G_{1:k}  \otimes I_n)^T  c^{*(k +1)} -   ( \hat G_{1:k } \otimes I_n  )  ^T c^{* (k+1 ) }    \|_F^2 .
\end{align}
\
\\
{\bf Step 1A.} Define the matrix formed by principal left singular vectors:
 \begin{align} \label{eq:the population in induction}
\tilde G _{k+1} = \svd  \big(  (  \hat G_{1:k }\otimes I_n  )  ^T    c ^{*(k+1)},2 \big) \in \mathbb O^{2n\times 2}.
 \end{align}
where the $\svd$ function is defined in \eqref{svd_notation}. Since $\rank ((  \hat G_{1:k-1}\otimes I_n  )  ^T    c ^{*(k)}) =2 , $ it follows that 
\begin{align}\label{eq:projection indentity in induction} ( \hat G_{1:k } \otimes I_n  )  ^T c^{* (k+1 ) } = \tilde G _{k+1}\tilde G^T_{k+1}  ( \hat G_{1:k } \otimes I_n  )  ^T c^{* (k+1 ) } . \end{align}
Therefore by \eqref{eq:induction step 1 term 1} and \eqref{eq:projection indentity in induction},
\begin{align} \nonumber 
   \eqref{eq:variance induction term 1}
    = &\|   \hat G_{k+1} \hat G_{k+1} ^T (\hat G_{1:k+1}  \otimes I_n)^T  c^{*(k +1)} -   \tilde G _{k+1}\tilde G^T_{k+1}  ( \hat G_{1:k } \otimes I_n  )  ^T c^{* (k+1 ) }     \|_F^2 
    \\ \nonumber
    \leqslant & \|  \hat G_{k+1} \hat G_{k+1} ^T  - \tilde G _{k+1}\tilde G^T_{k+1}  \|_{\op}  ^2 \| \hat G_{1:k } \otimes I_n    \|_{\op } ^2  \| c^{*(k+1)}\|_F ^2 
    \\ \label{eq:induction step 1 term 1 bound 1} 
    = & \|  \hat G_{k+1} \hat G_{k+1} ^T  - \tilde G _{k+1}\tilde G^T_{k+1}  \|_{\op} ^2    \| c^{* }\|_F ^2.
\end{align}
Perturbation of projection matrices $\|  \hat G_{k+1} \hat G_{k+1} ^T  - \tilde G _{k+1}\tilde G^T_{k+1}  \|_{\op}$ will again be estimated using \Cref{corollary:Davis Khan 2}. More specifically, in \Cref{corollary:Davis Khan 2}, we consider $\hat{A} =  (  \hat G_{1:k}\otimes I_n  )  ^T  \hat c^{(k+1)}$ and $A = \hat G_{1:k}\otimes I_n  )  ^T c ^{*(k+1)}$, and $\hat G_{k+1} \hat G_{k+1} ^T$ and $\tilde G _{k+1}\tilde G^T_{k+1}$ are respectively the projection matrices of $ \hat A$ and $ A$. In the subsequent {\bf Step 1B}, we verify the conditions that should be satisfied in \Cref{corollary:Davis Khan 2}. 
\
\\
\\
{\bf Step 1B.}  
By the induction assumption \eqref{eq:main variance induction general k}, we have 
\begin{equation}   \label{eq:induction bound exact order 2}
\begin{split}
    \| \hat G_{1:k } \hat G_{1:k}^T  c ^{*(k)}  -   c^{* (k)}  \|_F^2 
    \leqslant & \     C_3  \frac{n\log(dn)\|p^*\|_\infty \|c^{*} \|_F^2}{N} \sum_{j=1}^k\frac{1}{\sigma_2^2(c^{*(j)})} 
    \\ 
    \leqslant & \ \frac{C_3n\log(dn) (1+B^*d) C_c d  }{N}\sum_{j=1}^{d} \frac{   d }{ c^2_\sigma j(d-j)   }  \leqslant \left( \frac{B^*}{b^*} \right)^4\frac{ C_4  n d^2  \log^2(dn  ) }{N}  ,
    \end{split}
\end{equation}
 where $C_4$ is a universal constant, $B^*$ is given in \Cref{assumption:additive model}. The second inequality above follows from 
upper bound for $\|c^*\|_F^2$ in \Cref{corollary:frobenius norm of coefficient tensor} and the lower bound for singular value in \Cref{lemma:lower bounds for singular spectral values}, and the third inequality follows from the fact that 
 $ \sum_{j=1}^{d} \frac{   1 }{ j(d-j)   } =\bigO( \frac{\log(d)}{d})$. Consequently,
 \begin{equation} \label{eq:induction bound exact order}
  \| ( \hat G_{1:k }\otimes I_n )  (\hat G_{1:k} \otimes I_n) ^T  c ^{*(k+1)}  -   c^{* (k+1)}  \|_F^2  
  =     \| \hat G_{1:k } \hat G_{1:k}^T  c ^{*(k)}  -   c^{* (k)}  \|_F^2  
     \leqslant \frac{ \sigma_2^2(c^{* (k+1)})}{4} ,
\end{equation}  
where  the inequality follows from the lower bound for the singular value $c^{* (k)} $ given in \Cref{lemma:lower bounds for singular spectral values} and the assumption \eqref{eq:snr} on sample size $N$. Therefore, by \Cref{lemma:left projection preserves matrix structure}, 
\begin{align}
   \bullet  \  &\rank ((  \hat G_{1:k }\otimes I_n  )  ^T    c ^{*(k+1)}) =2 ; 
    \\
   \bullet  \ \label{eq:induction row space invariant}  &\text{$V^*_{k+1} V_{k+1}^{*T}$ is the projection matrix on to $ \row ( (  \hat G_{1:k}\otimes I_n  )  ^T c ^{*(k+1)}) $;  }
    \\
    \label{eq:lower bound of singular values after projection} \bullet  \ &\sigma_2((\hat G_{1:k} \otimes I_n) ^T  c ^{*(k+1)} )\geqslant  \frac{\sigma _2( c^{*(k+1)})}{2}.
\end{align} 

 We start to apply \Cref{corollary:Davis Khan 2}
to $\hat G_{k+1}$ and $\tilde G_{k+1}$,  observe that 
\begin{align}\nonumber
   &  \| (  \hat G_{1:k}\otimes I_n  )  ^T c ^{*(k+1)} - (  \hat G_{1:k}\otimes I_n  )  ^T  \hat c^{(k+1)}  \|_F \leqslant    \|  \hat G_{1:k}\otimes I_n  \|_{\op } \|c ^{*(k+1)}  -\hat c^{(k+1)}\| _F 
    \\
 =  & \label{eq:bound for term 1 of Davis Khan}
       \|       \hat c -c^{* }     \|_F    \leqslant    C_2  \bigg(\sqrt { \frac{     d^K n^{2K}    } { N} }    \bigg) \leqslant 8   \sigma_2(c^{*(k+1)} )  \leqslant 4  \sigma_2((\hat G_{1:k} \otimes I_n) ^T  c ^{*(k)} ) ,
\end{align}   
where the equality follows from the fact that   $\hat G_{1:k }\otimes I_n \in \mathbb O^{n^{k+1}\times 2n}  $, the second inequality follows from  
\eqref{eq:s cluster frobenius norm}, the third inequality follows from 
\eqref{eq:snr}, and the last inequality follows from \eqref{eq:lower bound of singular values after projection}.  
Since \eqref{eq:bound for term 1 of Davis Khan} and  \eqref{eq:induction row space invariant} hold,  by \Cref{corollary:Davis Khan 2}
\begin{align}\nonumber 
    & \| \tilde G_{k+1}\tilde G_{k+1}^T - \hat G_{k+1}\hat G_{k+1} \|_{\op} 
    \\
    \leqslant  \label{eq:davis khan in induction}&\frac{ \| (  \hat G_{1:k}\otimes I_n  )  ^T  ( c ^{*(k+1)} -\hat c^{(k+1) }) V_{k+1}^*   \|_{\op}  }{ \sigma_2  ( ( \hat G_{1:k } \otimes I_n  )  ^T c^{* (k+1 ) } ) } 
    \\ + &  \label{eq:davis khan in induction 2} \frac{ \| (  \hat G_{1:k}\otimes I_n  )  ^T  ( c ^{*(k+1)} -\hat c^{(k+1) })     \|_{\op}^2  }{ \sigma_2^2 ( ( \hat G_{1:k } \otimes I_n  )  ^T c^{* (k+1 ) } )^2  }  
\end{align}
We will bound \eqref{eq:davis khan in induction} in { \bf Step 1C} and \eqref{eq:davis khan in induction 2} in { \bf Step 1D}.
\ 
\\
\\
{ \bf Step 1C.} For  \eqref{eq:davis khan in induction},  we have
\begin{align} \nonumber 
   & \|  (  \hat G_{1:k}\otimes I_n  )  ^T  ( c ^{*(k+1)} -\hat c^{(k+1) }) V_{k+1}^*  \|_{\op} 
   \\  \nonumber 
   = &  \|  (  \hat G_{1:k}\otimes I_n   )    (  \hat G_{1:k}\otimes I_n   )  ^T  ( c ^{*(k+1)} -\hat c^{(k+1) }) V_{k+1}^*  \|_{\op} 
    \\ \nonumber  
    = & \|  (  \hat G_{1:k}  \hat G_{1:k} ^T \otimes I_n   )      ( c ^{*(k+1)} -\hat c^{(k+1) }) V_{k+1}^*  \|_{\op} 
    \\ \label{eq:induction step 1c term 1}
    \leqslant  & \|   (  \{   U_{k}^* U_{k}^{*T}  -  \hat G_{1:k}  \hat G_{1:k}^T  \}     \otimes I_n    )      ( c ^{*(k+1)} -\hat c^{(k+1) }) V_{k+1}^*  \|_{\op}  
    \\ \label{eq:induction step 1c term 2}
    + & \|    (  U_{k}^* U_{k}^{*T}   \otimes I_n )      ( c ^{*(k+1)} -\hat c^{(k+1) }) V_{k+1}^*  \|_{\op}  .
\end{align}
Note that
\begin{equation}
 \begin{split}
\label{eq:eq:induction step 1c term 2 bound}
\eqref{eq:induction step 1c term 2}  = &  \|    (  U_{k}^*\otimes I_n)( U^*_{k} \otimes I_n )^{T}      ( c ^{*(k+1)} -\hat c^{(k+1) }) V_{k+1}^*  \|_{\op}       
\\ 
 = & \|   ( U^*_{k} \otimes I_n )^{T}      ( c ^{*(k+1)} -\hat c^{(k+1) }) V_{k+1}^*  \|_{\op} \\
 \leqslant &  C_1 \sqrt{\frac{n\log(dn)\|p^*\|_\infty}{N} } \leqslant       \frac{\sqrt{ C_3 }}{8 }  \sqrt{\frac{n\log(dn)\|p^*\|_\infty}{N} } ,
  \end{split}
\end{equation} 
where the last inequality follows from \eqref{eq:sup projected norm} and we choose $C_3$ in \eqref{core1_est} to be at least $64C_1^2$. In addition, 
\begin{align} \nonumber 
    \eqref{eq:induction step 1c term 1}  \leqslant &\|     \{   U_{k}^* U_{k}^{*T}  -  \hat G_{1:k}  \hat G_{1:k} ^T    \}     \otimes I_n       \| _{\op}   \|      c ^{*(k+1)} -\hat c^{(k+1) }  \|_{\op} \| V_{k+1}^*  \|_{\op}  
    \\  \nonumber 
    = &  \|        U_{k}^* U_{k}^{*T}  -  \hat G_{1:k}  \hat G_{1:k}     ^T     \| _{\op}   \|      c ^{* } -\hat c   \|_{\op}  
    \leqslant   \|        U_{k}^* U_{k}^{*T}  -  \hat G_{1:k}  \hat G_{1:k}  ^T        \| _{\op}
      C_2  \bigg(\sqrt { \frac{     d^K n^{2K}    } { N} }    \bigg)
      \\ \label{eq:induction step 1c term 1 third bound}
      \leqslant  &C_4  \frac{\sqrt{nB^{*3}\|p^*\|_\infty}}{ b^{*2}\sigma _2( c^{*(k)})} \sqrt{\frac{d\log(dn)}{N}} 
      C_2  \bigg(\sqrt { \frac{     d^K n^{2K}    } { N} }    \bigg) 
        \leqslant  \frac{\sqrt{ C_3 }}{8 }  \sqrt{\frac{n\log(dn)\|p^*\|_\infty}{N} },
\end{align}
where the second  inequality follows from \eqref{eq:s cluster frobenius norm}, the third inequality in which $C_4$ is a universal constant follows from \Cref{lemma:subspace consistency give tensor consistency}, and the last inequality holds for sufficiently large $N$ in \eqref{eq:snr}.

Therefore \eqref{eq:induction step 1c term 1}, \eqref{eq:induction step 1c term 2}, \eqref{eq:eq:induction step 1c term 2 bound}, and \eqref{eq:induction step 1c term 1 third bound} together lead to 
\begin{align} \label{eq:davis khan in induction term 1 third bound}
   & \|  (  \hat G_{1:k}\otimes I_n  )  ^T  ( c ^{*(k+1)} -\hat c^{(k+1) }) V_{k+1}^*  \|_{\op} \leqslant  \frac{\sqrt{ C_3 }}{4 }  \sqrt{\frac{n\log(dn)\|p^*\|_\infty}{N} } .
\end{align}
So by \eqref{eq:davis khan in induction term 1 third bound} and \eqref{eq:lower bound of singular values after projection},
\begin{align} \label{eq:davis khan in induction term 1}
   \eqref{eq:davis khan in induction} \leqslant  \frac{\sqrt{ C_3 }}{2 }  \frac{\sqrt{n\log(dn)\|p^*\|_\infty}}{\sqrt{N} \sigma_2( c^{*(k+1)} )}    .
\end{align} 
\
\\
\\
{\bf Step 1D.} Observe that 
\begin{align} \label{eq:davis khan in induction 2 bound}
    \eqref{eq:davis khan in induction 2} 
    \leqslant \frac{ \|  c ^{*(k+1)} -\hat c^{(k+1) }      \|_{\op}^2  }{ \sigma_2^2 ( ( \hat G_{1:k } \otimes I_n  )  ^T c^{* (k+1 ) } )^2  }   \leqslant \frac{ 4C_2^2 d^K n^{2K}}{ N \sigma^2_2(c^{*(k+1)})  }\leqslant  \frac{\sqrt{ C_3 }}{2 }  \frac{\sqrt{n\log(dn)\|p^*\|_\infty}}{\sqrt{N} \sigma_2( c^{*(k+1)} )}  ,
\end{align}
where second inequality follows from \eqref{eq:s cluster frobenius norm} and \eqref{eq:lower bound of singular values after projection}, and the last inequality hold for sufficiently large $N$ in \eqref{eq:snr}. Consequently, 
\eqref{eq:davis khan in induction}, \eqref{eq:davis khan in induction 2}, \eqref{eq:davis khan in induction term 1}, and \eqref{eq:davis khan in induction 2 bound} together lead to 
\begin{align} \label{eq:induction subspace deviation bound}
    \| \tilde G_{k+1}\tilde G_{k+1}^T - \hat G_{k+1}\hat G_{k+1} \|_{\op} \leqslant  \sqrt{ C_3 }  \frac{\sqrt{n\log(dn)\|p^*\|_\infty}}{\sqrt{N} \sigma_2( c^{*(k+1)} )} .
\end{align}
  \
 \\
{\bf Step 1E.} We complete the induction in this step. Note that 
 \begin{align} \nonumber 
   \eqref{eq:variance induction term 1}
    \leqslant  & \|  \hat G_{k+1} \hat G_{k+1} ^T  - \tilde G _{k+1}\tilde G^T_{k+1}  \|_{\op} ^2    \| c^{* }\|_F ^2 \leqslant  C_3   \frac{n\log(dn)\|p^*\|_\infty}{N \sigma^2_2( c^{*(k+1)} )}    \|c^*\|_F^2  ,
\end{align} 
 where the first inequality follows from   \eqref{eq:induction step 1 term 1 bound 1}, and the last inequality follows from \eqref{eq:induction subspace deviation bound}. \eqref{eq:variance induction term 1} and \eqref{eq:variance induction term 2} together lead to 
 \begin{align*} 
    \|    \hat G_{1:k+1} \hat G_{1:k+1}^T   c^{* (k+1)}  -c^{* (k+1)}\|_F^2 
    \leqslant    C_3  \frac{n\log(dn)\|p^*\|_\infty \|c^{*} \|_F^2}{N} \sum_{j=1}^{k+1}\frac{1}{\sigma_2^2(c^{*(j)})}  .
\end{align*} 
  This completes the induction. Note that by induction
  \begin{align}\label{eq:main variance induction general d minus 1}
\eqref{eq:main variance term 2} =  \|    \hat G_{1:d - 1} \hat G_{1:d -1}^T   c^{* ( d- 1)}  -c^{* (d- 1)}\|_F^2 \leqslant  C_3  \frac{n\log(dn)\|p^*\|_\infty \|c^{*} \|_F^2}{N} \sum_{j=1}^{d-1}\frac{1}{\sigma_2^2(c^{*(j)})} .
\end{align}
\
\\
{\bf Step 2.} We have 
\begin{align} \nonumber 
    \eqref{eq:main variance term 1} =&  \| \hat G_{1:d-1 } \hat G_{1:d -1}^T  ( \hat c^{(d-1)}  -    c^{* (d-1)}  ) \|_F^2 
    \\  \label{eq:coefficient tensor step 2 term 1}
    \leqslant  & 2 \| ( \hat G_{1:d-1 } \hat G_{1:d -1}^T  - U^*_{d-1} U^{*T}_{d-1} )    ( \hat c^{(d-1)}  -    c^{* (d-1)}  )\|_F^2  
    \\
    \label{eq:coefficient tensor step 2 term 2}
    +   &  2 \| U^*_{d-1} U^{*T}_{d-1}     ( \hat c^{(d-1)}  -    c^{* (d-1)}  ) \| _F^2 .
\end{align}
Note that 
 \begin{align} \nonumber
     \eqref{eq:coefficient tensor step 2 term 1} \leqslant & 2 \| \hat G_{1:d-1 } \hat G_{1:d -1}^T  - U^*_{d-1} U^{*T}_{d-1} \| _{\op }^2  \| \hat c^{(d-1)}  -    c^{* (d-1)} \| _F ^2 
    \\ \nonumber
    =  &2 \| \hat G_{1:d-1 } \hat G_{1:d -1}^T  - U^*_{d-1} U^{*T}_{d-1} \| _{\op }^2  \| \hat c   -    c^{*  } \| _F ^2 \leqslant      C_5  \frac{nB^{*3}\|p^*\|_\infty}{ b^{*4}\sigma^2 _2( c^{*(d-1)})} \frac{d\log(dn)}{N}  \| \hat c   -    c^{*  } \| _F ^2
    \\ \label{eq:eq:coefficient tensor step 2 term 2 bound 2}
    \leqslant  &  C_5  \frac{nB^{*3}\|p^*\|_\infty}{ b^{*4}\sigma^2 _2( c^{*(d-1)})} \frac{d\log(dn)}{N} 
        C_2 ^2  \bigg(  \frac{     d^K n^{2K}    } { N}      \bigg) \leqslant  C_6  \frac{n\log(dn)\|p^*\|_\infty }{N} ,
\end{align}
where $C_5,C_6$ are universal constants. The second inequality follows from \Cref{lemma:subspace consistency give tensor consistency} given \eqref{eq:main variance induction general d minus 1}, the third inequality follows from 
\eqref{eq:s cluster frobenius norm}, and the fourth inequality follows from \eqref{eq:snr}. In addition,
\begin{align} \label{eq:coefficient tensor step 2 term 2 bound}
    \eqref{eq:coefficient tensor step 2 term 2} = 2 \| (U^*_{d-1} U^{*T}_{d-1}    \otimes I_n)  ( \hat c^{(d)}  -    c^{* (d )}  ) \| _F^2 \leqslant 2C_1^2 \frac{n\log(dn)\|p^*\|_\infty}{N} .
\end{align}
where the inequality follows from \eqref{eq:sup projected norm}. 
 \\
\\
 Therefore
\begin{align} \label{eq:main variance induction general d minus 1 term 2}
     \eqref{eq:main variance term 1} \leqslant C_7    \frac{n\log(dn)\|p^*\|_\infty}{N}   \leqslant 2C_7 \frac{n\log(dn)\|p^*\|_\infty\|c^*\|^2_F }{\sigma^2 _2( c^{*(d-1)})N}, 
\end{align}
for some universal constant $C_7$, where the first inequality follows from  \eqref{eq:eq:coefficient tensor step 2 term 2 bound 2}  and \eqref{eq:coefficient tensor step 2 term 2 bound}, and the second inequality follows from the fact that 
 $\sigma_2   ( c^{*(d-1}) \leqslant \|c^{*(d-1) }\|_F   =  \|c^{*  }\|_F$. At this point, we combine 
\eqref{eq:main variance induction general d minus 1} and
\eqref{eq:main variance induction general d minus 1 term 2} and obtain the final estimation for the relative error:
\begin{displaymath}
    \frac{\|c^* - \hat G_{1:d }  \|_F}{\|c^{*} \|_F} = {\bigO}_{\p} \left( \sqrt{\frac{n\log(dn)\|p^*\|_\infty }{N} \sum_{j=1}^{d-1}\frac{1}{\sigma_2^2(c^{*(j)})} }\right).
\end{displaymath}
By further applying $\|p^*\|_\infty \leqslant 1+ B^*d$ from \Cref{assumption:additive model}, we arrive at the results in \Cref{thm: main theorem}.
\end{proof}

\subsection{Proof of \Cref{thm: main corollary}}
\label{app:main proof-2}
\begin{proof}
 Let  $ \bar{p}^*  : \mathbb R^d \to \mathbb R$  be the orthogonal projection of $p^*$:
 $$  \bar{p}^*   (x_1,\ldots, x_d) = \sum_{\mi_1= 1 }^n \cdots  \sum_{\mi_d = 1 }^n    c^*(\mi_1,\ldots, \mi_d) \phi_{\mi_1} (x_1) \cdots \phi_{\mi_d} (x_d) .$$
 Then 
    \begin{align*}
         \| p^* -  \tilde{p} \|_{\lt([0,1]^d) } ^2  \leqslant  &  \   \| p^* -  \bar{p}^*\|_{\lt([0,1]^d) }^2 + \|  \bar{p}^* - \tilde{p} \|_{\lt([0,1]^d) }  ^2  
         \\
         = & \  \| p^* - \bar{p}^* \|_{\lt([0,1]^d) }^2 + \| c^*  - 
 \hat{G}_{1:d}  \|_{\lt([0,1]^d) } ^2   
         \\
         \leqslant  & \   C_8 B^{*2} d  n^{-2\alphaspace } +C_9 \frac{d n B^*  \log(dn)  }{N}\sum_{j=1}^{d-1} \frac{ \|c^*\|_F^2  }{\sigma^2_2(c^{*(j)})  }
    \end{align*}
where $C_8$ and $C_9$ are some universal constants. $\hat{G}_{1:d}$ is the output tensor train from TT-SVD (\Cref{algorithm: TT-SVD}). The last inequality follows from the error estimation for the approximation error $ \| p^* - \bar{p}^* \|_{\lt }$ by \Cref{lemma:bias} in the appendix and $\| c^*  - 
 \hat{G}_{1:d}  \|_{\lt }$ by \Cref{thm: main theorem}. Note that by \Cref{corollary:frobenius norm of coefficient tensor},
 $ \| c^*\|_F ^ 2 =\bigO(B^*d) $.
 In addition, by \eqref{eq:normalize the marginal},
 $  \|p^*\|_{\lt([0,1]^d) } ^2 = 1+ \sum_{j=1}^d \| q_j^*\|  _{\lt([0,1]) }^2.$ 
 Therefore,
 $$ 1+b^{*2}d \leqslant \| p^*\| ^2 _{\lt([0,1]^d) } \leqslant 1+ B^{*2}d. $$ by \Cref{assumption:additive model}. Consequently, 
        \begin{align*}
            \frac{ \| p^* - \tilde{p}\|_{\lt([0,1]^d) } ^2  }{\|p^*\|_{\lt([0,1]^d) } ^2 } \leqslant  \  &  C_{10} \left(\frac{B^*}{b^*}\right)^2 \bigg(  n^{-2\alphaspace} + \frac{ dn  \log(dn)  }{N}\sum_{j=1}^{d-1} \frac{ 1 }{\sigma^2_2(c^{*(j)})  }    \bigg). 
        \\
        \leqslant  \  & C_{11} \bigg( \left(\frac{B^*}{b^*}\right)^2 n^{-2\alphaspace} + \frac{B^{*4}}{b^{*6}}\frac{d n  \log^2(dn)  }{N}      \bigg),
        \end{align*} 
        where the last equality follows from \Cref{lemma:lower bounds for singular spectral values} and the assumption that $B^* \geqslant 1$, and the third inequality follows from the fact that 
 $ \sum_{j=1}^{d} \frac{   1 }{ j(d-j)   } =\bigO( \frac{\log(d)}{d})$.

\end{proof}

\section{Useful results for the proof of main theorems}

\subsection{TT-SVD}

\begin{lemma} \label{eq:orthogonality of tensor core composition}
    Let $\hat G_{1:d }$ be the output tensor cores from Algorithm \ref{algorithm: TT-SVD} with target rank $r$. Then 
    $\hat G_{1:k } \in \mathbb O^{n^k\times r}  $ for $k\in \{1,\cdots, d-1 \}$,
\end{lemma}

\begin{proof}
    We proceed by induction. 
    Since  $ \hat G_1 = \svd( \hat c^{(1)}, r) \in \mathbb R^{n\times r},$ it follows that   $\hat G_1 \in \mathbb O^{n \times r} $.  Assume 
    $\hat G_{1:k-1 } \in \mathbb O^{n^{k-1}\times r}  $,   
    by \Cref{lemma:composition of orthogonal matrices} we have
    $$ \hat G_{1:k }=  \big(\hat G_{1:k-1} \otimes I_n \big) \hat G_{k }  \in \mathbb R^{n^{k }\times r} \in \mathbb O^{n^k\times r}. $$
\end{proof}

\begin{lemma} \label{lemma:tensor core orthogonal decomposition}
Let $\hat G_{1:d }$ be the output tensor cores from Algorithm \ref{algorithm: TT-SVD} with target rank $r$. For $k\in \{1,\cdots, d-1 \}$, it holds that 
    \begin{align*}
        &\|\hat G_{1:k}\hat G_{1:k} ^T c^{*(k)} -c^{*(k)}  \|_F^2 
        \\
        = &\|\hat G_{1:k}\hat G_{1:k} ^T c^{*(k)}   -( \hat G_{1:k-1}\otimes I_n )  ( \hat G_{1:k-1} \otimes I_n )  ^T c^{* (k ) }  \|_F^2  +\|\hat G_{1:k-1}\hat G_{1:k-1} ^T c^{*(k-1)} -c^{*(k-1)}  \|_F ^2 .
    \end{align*}     
\end{lemma}
\begin{proof}
    Observe  that 
    \begin{align}\nonumber 
        &\|\hat G_{1:k}\hat G_{1:k} ^T c^{*(k)} -c^{*(k)}  \|_F^2 
        \\ \nonumber 
        = & \|\hat G_{1:k}\hat G_{1:k} ^T c^{*(k)}   -( \hat G_{1:k-1}\otimes I_n )  ( \hat G_{1:k-1} \otimes I_n )  ^T c^{* (k ) }  \|_F^2 
        \\ \label{eq:orthonormal decomposition 1}
        + &   \| ( \hat G_{1:k-1}\otimes I_n )  ( \hat G_{1:k-1} \otimes I_n )  ^T c^{* (k)  }   -c^{*( k)  }  \|_F ^2  
       \\ \label{eq:orthonormal decomposition 2}
       +& \langle \hat G_{1:k}\hat G_{1:k} ^T c^{*(k)}     , ( \hat G_{1:k-1}\otimes I_n )  ( \hat G_{1:k-1} \otimes I_n )  ^T c^{* (k)  }   -c^{*( k)  }  \rangle.
       \\
         \label{eq:orthonormal decomposition 3}
         - & \langle     ( \hat G_{1:k-1}\otimes I_n )  ( \hat G_{1:k-1} \otimes I_n )  ^T c^{* (k ) }  , ( \hat G_{1:k-1}\otimes I_n )  ( \hat G_{1:k-1} \otimes I_n )  ^T c^{* (k)  }   -c^{*( k)  }  \rangle.
    \end{align}     
  {\bf Step 1.}  Note that 
        \begin{align*}
          (  \hat G_{1:k-1}\otimes I_n )  ( \hat G_{1:k-1} \otimes I_n )  ^T c^{* (k)  }   
          = \{ (\hat G_{1:k-1} \hat G_{1:k-1}^T ) \otimes I_n\}c^{* (k)  } =  \hat G_{1:k-1} \hat G_{1:k-1}^T   c^{* (k-1)  }.
    \end{align*}
    Therefore 
     $ 
          \eqref{eq:orthonormal decomposition 1}  =\|\hat G_{1:k-1} \hat G_{1:k-1}^T   c^{* (k-1)  }   -c^{*( k)  } \|_F^2  .$
          \\
          \\
   {\bf Step 2.}  By \Cref{eq:orthogonality of tensor core composition}, 
   $ \hat G_{1:k-1} \in \mathbb O ^{n^{k-1}\times r}$.
   Then the kronecker product of two orthogonal matrices is still orthogonal,
   $\hat G_{1:k-1}\otimes I_n \in \mathbb O^{n^{k}\times rn}$. The details are in \Cref{lemma:kronecker product orthogonality}.
   It follows that
   \begin{align}  \nonumber 
      & (  \hat G_{1:k-1}\otimes I_n )  ( \hat G_{1:k-1} \otimes I_n )  ^T c^{* (k)  }   -c^{*( k) } 
       =  \big( ( \hat G_{1:k-1}\otimes I_n )  ( \hat G_{1:k-1} \otimes I_n )^T    - I_{n^{k  }} \big)  c^{*( k) } 
       \\ \label{eq:tensor core coposition complement}
       =  &  ( \hat G_{1:k-1}\otimes I_n )_\perp   
       \big[ ( \hat G_{1:k-1} \otimes I_n ) _\perp  \big] ^T c^{*( k) }  .
   \end{align}  
Therefore 
\begin{align*}
\eqref{eq:orthonormal decomposition 3}     =
 \langle     ( \hat G_{1:k-1}\otimes I_n )  ( \hat G_{1:k-1} \otimes I_n )  ^T c^{* (k ) }  ,  ( \hat G_{1:k-1}\otimes I_n )_\perp   
       \big[ ( \hat G_{1:k-1} \otimes I_n ) _\perp  \big] ^T c^{*( k) }  \rangle = 0.
\end{align*}
    {\bf Step 3.} By  \eqref{eq:tensor core coposition complement},  it follows that
\begin{align*}
    \eqref{eq:orthonormal decomposition 2} = &\langle ( \hat G_{1:k-1}\otimes I_n ) \hat G_k \hat G_{1:k} ^T c^{*(k)}     , ( \hat G_{1:k-1}\otimes I_n )_\perp   
       \big[ ( \hat G_{1:k-1} \otimes I_n ) _\perp  \big] ^T c^{*( k) }  \rangle 
    \\
    = & \langle  \big[ ( \hat G_{1:k-1} \otimes I_n ) _\perp  \big] ^T   (\hat G_{1:k-1}\otimes I_n ) \hat G_k \hat G_{1:k} ^T c^{*(k)}     ,    
        \big[ ( \hat G_{1:k-1} \otimes I_n ) _\perp  \big] ^T c^{*( k) }   \rangle = 0.
\end{align*}
This immediately leads to the desired result. 
\end{proof}

\begin{lemma} \label{lemma:tensor core orthogonal decomposition 2} Let $\hat G_{1:d }$ be the output tensor cores from Algorithm \ref{algorithm: TT-SVD} with target rank $r$. Then
   \begin{align*}
     \|c^* - \hat G_{1:d }  \|_F^2  
       =  
       \|  c^{* (d-1)} -\hat G_{1:d-1 } \hat G_{1:d -1}^T   c ^{* (d-1)}    \|_F^2 +
    \|   \hat G_{1:d-1 } \hat G_{1:d -1}^T   c^{* (d-1)}  - \hat G_{1:d-1 } \hat G_{1:d -1}^T \hat c^{(d-1)}  \|_F^2  
         . 
\end{align*} 
\end{lemma}
\begin{proof} 
Since $\hat G_{1:d-1 } \hat G_{1:d -1}^T \hat c ^{(d-1)} $ is a reshape matrix of $\hat G_{1:d }     $, we have 
$$ \|c^* - \hat G_{1:d }  \|_F^2 =  \|c^{* (d-1)}  - \hat G_{1:d-1 } \hat G_{1:d -1}^T \hat c ^{(d-1)}     \|_F^2 .
 $$
It suffices to show that 
\begin{equation} \label{eq:tensor core orthogonal decomposition 2 cross term}
    \langle\hat G_{1:d-1 } \hat G_{1:d -1}^T   c ^{* (d-1)} -  c^{* (d-1)}  , \  \hat G_{1:d-1 } \hat G_{1:d -1}^T \hat c^{(d-1)}  -\hat G_{1:d-1 } \hat G_{1:d -1}^T   c^{* (d-1)}      \rangle  = 0.
\end{equation}  
By \Cref{eq:orthogonality of tensor core composition},  $ \hat G_{1:d-1} \in \mathbb O ^{n^{d-1}\times r}$. So 
\begin{align*}    \hat G_{1:d-1 } \hat G_{1:d -1}^T   c ^{* (d-1)} -  c^{* (d-1)} =  (\hat G_{1:d-1 })_\perp  (\hat G_{1:d -1}^T)_\perp c^{* (d-1)}.
\end{align*}  
Therefore $(\hat G_{1:d-1 }) ^T_\perp  \hat G_{1:d-1 }=0$ and so
\begin{align*}
    \eqref{eq:tensor core orthogonal decomposition 2 cross term} =
    &\langle  (\hat G_{1:d-1 })_\perp  (\hat G_{1:d -1})_\perp^T c^{* (d-1)} , \ \hat G_{1:d-1 } \hat G_{1:d -1}^T ( \hat c^{(d-1)}  -   c^{* (d-1)} )     \rangle  
    \\
    =& \langle   (\hat G_{1:d -1})_\perp^T c^{* (d-1)} , \ (\hat G_{1:d-1 }) ^T_\perp  \hat G_{1:d-1 } \hat G_{1:d -1}^T ( \hat c^{(d-1)}  -   c^{* (d-1)} )     \rangle  =0.
\end{align*}
\end{proof}

 \begin{lemma} \label{lemma:subspace consistency give tensor consistency} Suppose all the conditions in \Cref{thm: main theorem} holds.  Let $\hat G_{1:d }$ be the output tensor cores from Algorithm \ref{algorithm: TT-SVD} with target rank $r=2$. Suppose 
       that 
 \begin{equation}   \label{eq:general singular value lower bound}
   \| \hat G_{1:k } \hat G_{1:k}^T  c ^{*(k)}  -   c^{* (k)}  \|_F^2  
    \leqslant    C  \frac{n\log(dn)\|p^*\|_\infty \|c^{*} \|_F^2}{N} \sum_{j=1}^k\frac{1}{\sigma_2^2(c^{*(j)})} 
\end{equation}  
where $C$ is a universal constant. Let  $U^*_k$ be   defined as in \Cref{lemma:additive row and column spaces}. Note that   by \Cref{lemma:additive row and column spaces},
 $U_k^*U_k^{*T}   \in \mathbb R^{n^k \times n^k} $  is the projection matrix onto the space $\col(c^{* (k)})$.
Then 
 \begin{equation}\label{{eq:subspace consistency give tensor consistency}}
       \| \hat G_{1:k } \hat G_{1:k}^T  - U^*_k U^{*T }_k\|_{\op}    
        \leqslant C'  \frac{\sqrt{nB^{*3}\|p^*\|_\infty}}{ b^{*2}\sigma _2( c^{*(k)})} \sqrt{\frac{d\log(dn)}{N}},
 \end{equation}  
 where $C'$ a universal constant.
 \end{lemma}
\begin{proof}
Note that by \Cref{lemma:lower bounds for singular spectral values}, for $k \in \{1,\cdots, d-1 \}$,
$
    \sigma_2^2(c^{*(k)}) \geqslant \frac{c_\sigma^2 k(d-k)}{d}.
$ In addition, by \Cref{corollary:frobenius norm of coefficient tensor}, $\| c^*\|_F^2\leqslant C_c d $ where $C_c = B^* + 2$. Therefore by \eqref{eq:general singular value lower bound}
\begin{equation}   \label{eq:induction upper bound for frobenius norm}
   \| \hat G_{1:k } \hat G_{1:k}^T  c ^{*(k)}  -   c^{* (k)}  \|_F^2  
    \leqslant    C  \frac{n\log(dn)\|p^*\|_\infty \|c^{*} \|_F^2}{N} \frac{d}{c^2_\sigma} \sum_{j=1}^k\frac{1}{j(d-j)} \leqslant C'' \frac{C_cn\|p^*\|_\infty}{c^2_\sigma} \frac{d\log(dn)}{N}
\end{equation} 
 where we have used the fact that 
 $ \sum_{j=1}^{d} \frac{   1 }{ j(d-j)   } =\bigO( \frac{\log(d)}{d})$ and $C''$ is some universal constant.

 By \Cref{lemma:left projection preserves matrix structure},  $ \rank (  \hat G_{1:k } \hat G_{1:k}^T  c ^{*(k)} ) =2$,
  and   
 $\hat G_{1:k } \hat G_{1:k}^T$ is the projection matrix on $\col (\hat G_{1:k } \hat G_{1:k}^T  c ^{*(k)} ) $. 
 Consequently,
\begin{equation} \label{eq:induction step 1c term 1 second bound}
 \begin{split}
     \| \hat G_{1:k } \hat G_{1:k}^T  - U^*_k U^{*T }_k\|_{\op}    \leqslant    \frac{ 2\| \hat G_{1:k } \hat G_{1:k}^T  c ^{*(k)}  -   c^{* (k)}  \|_{\op}     }{ \sigma _2( c^{*(k)})  -\| \hat G_{1:k } \hat G_{1:k}^T  c ^{*(k)}  -   c^{* (k)}  \|_{\op}   } 
        \leqslant \ &   \frac{ 2\| \hat G_{1:k } \hat G_{1:k}^T  c ^{*(k)}  -   c^{* (k)}  \|_{\op}     }{ \frac{1}{2}\sigma _2( c^{*(k)}) }  \\
         \leqslant \ & 4  \frac{\sqrt{C''C_cn\|p^*\|_\infty}}{c_\sigma \sigma _2( c^{*(k)})} \sqrt{\frac{d\log(dn)}{N}}
        \end{split}
 \end{equation}
       where the first inequality follows from   \Cref{lemma:Wedin} and $\rank(c^{*(k)})=2$, the second inequality follows from the assumption on sufficiently large sample size $N$ as stated in \eqref{eq:snr}, and the last inequality follows from \eqref{eq:induction upper bound for frobenius norm}. Furthermore, by the expressions of $C_c$ in \Cref{corollary:frobenius norm of coefficient tensor} and $c_\sigma$ in \Cref{lemma:lower bounds for singular spectral values}, we have $\frac{\sqrt{C_c}}{c_\sigma} \sim \sqrt{\frac{B^{*3}}{b^{*4}}}$ can be bounded a universal constant, leading to the final estimation \eqref{{eq:subspace consistency give tensor consistency}}. 
\end{proof}

\subsection{Ground truth density function}



\begin{lemma}
    \label{lemma:bias}  Under \Cref{assumption:additive model}, let 
 $$  \bar{p}^*   (x_1,\cdots, x_d) = \sum_{\mi_1= 1 }^n \cdots  \sum_{\mi_d = 1 }^n    c^*(\mi_1,\cdots, \mi_d) \phi_{\mi_1} (x_1) \cdots \phi_{\mi_d} (x_d) $$
where $c^* \in \mathbb R^{ n ^d}$ is the coefficient tensor formulated as \eqref{eq:additive model for coefficient tensor}. There exists an absolute constant $C$ such that 
 $$ \| p^* - \bar{p}^* \|_{\lt([0,1]^d)} \leqslant C n^{-  \alphaspace} \sum_{j=1}^d \| q^*_j \| _{H^{  \alphaspace} ([0,1])} .$$
\end{lemma}
\begin{proof}
  Direct calculations show that 
  $$ \bar{p}^*   (x_1,\cdots, x_d) = 1 + \sum_{l_1=1}^n c^*_{1}(l_1) \phi_{l_1} (x_1)  + \cdots + \sum_{l_d= 1}^n c^*_{d}(l_d) \phi_{l_d} (x_d)   .$$
 where $c^*_j \in \mathbb{R}^{n}$ is defined in \eqref{eq:explicit  expression of coefficients}. Therefore 
  $$  \| p^* - \bar{p}^* \|_{\lt([0,1]^d)} \leqslant \sum_{j=1}^d 
  \bigg \| \sum_{l_j =1}^n c^*_{j}(l_j) \phi_{l_j}   - q_j^*   \bigg \|_{\lt ([0,1])}. $$ 
Denote $\mathcal S_n$ the class of polynomial up degree $n$ on $[0,1]$, and  
  $\mathcal P_n$   the projection operator from $\lt([0,1])$ to $\mathcal S_n$  respect to $\lt$ norm.    Since Legendre polynomials  are orthogonal, \eqref{eq:explicit  expression of coefficients} and  \eqref{eq:normalize the marginal}
  \begin{align}\label{eq:identify projected marginal} \sum_{l_j =1}^n c^*_{j}(l_j) \phi_{l_j} =\mathcal P_n(q_j^*).    
  \end{align}
  By the  Legendre polynomial     approximation theory  shown in the appendix of  \cite{khoo2024nonparametric}, 
  it follows that 
   \begin{align}\label{eq:approximatiopn error for marginals}  \bigg \| \sum_{l_j =1}^n c^*_{j}(l_j) \phi_{l_j}   - q_j^*   \bigg \|_{\lt ([0,1])} \leqslant C\| q_j^*\|_{H^{  \alphaspace}([0,1])} n^{-  \alphaspace}
   \end{align}
  for some absolute constant $C$ which does not depend on $c^*_j$ or $q^*_j$.  The desired result follows immediately.
\end{proof}

\begin{corollary} \label{corollary:frobenius norm of coefficient tensor}
Under \Cref{assumption:additive model}, for sufficiently large $n$, there exists a constant $C$ such that 
 $$ \| c^*\|_F^2 \leqslant  C_c d , \quad \text{where} \quad C_c =  B^*+2$$
 where $c^* \in \mathbb R^{ n ^d}$ is the coefficient tensor formulated as \eqref{eq:additive model for coefficient tensor}.
\end{corollary}
\begin{proof}
Observe that 
 \begin{align*}
     \| c^*\|_F^2  = &\| \bvec{1}^d +    c_1^* \otimes \bvec{1}^{d-1} + \bvec{1}  \otimes c_2^* \otimes \bvec{1}^{d-2} + \cdots + \bvec{1}^{d-1}\otimes c_d^* \|_F^2 
     \\
     = & \| \bvec{1}^d  \| _F^2 +  \| c_1^* \otimes \bvec{1}^{d-1} 
 \|_F^2 + \cdots  +   \|\bvec{1}^{d-1}\otimes c_d^*\|_F^2
 =   1 + \sum_{j=1}^d \| c_j^*\|^2.
 \end{align*} 
where the first equality follows from 
 \eqref{eq:additive model for coefficient tensor} , and the second equality follows from orthogonality.  In addition,  note that 
 $$  \| c_j^*\|= \big\|\sum_{l_j =1}^n c^*_{j}(l_j) \phi_{l_j} \big \| _{\lt ([0,1])} \leqslant  \| q_j^*\| _{\lt ([0,1])}  + C_1\| q_j^*\|_{H^{  \alphaspace}([0,1])} n^{-  \alphaspace} 
 \leqslant B^*+1,  $$
 where the first inequality follows from 
 \eqref{eq:approximatiopn error for marginals}, and the last inequality follows from \Cref{assumption:additive model} and $n$ being sufficiently large. The desired result follows immediately.
\end{proof}

\subsection{Row and column spaces of the coefficient tensor}
\begin{lemma}
    Let $k\in \{1,\cdots, d-1 \}$. Let $c^{* (k)} \in \mathbb R^{n^k}\times \mathbb R^{n^{d-k}} $ be the $k$-th unfolding of the tensor coefficient $c^*$ as in \eqref{eq:additive model for coefficient tensor} for $k=1,\cdots,d-1$. 
The row and the column space of $c^{* (k)} $ are given by   
    \begin{align}\label{eq:column space basis}
    \col ( c^{* (k)}) = &\Span \{ \bvec{1}^{k} , \    c_1^* \otimes \bvec{1}^{ k-1} + \bvec{1} \otimes  c_2^* \otimes \bvec{1}^{ k-2} + \cdots +\bvec{1}^{k-1} \otimes c_k^* \}, 
    \\ \label{eq:row space basis}
      \row ( c^{* (k)}) =&  \Span \{ \bvec{1}^{d-k} , \  c_{k+1} ^* \otimes \bvec{1}^{d-k-1} + \bvec{1} \otimes  c_{k+2} ^* \otimes \bvec{1}^{d-k-2} + \cdots +\bvec{1}^{d-k-1} \otimes c_d^* \}.   \end{align}
Consequently, $\rank(c^{* (k)}) \leqslant 2$. Furthermore, if $q^*_j \not = 0  $ for all $j\in \{ 1,\cdots, d\}$,  
      $\rank(c^{* (k)}) =2$.
\end{lemma}
\begin{proof}
    It is direct to show that 
    \begin{equation}\label{eq:reshape row column spaces}
    \begin{split}
   c^{*(k)}  = & ( \bvec{1}^k )  \otimes   \{   \bvec{1}^{d-k} \}   + (   c_1^*  \otimes  \bvec{1}^{k-1} ) \otimes  \{ \bvec{1}^{d-k} \}  + 
\cdots + ( \bvec{1}^{k-1} \otimes c_k^*)  \otimes \{ \bvec{1}^{d-k}  \}
\\ 
+ & (\bvec{1}^{k })  \otimes  \{ c_{k+1} ^* \otimes \bvec{1}^{d-k-1}\} + \cdots  + (\bvec{1}^{k })  \otimes  \{ \bvec{1}^{d-k-1} \otimes c_d^*\},
\end{split}
\end{equation}
where   terms with $( \  )$ are  viewed as column vectors in $\mathbb R^{n^k}$, and  terms with $\{ \  \} $ are  viewed as   row vectors in $\mathbb R^{n^{d-k}}$.
To  show  \eqref{eq:column space basis} holds for the column space, it suffices to rewrite  
\begin{align*}
    c^{*(k)} = & ( \bvec{1}^k )  \otimes   \big \{   \bvec{1}^{d-k} + 
   c_{k+1} ^* \otimes \bvec{1}^{d-k-1} + \bvec{1} \otimes  c_{k+2} ^* \otimes \bvec{1}^{d-k-2} + \cdots +\bvec{1}^{d-k-1} \otimes c_d^*   
   \big\}  
   \\ +&  \big(c_1^* \otimes \bvec{1}^{ k-1} + \bvec{1} \otimes  c_2^* \otimes \bvec{1}^{ k-2} + \cdots +\bvec{1}^{k-1} \otimes c_k^* \big) \otimes \{  \bvec{1} ^{d-k}\}.
\end{align*}
Therefore  \eqref{eq:column space basis} holds. The justification of   \eqref{eq:row space basis}  is similar and   omitted for brevity.
\end{proof}

\begin{lemma} \label{lemma:additive row and column spaces}
Let $c^*$ be coefficient tensor as in \eqref{eq:additive model for coefficient tensor} and $c^*_j$ be as in \eqref{eq:explicit  expression of coefficients}. Denote   $ \bvec{e}_1 = [1,0]^T  $, $ \bvec{e}_2 = [0,1]^T  $ and let the matrices
$U_k^* \in \mathbb R^{n^k \times 2    }  $ and $V_k^*\in \mathbb R^{n^{d-k} \times  2 }     $ be such that 
\begin{align} \label{eq:additive column space}
    U_k^* =&  \bvec{1}^k \otimes \bvec{e}_1  +   \frac{  c_1^* \otimes \bvec{1}^{ k-1} + \bvec{1} \otimes  c_2^* \otimes \bvec{1}^{ k-2} + \cdots +\bvec{1}^{k-1} \otimes c_k^*  }{\sqrt{  \| c_1^*\|^2 + \cdots +  \| c_k^*\|^2}} \otimes  \bvec{e}_2,
    \\ \label{eq:additive row space}
    V_k^* = &   \bvec{1}^{d-k} \otimes \bvec{e}_1  +    \frac{     c_{k+1} ^* \otimes \bvec{1}^{d-k-1} + \bvec{1} \otimes  c_{k+2} ^* \otimes \bvec{1}^{d-k-2} + \cdots +\bvec{1}^{d-k-1} \otimes c_d^*  }{\sqrt{  \| c_{k+1}^*\|^2 + \cdots +  \| c_d^*\|^2}} \otimes \bvec{e}_2 .
\end{align} 
Then 
 $U_k^*U_k^{*T}   \in \mathbb R^{n^k \times n^k} $  is the projection matrix onto the space $\col(c^{* (k)})$, and $V_k^*V_k^{*T}  \in \mathbb R^{n^{d-k} \times n^{d-k}} $  is the projection matrix onto  the space $\row(c^{* (k)})$.
\end{lemma}
\begin{proof}
    Note that the set of vectors  $\{ \bvec{1}^k, c_1^* \otimes \bvec{1}^{ k-1}  ,  \bvec{1} \otimes  c_2^* \otimes \bvec{1}^{ k-2}  ,  \bvec{1}^{k-1} \otimes c_k^*\}$ are    orthogonal. Direct calculations show that 
    $$ U_k^{*T}  U_k^* =\bvec{e}_1 \otimes \bvec{e}_1 + 
    \bvec{e}_2  \otimes \bvec{e}_2 =I_2 \in \mathbb R^{2\times 2}.$$
    Therefore $U_k^* \in \mathbb O^{n^k \times 2    }$. Since $\col (U_k^*  ) =\col ( c^{* (k)})  $, $U_k^*U_k^{*T} \in \mathbb R^{n^k \times n^k} $  is the projection matrix onto the space $\col(c^{* (k)})$.
    The justification of  $V_k^* V_k^{*T} $   is similar and   omitted for brevity.
\end{proof}
For simplicity, in \eqref{eq:additive column space} and \eqref{eq:additive row space} we introduce the notations 
    \begin{align}\label{eq:reshape matrix colum basis tensor}
        &u^*_k =    \frac{  c_1^* \otimes \bvec{1}^{ k-1} + \bvec{1} \otimes  c_2^* \otimes \bvec{1}^{ k-2} + \cdots +\bvec{1}^{k-1} \otimes c_k^*  }{\sqrt{  \| c_1^*\|^2 + \cdots +  \| c_k^*\|^2}}  \ \in \mathbb R^{n^k} ,
        \\
        \label{eq:reshape matrix row basis tensor} 
       & v^*_k =  \frac{     c_{k+1} ^* \otimes \bvec{1}^{d-k-1} + \bvec{1} \otimes  c_{k+2} ^* \otimes \bvec{1}^{d-k-2} + \cdots +\bvec{1}^{d-k-1} \otimes c_d^*  }{\sqrt{  \| c_{k+1}^*\|^2 + \cdots +  \| c_d^*\|^2}}  \ \in \mathbb R^{n^{d-k}}  .
       \end{align}

\begin{lemma} \label{lemma:norms of projection basis functions}
 Let $u^*_k$ and $v^*_k$ be as in \eqref{eq:reshape matrix colum basis tensor} and \eqref{eq:reshape matrix row basis tensor} respectively. Denote
    \begin{align}
 \label{eq:reshape matrix colum basis function}
      &\Phi^* _k(x_1,\cdots, x_k) = \frac{ \sum_{\mi_1=1}^n   c_{1}  ^* (\mi_1)  \phi_{\mi_1}(x_1) +  \cdots +\sum_{\mi_k=1}^n   c_{k} ^*(\mi_k)  \phi_{\mi_k}(x_k) }{\sqrt{  \| c_1^*\|^2 + \cdots +  \| c_k^*\|^2} }  ,
     \\ \label{eq:reshape matrix row basis function}
     &\Psi^* _k(x_{k+1},\cdots, x_d) =   \frac{ \sum_{\mi_{k+1}=1}^n   c_{k+1}  ^*(\mi_{k+1})   \phi_{\mi_{k+1}}(x_{k+1}) +  \cdots +\sum_{\mi_d=1}^n   c^* _{d}(\mi_d) \phi_{\mi_d}(x_d)}{\sqrt{  \| c_{k+1}^*\|^2 + \cdots +  \| c_d^*\|^2}}.
 \end{align}
 Then $u^*_k $ is the coefficient tensor of $ \Phi^* _k$, and $v^*_k $ is the coefficient tensor of $ \Psi^* _k$.
Then $$\| \Phi_k^* \|_{\lt([0,1]^{k})} = \| u^*_k\|_F =1 \quad \text{and} \quad  
 \| \Psi_k^* \|_{\lt([0,1]^{d-k})} = \| v^*_k\|_F =1.$$  
In addition, suppose    \Cref{assumption:additive model} holds.
  Then for sufficient large $n$ such that 
$$ \| \Phi_k^* \|_\infty  \leqslant   C'\sqrt{k} \quad \text{and} \quad \| \Psi_k^* \|_\infty \leqslant  C'\sqrt{d-k} \quad \text{where} \quad C' = \frac{4B^*}{b^*} . $$
\end{lemma}
\begin{proof} To see that  $u^*_{k}$ is the coefficient tensor of $\Phi^* _k(x_1,\cdots, x_k) $, it suffices to observe that the $(l_1, \cdots, l_k)$-coordinate of  $u^*_{k}$  is 
$  \langle\Phi^* _k , \phi^1_{l_1} \cdots\phi^k_{l_k}  \rangle  .  $ 
Since the coefficient tensor is distance preserving, it holds that $$ \| \Phi_k^* \|_{\lt([0,1]^{k})} = \| u^*_k\|_F =1 .$$
\
\\
Next, we bound $ \| \Phi^*_k\|_\infty $. 
By \eqref{eq:identify projected marginal} and \Cref{theorem:infty norm approximation of polynomial}, it follows that there exists a constant $C_1$ such that 
$$ \big\| \sum_{\mi_j=1}^n   c_{j }  ^*(\mi_1)   \phi_{\mi_j}   \big\|_\infty = \| \mathcal P_n (q_j^*)\|_\infty  \leqslant \| q_j^* \|_\infty  +C_1 n^{-  \alphaspace} .$$
Therefore, the numerator of \eqref{eq:reshape matrix colum basis function} is bounded by  
\begin{align} \label{eq:upper bound in L infty bound}
  \bigg\|   \sum_{\mi_1=1}^n   c_{1}  ^*(\mi_1)   \phi_{\mi_1} +  \cdots +\sum_{\mi_k=1}^n   c_{k } ^*(\mi_k)  \phi_{\mi_k}  \bigg\|_\infty  \leqslant \sum_{j=1}^k \|q_j^*\|_\infty +C_1 k n^{-  \alphaspace} \leqslant 2kB^*
\end{align}
   upon using \Cref{assumption:additive model} and sufficiently large $n$. Observe that 
    $$\| c_j^*\| =\big\| \sum_{\mi_j=1}^n   c_{j}  ^*( \mi_j)   \phi_{\mi_j}   \big\|_{\lt([0,1])} = \| \mathcal P_n(q_j^*) \|_{\lt([0,1])} \geqslant \| q_j^*\|_{\lt([0,1])} - C_2 n^{-  \alphaspace} \geqslant b^*- C_2 n^{-  \alphaspace} \geqslant b^*/2  $$
    from some constant $C_2$, where   the second  equality follows from \eqref{eq:identify projected marginal},   the first  inequality follows from \Cref{corollary:legendre}, the second inequality follows from the \Cref{assumption:additive model} that   $\min_{j\in\{1,\cdots, d\}} \| q^*_j\|_{\lt([0,1])} \geqslant b^*    $, and the last inequality holds for sufficiently large  $n$. Therefore, the denominator
    \begin{align}\label{eq:lower bound L infty bound}
    \sqrt {\| c_1^*\|^2 + \cdots +  \| c_k^*\|^2 }  \geqslant  b^*\sqrt{k}/2
    \end{align}
    Note that \eqref{eq:upper bound in L infty bound} and \eqref{eq:lower bound L infty bound} immediately imply $\| \Phi_k^* \|_\infty  \leqslant   4\frac{B^*}{b^*}\sqrt{k} $. The analysis for $\| \Psi_k^* \|_\infty$ is the same and therefore omitted. 
\end{proof}

\subsection{Singular Values of the reshaped coefficient tensor}
\begin{lemma} \label{lemma:lower bounds for singular spectral values}
Let $c^*$  be defined as in \eqref{eq:additive model for coefficient tensor}. Denote 
$$a_k^* =\sum_{i=1}^k \|c_i^*\|^2 \quad \text{and} \quad b_k^*=\sum_{i=k+1}^d \|c_i^*\|^2.$$  Then the $k$-th unfolding of coefficient tensor $c^*$ has the smallest singular value formulated as 
\begin{displaymath}
    \sigma_2^2(c^{*(k)}) = \frac{1+a^*_k+b^*_k-\sqrt{(1+ a^*_k+b^*_k)^2 - 4a^*_kb^*_k}}{2}
\end{displaymath}
$\text{~for~} k = 1,\cdots, d$. Suppose in addition that \Cref{assumption:additive model} holds, for sufficient large $n$, $\sigma_2(c^{*(k)})$ is lower bounded by 
    \begin{equation} \label{eq:lower bound singular values}
        \sigma_2(c^{*(k)}) \geqslant  c_\sigma \sqrt{\frac{  k(d-k)}{d} } \quad \text{where} \quad c_\sigma = \frac{b^{*2}}{4\sqrt{1+4B^{*2}}}.
    \end{equation}    
\end{lemma}

\begin{proof}
Let $k\in\{1,\cdots, d \}$. Denote 
$$f_k^* = c^*_1 \otimes \bvec{1}^{k-1} + \bvec{1} \otimes c^*_2 \otimes \bvec{1}^{k-2} + \cdots  + \bvec{1}^{k-1} \otimes c^*_k,   \text{ and }  g_k^*= c^*_{k+1} \otimes \bvec{1}^{d-k-1} + \cdots + \bvec{1}^{d-k} \otimes c^*_d .$$
Then one can easily check that $a^*_k = \|f^*\|_F^2$ and $b^*_k = \|g^*\|_F^2$ due to orthogonality, and we can write 
\begin{align*}
   c^* =   ( \bvec{1}^k )  \otimes   \big \{   \bvec{1}^{d-k} + 
    g_k^*  
   \big\}  + \big( f_k^* \big) \otimes \{  \bvec{1} ^{d-k} \} .
\end{align*}
where   terms with $( \  )$ are  viewed as column vectors in $\mathbb R^{n^k}$, and  terms with $\{ \  \} $ are  viewed as   row vectors in $\mathbb R^{n^{d-k}}$.
Therefore
$$c^{*(k)} = \begin{bmatrix}
    \bvec{1}^k  &  f_k^*
\end{bmatrix}
\begin{bmatrix}
\bvec{1}^{d-k} + g_k^* \\ \bvec{1}^{d-k}
\end{bmatrix}.
$$
Observe that 
\begin{equation*}
    c^{*(k)}(c^{*(k)})^T = \begin{bmatrix}
    \bvec{1}^k &  f^*_k
\end{bmatrix}
\begin{bmatrix}
1 + \|g^*_k\|_F^2 & 1 \\ 1 &1 
\end{bmatrix}
\begin{bmatrix}
    \bvec{1}^k \\  f^*_k
\end{bmatrix} = 
\begin{bmatrix}
    \bvec{1}^k &  f^*_k/\|f^*_k\|_F
\end{bmatrix}
\begin{bmatrix}
1 + \|g^*_k\|_F^2 & \|f^*_k\|_F \\ \|f^*_k\|_F & \|f^*_k\|_F^2 
\end{bmatrix}
\begin{bmatrix}
    \bvec{1}^k \\  f^*_k/\|f^*_k\|_F
\end{bmatrix}
\end{equation*}
Note that $ \begin{bmatrix} 
    \bvec{1}^k &  f^*_k/\|f^*_k\|_F
\end{bmatrix} $ is an orthogonal matrix in $\mathbb O^{n^k\times 2}$. Therefore the eigenvalues of $c^{*(k)}(c^{*(k)})^T$ coincide with  the eigenvalues of  $\begin{bmatrix}
1 + \|g^* _k\|_F^2 & \|f ^* _k\|_F \\ \|f^* _k\|_F & \|f^*_k\|_F^2 
\end{bmatrix} \in \mathbb R^{2\times 2}$. Direct calculations show that 
$$\sigma_2^2(c^{*(k)}) = \lambda_2( c^{*(k)}(c^{*(k)})^ T   ) =\frac{1+a^*_k+b^*_k-\sqrt{(1+ a^*_k+b^*_k)^2 - 4a^*_kb^*_k}}{2}  .$$
For \eqref{eq:lower bound singular values}, note that by \eqref{eq:approximatiopn error for marginals}, there exists a constant $C_1$
   \begin{align} \bigg \| \sum_{l_j =1}^n c^*_{j}(l_j) \phi_{l_j}   - q_j^*   \bigg \|_{\lt ([0,1])} \leqslant C_1  n^{-  \alphaspace}
   \end{align}
Since $ \big  \| \sum_{l_j =1}^n c^*_{j}(l_j) \phi_{l_j}\big \| _{\lt ([0,1])}  =\|c^*_j\| $, it follows that 
$$    \| q_j^*  \|_{\lt ([0,1])} - C_1n^{-  \alphaspace } \leqslant \|c^*_j\|  \leqslant \| q_j^* \|_{\lt ([0,1])} + C_1 n^{-  \alphaspace }  $$
It follows from   \Cref{assumption:additive model} and 
\eqref{eq:approximatiopn error for marginals} that 
   $$ b^* \leqslant  \| q^*_j\|_{\lt([0,1])} \leqslant B^*.$$
So for sufficiently large $n$, it follows that 
 $ \frac{b^*}{2} \leqslant \|c^*_j\| \leqslant 2B^*.$ 
 Consequently, 
 $$  \frac{b^{*2} k }{4} \leqslant a^*_k \leqslant 4B^{*2} k \quad \text{and} \quad  \frac{b^{*2} (d-k) }{4} \leqslant b^*_k  \leqslant 4B^{*2}  (d-k) .$$
Therefore 
\begin{align*}
    \sigma_2^2(c^{*(k)})     =  & \frac{1+a^*_k+b^*_k-\sqrt{(1+ a^*_k+b^*_k)^2 - 4a^*_kb^*_k}}{2}
    \\
    =& \frac{(1+a^*_k+b^*_k-\sqrt{(1+ a^*_k+b^*_k)^2 - 4a^*_kb^*_k} ) (1+a^*_k+b^*_k+\sqrt{(1+ a^*_k+b^*_k)^2 - 4a^*_kb^*_k )}}{2(1+a^*_k+b^*_k +\sqrt{(1+ a^*_k+b^*_k)^2 - 4a^*_kb^*_k)}} 
    \\ 
    \geqslant &  \frac{ 4 a^*_k b_k^*}{4(1+a^*_k+b^*_k   )} \geqslant \frac{\frac{1}{16}b^{*4} k(d-k)}{1+4B^{*2}d} \geqslant  \frac{\frac{1}{16}b^{*4} k(d-k)}{(1+4B^{*2})d}
\end{align*}

\end{proof}

\subsection{Variance of cluster basis estimators}
\begin{lemma}\label{lemma:maximum projection norm} Suppose    \Cref{assumption:additive model} holds.
 Let $c^*$ be the coefficient tensor of $p^*$, and $\hat c$ be the $s$-cluster 
 coefficient tensor estimator defined in 
 \Cref{definition:cluster coefficient tensor estimator} with $s\geqslant 2$. 
Let $U^*_k$ and $V_k^*$  be defined as in \Cref{lemma:additive row and column spaces}. Then, 
\begin{equation}\label{eq:maximum projection norm}
\max_{k\in\{ 1,\cdots, d\}}\|  ( U_{k-1}^{*T}   \otimes I_n )  (c^{*(k)} -\hat c^{(k)}[s]) V_k ^* \|_F = {\bigO } _\p\left(\sqrt{\frac{n\log(dn)\|p^*\|_\infty}{N} } \right) .
\end{equation}


\end{lemma}

\begin{proof}
    We suppose $s \geqslant 3$ in the proof, where the case of $s=2$  is a simplified version from the same calculations. Given any positive number $t$, we have  
\begin{equation}\label{eq:Bernstein bound for projection term detail} 
    \begin{split}
        &\p\left( \max_{k\in\{ 1,\cdots, d\}}\|  ( U_{k-1}^{*T}   \otimes I_n )  (c^{*(k)} -\hat c^{(k)}[s]) V_k ^* \|^2_F  \geqslant 4n t^2\right)\\
    = \ &     \p\left(\|  ( U_{k-1}^{*T}   \otimes I_n )  (c^{*(k)} -\hat c^{(k)}[s]) V_k ^* \|^2_F  \geqslant 4n t^2 \text{~for least one~} k \in\{ 1,\cdots, d\}  \right)\\
   \leqslant \ &  \sum_{k=1}^d \p\left(\|  ( U_{k-1}^{*T}   \otimes I_n )  (c^{*(k)} -\hat c^{(k)}[s]) V_k ^* \|^2_F  \geqslant 4n t^2  \right) \qquad \qquad  \text{(by union bound)}\\
   = \ & \sum_{k=1}^d \p\left( \eqref{eq:left and right projection} 
 \geqslant 4n t^2  \right)\\
 \leqslant \ & \sum_{k=1}^d \sum_{l=1}^n \p\left( | \langle p^*-\hat p , \Phi^* _{k-1} \phi_\mi^k  \Psi^*_k\rangle| \geqslant t \right) + \p\left( | \langle p^*-\hat p ,  \phi_\mi^k \rangle| \geqslant t \right) \\
 &\hspace{50pt}+ \p\left( |\langle p^*-\hat p , \Phi^*_{k-1} \phi_\mi^k \rangle| \geqslant t \right)  + \p\left( |\langle p^*-\hat p , \phi_\mi^k  \Psi^*_k\rangle| \geqslant t \right).
    \end{split}
\end{equation}    
We now bound the four terms in the double sums above by \Cref{lemma:density Bernstein}, which is based on Bernstein inequality.    
    
We first bound $$\p\left( | \langle p^*-\hat p , \Phi^* _{k-1} \phi_\mi^k  \Psi^*_k\rangle| \geqslant t \right).$$ Note that  by \Cref{lemma:norms of projection basis functions}, it holds that 
    \begin{align*}
        \| \Phi^* _{k-1} \phi_\mi^k  \Psi^*_k  \| _{\lt([0,1]^d)} 
     = 
    \|\Phi^* _{k-1} \|_{\lt([0,1]^{k-1} ) } \|\phi_\mi \|_{\lt([0,1])}  \|\Psi^* _{k+1} \|_{\lt([0,1]^{d-k} ) } =1 ,
    \end{align*}  
    and that 
    \begin{align*}
         \| \Phi^* _{k-1}\phi_\mi^k  \Psi^*_k  \| _{\infty} 
     =  
    \|\Phi^* _{k-1} \|_{\infty}   \|\phi_\mi \|_{\infty}   \|\Psi^* _{k+1} \|_{\infty}  \leqslant C'^2 \sqrt{   n (k-1)(d-k)} .
    \end{align*}  
where $C'= \frac{4B^*}{b^*}$. 
    Therefore, by \Cref{lemma:density Bernstein},
    \begin{equation}
    \begin{split}
   \p\left( | \langle p^*-\hat p , \Phi^* _{k-1} \phi_\mi^k  \Psi^*_k\rangle| \geqslant t \right) \leqslant  & \ 2  \exp\bigg(  \frac{ - Nt^2}{  \|p^*\|_{\infty} \| \Phi^*_{k-1} \phi_\mi^k  \Psi^*_k\|^2_{\lt([0,1]^d)}   +   \| \Phi^*_{k-1} \phi_\mi^k  \Psi^*_k\|_\infty t }\bigg) \\
    \leqslant & \ 2  \exp\bigg(  \frac{ - Nt^2}{  \| p^*\|_\infty    +  C'^2  \sqrt{   n (k-1)(d-k)} t }\bigg) =: \frac{\epsilon}{4dn} .
    \end{split}
    \end{equation}
where we assume the right-hand side is $\frac{\epsilon}{4dn}$ for some sufficiently small constant $\epsilon
$. Solving the last line above yields the quadratic equation 
\begin{displaymath}
    Nt^2 -  C'^2  \sqrt{   n (k-1)(d-k)} \log(\frac{8dn}{\epsilon}) t - \log(\frac{8dn}{\epsilon})\| p^*\|_\infty = 0,
\end{displaymath}
whose larger root can be bounded by 
\begin{align*}
    t_+ = \ &  \frac{C'^2  \sqrt{   n (k-1)(d-k)} \log(\frac{8dn}{\epsilon})  + \sqrt{C'^4 n(k-1)(d-k)[\log(\frac{8dn}{\epsilon})]^2 + 4N \log(\frac{8dn}{\epsilon})\| p^*\|_\infty }  }{2N} \\ 
  \leqslant \ & C_\epsilon \sqrt{\frac{\log(dn)\| p^*\|_\infty}{N}} =:t_\epsilon
\end{align*}
upon choosing sufficiently large $N$. Here $C_\epsilon \sim \log(1/\epsilon)$ is some sufficiently large positive constant which only depends on $\epsilon$. At this point, we have 
\begin{equation*}
  \p\left( | \langle p^*-\hat p , \Phi^* _{k-1} \phi_\mi^k  \Psi^*_k\rangle| \geqslant t_\epsilon \right) \leqslant  \frac{\epsilon}{4dn}.
\end{equation*}
For the other three terms in the double sums, by \Cref{lemma:norms of projection basis functions} we have 
\begin{align*}
    \|\phi_\mi^k\|_\infty \leqslant \sqrt{n}, \quad \|\Phi^* _{k-1} \phi_\mi^k\|_\infty \leqslant \sqrt{   n  (k-1)  }, \quad   \| \phi_\mi^k \Psi^*_k\|_\infty \leqslant \sqrt{   n  (d-k)  }.
\end{align*}
Following the same calculations, we can show that with high probability
\begin{align*}
   & \p\left( \big |  \langle p^*-\hat p , \phi_\mi^k \rangle  \big| \geqslant t_\epsilon \right) \leqslant \frac{\epsilon}{4dn} \\ 
   &\p\left( \big |  \langle p^*-\hat p , \Phi^* _{k-1} \phi_\mi^k \rangle   \big| \geqslant t_\epsilon \right) \leqslant \frac{\epsilon}{4dn} , \\ 
   &\p\left( \big |  \langle p^*-\hat p ,  \phi_\mi^k \Psi^*_k\rangle    \big|  \geqslant  t_\epsilon \right) \leqslant  \frac{\epsilon}{4dn} 
\end{align*}
with appropriate choice of $N$. Inserting the estimates above into \eqref{eq:Bernstein bound for projection term detail}, we conclude that 
\begin{align*}
    \p\left( \max_{k\in\{ 1,\cdots, d\}}\|  ( U_{k-1}^{*T}   \otimes I_n )  (c^{*(k)} -\hat c^{(k)}[s]) V_k ^* \|^2_F  \geqslant 4n t_\epsilon^2\right) \leqslant \epsilon 
\end{align*}
for any given positive $\epsilon$. This indicates that, with high probability, 
\begin{align*}
    \max_{k\in\{ 1,\cdots, d\}}\|  ( U_{k-1}^{*T}   \otimes I_n )  (c^{*(k)} -\hat c^{(k)}[s]) V_k ^* \|_F \leqslant 2\sqrt{n} t_\epsilon =  2 C_\epsilon \sqrt{\frac{n\log(dn)\| p^*\|_\infty}{N}}
\end{align*}

\end{proof}

\begin{lemma}\label{lemma:Frobenius norm s-cluster} Suppose    \Cref{assumption:additive model} holds.
 Let $c^*$ be the coefficient tensor of $p^*$, and $\hat c$ be the $s$-cluster 
 coefficient tensor estimator defined in \Cref{definition:cluster coefficient tensor estimator} with $s\geqslant 2$. 
Then 
\begin{align} \label{eq:Frobenius norm s-cluster}
       \|      c^{* } -\hat c     \|_F    = {\bigO } _\p \bigg(\sqrt { \frac{     d^s n^{2s}    } { N} }    \bigg) .
\end{align}
\end{lemma}
\begin{proof}
        It suffices to show that  
        $$ \E \big( \|      c^{* } -\hat c     \|_F  ^2  \big)  = {\bigO }  \bigg(  \frac{     d^s n^{2s}     } { N}      \bigg). $$
        To this end, let $\mathcal B_{j,d}$ be the $j$-cluster basis functions defined in \eqref{k-cluster basis}. Then by
         \begin{align*}
           \E  (\|      c^{* } -\hat c     \|_F  ^2  )  =
            \sum_{j=0}^s  \bigg \{    \sum_{\varphi_0\in \mathcal B_{1,d}}\E ( \langle \hat p - p^*, \varphi_0\rangle ^2)  + \cdots +   \sum_{\varphi_s\in \mathcal B_{s,d}}\E ( \langle \hat p -p^*  , \varphi_s\rangle  ^2 )   \bigg\}  
        \end{align*}
      Let $g_s(x) =\prod_{k=1}^s \phi_{i_{j_k}}(x_{j_k}) $ be a generic element in $\mathcal B_{s,d}$. Then since $\{ \phi_\mi\}_{\mi=1}^\infty$ are Legendre polynomials, it follows that  with $\|\phi_\mi\|_\infty \leqslant \sqrt \mi$ and therefore $\|g_s\|_\infty  \leqslant n^{s/2}$. Consequently,
        \begin{align*}
            \E   (      \langle  \hat p - p^*  , g _s       \rangle ^2  ) = &\var \bigg(    \frac{1}{N}\sum_{i=1}^N g_s (x^{(i)} ) \bigg)= \frac{\var (      g_s (x^{(1)} )  ) }{N}   \\
             \leqslant & \frac{ \int g_s^2(x) p^*(x) \dd x }{N}  
             \leqslant \frac{  \|g_s\|^2_\infty \int p^*(x)\dd x }{N} \leqslant \frac{n^{s}}{N},
        \end{align*}
where the last equality follows from the fact that $p^*$ is a density and so $  \int p^*(x)dx =1$.
Since the cardinality of $\mathcal B_{s,d}$ is at most $n^sd^s$, it follows that 
$ \sum_{\varphi_s\in \mathcal B_{s,d}}\E ( \langle \hat p -p^*  , \varphi_s\rangle  ^2 )  \leqslant    \frac{d^s n^{2s}}{N}  .  $ 
Since $s$ is a bounded integers, it follows that 
 \begin{align*}
           \E  (\|      c^{* } -\hat c     \|_F  ^2  )  \leqslant \sum_{j=1}^s \frac{d^j n^{2j}}{N}   = \bigO \bigg(  \frac{     d^s n^{2s}    } { N}      \bigg).
        \end{align*}

\end{proof}

\begin{lemma} Let $c^*$ be the coefficient tensor of $p^*$, and $\hat c[s]$ be the $s$-cluster 
 coefficient tensor estimator defined in 
 \Cref{definition:cluster coefficient tensor estimator} with $s\geqslant 2$. 
Let $U^*_k$ and $V_k^*$  be defined as in \eqref{eq:additive column space} and  \eqref{eq:additive row space} respectively. $\Phi^*_k$ and $\Psi^*_k$ are defined as in \eqref{eq:reshape matrix colum basis function} and \eqref{eq:reshape matrix row basis function} respectively. The following equalities hold:
\begin{itemize}[leftmargin=*]
    \item If $s\geqslant 3$,   
    \begin{equation}\label{eq:left and right projection} 
    \begin{split}&  \|  ( U_{k-1}^{*T}   \otimes I_n )  (c^{*(k)} -\hat c[s]^{(k)}) V_k ^* \|_F^2 
         \\ = & \   
      \sum_{\mi=1}^n  \bigg(  \langle p^*-\hat p , \phi_\mi^k \rangle ^2 
       +    \langle p^*-\hat p , \Phi^* _{k-1} \phi_\mi^k \rangle ^2  +    \langle p^*-\hat p ,   \phi_\mi^k \Psi_k^* \rangle ^2   
       +     
         \langle p^*-\hat p , \Phi^* _{k-1} \phi_\mi^k \Psi^*_k\rangle ^2 \bigg) ;
         \end{split}
      \end{equation}
      \item If   $s=2$, 
      \begin{equation}  \label{eq:left and right projection two clusters}
        \begin{split}     \|  ( U_{k-1}^{*T}   \otimes I_n )  (c^{*(k)} -\hat c[s]^{(k)}) V_k ^* \|_F^2   
         =   
      \sum_{\mi=1}^n  \bigg(  \langle p^*-\hat p , \phi_\mi^k  \rangle ^2 
         +    \langle p^*-\hat p , \Phi^* _{k-1} \phi_\mi^k \rangle ^2  +    \langle p^*-\hat p ,   \phi_\mi^k \Psi_k^* \rangle ^2 \bigg).
      \end{split}
      \end{equation}
\end{itemize}  
       
\end{lemma}

\begin{proof} For $\mi \in \{1,\cdots,n \}$, denote 
$\chi_\mi \in \mathbb R^n$ be the $\mi$-th canonical basis. That is,   for $k \in \{1,\cdots,n \}$
$$\chi_{\mi}(k) =\begin{cases}
    0, &\text{if } \mi\not= k  ;
    \\
    1, &\text{if }\mi  = k  .
\end{cases}  $$
Denote $\mathcal T $   the reshaped  $\mathbb R^{2\times n \times 2}$ tensor of  the matrix  $ ( U_{k-1}^{*T}   \otimes I_n )  (c^{*(k)} -\hat c[s]^{(k)}) V_k ^*$. It follows from  \eqref{eq:additive column space},  \eqref{eq:additive row space},
\eqref{eq:reshape matrix colum basis tensor} and \eqref{eq:reshape matrix row basis tensor} that,  for any   $   \mi \in \{ 1,\cdots, n\}    $   the $(a,\mi, b)$-coordinate of $\mathcal T$ satisfies 
    \begin{align*}
       \mathcal T _{a,\mi,b} =\begin{cases}
           \langle c^{* } -\hat c[s]  , \  \bvec{1}^{k-1} \otimes \chi_\mi \otimes \bvec{1}^{d-k} \rangle, &\text{if }a=1, b=1;
           \\
            \langle c^{* } -\hat c[s]  , \  u^*_{k-1} \otimes \chi_\mi \otimes \bvec{1}^{d-k} \rangle, 
            &\text{if }a=2,  b=1;
             \\
            \langle c^{* } -\hat c[s]  , \  \bvec{1}^{k-1} \otimes \chi_\mi \otimes v^*_{k} \rangle, 
            &\text{if }a=1,  b=2;
             \\
            \langle c^{* } -\hat c[s]  , \ u^*_{k-1} \otimes \chi_\mi \otimes v^*_{k} \rangle, 
            &\text{if }a=2,  b=2.
       \end{cases}
    \end{align*} 
By \eqref{eq:reshape matrix colum basis function} and \eqref{eq:reshape matrix row basis function},
\begin{itemize}[leftmargin=*]
 
  \item $\bvec{1}^{k-1} \otimes \chi_\mi \otimes \bvec{1}^{d-k}$
  is the coefficient tensor of $\phi_\mi (x_k) \in \Span \{\mathcal B_{1,d} \}$;
   \item   $ u^*_{k-1} \otimes \chi_\mi \otimes \bvec{1}^{d-k}  $  is the coefficient tensor of $\Phi^* _{k-1} (x_1,\cdots, x_{k-1}) \phi_\mi (x_k) \in \Span \{  \mathcal B_{2,d} \}$;
  \item  $ \bvec{1}^{k-1} \otimes \chi_\mi \otimes v^*_{k}  $  is the coefficient tensor of $  \phi_\mi (x_k)\Psi^* _k(x_{k+1},\cdots, x_d)   \in \Span \{  \mathcal B_{2,d} \}$; 
\item $u^*_{k-1} \otimes \chi_\mi \otimes v^*_{k}$  is the coefficient tensor of $  \Phi^* _{k-1}(x_1,\cdots, x_{k-1})\phi_\mi (x_k)\Psi^* _k(x_{k+1},\cdots, x_d)   \in \Span \{  \mathcal B_{3,d} \}$.
\end{itemize}
{\bf Step 1.} Suppose $s\geqslant 3$. By \Cref{lemma:s cluster coefficient tensor compactible},
\begin{itemize}[leftmargin=*]
  \item   $\langle c^{* } -\hat c[s]  , \  \bvec{1}^{k-1} \otimes \chi_\mi \otimes \bvec{1}^{d-k} \rangle  = \langle p^* -\hat p , \phi_\mi^k \rangle$;
\item   $\langle c^{* } -\hat c[s]  , \ u^*_{k-1} \otimes \chi_\mi \otimes \bvec{1}^{d-k}  \rangle  = \langle p^* -\hat p , \Phi^* _{k-1}  \phi_\mi^k  \rangle$;
\item   $\langle c^{* } -\hat c[s]  , \  \bvec{1}^{k-1} \otimes \chi_\mi \otimes v^*_{k}  \rangle  = \langle p^* -\hat p , \phi_\mi^k \Psi^* _k \rangle$;
\item   $\langle c^{* } -\hat c[s]  , \  u^*_{k-1} \otimes \chi_\mi \otimes v^*_{k}  \rangle  = \langle p^* -\hat p , \Phi^* _{k-1}\phi_\mi^k \Psi^* _k  \rangle$.
\end{itemize}
  This immediately leads to  \eqref{eq:left and right projection}.
  \\
  \\
  {\bf Step 2.} Suppose $s=2$. Note that   
  $    \Phi^* _{k-1}\phi_\mi^k \Psi^* _k   \in \Span \{  \mathcal B_{3,d} \}.$   
  By \Cref{lemma:s cluster coefficient tensor compactible},   $ \langle p^* -\hat p , \Phi^* _{k-1}\phi_\mi^k\Psi^* _k  \rangle = 0$. 
  The rest of the calculations  are the same as   {\bf Step 1}. 
\end{proof}

\begin{lemma} \label{lemma:s cluster coefficient tensor compactible}
    Suppose $ p^* $ satisfies \eqref{eq:additive model}, and $\hat c[s]$ is the $s$-cluster coefficient tensor estimator of $p^*$ in \Cref{definition:cluster coefficient tensor estimator}, where   $2\leqslant s\leqslant d$.  
     Suppose $t \in \{1,\cdots,d \}$  and $F :\mathbb R^d \to \mathbb R$ is any function such that $ F \in  \Span \{  \mathcal B_{t,d}\}$. Let $f \in \mathbb R^{n^d}$ be the coefficient tensor of $F$ as      
     \begin{align*}
     f(\mi_1,\cdots, \mi_d) := \langle F, \phi_{\mi_1}^1  \cdots \phi_{\mi_d}^d \rangle \quad \text{for } 1\leqslant \mi_1, \cdots, \mi_d \leqslant n .
     \end{align*}   
     Then 
     $$ \langle c^* - \hat c[s], f\rangle =\begin{cases}
         0, &\text{if } t>s;
         \\
         \langle p^*-  \hat p  , F\rangle, &\text{if } t\leqslant s.
     \end{cases} $$
\begin{proof}
Observe that by orthogonality of $ \{ \mathcal B_{j,d}\}_{j=1}^d$
\begin{align}\nonumber 
\langle c^* - \hat c[s], f\rangle  = &
\sum_{\mi_1,\cdots, \mi_d =1}^n  (c^* - \hat c[s])(\mi_1,\cdots, \mi_d)\cdot  f(\mi_1,\cdots, \mi_d)  
\\ \nonumber 
= & \sum_{\mi_1,\cdots, \mi_d =1}^n \langle p ^* - \hat  p _s , \  \phi_{\mi_1}^1 \cdots  \phi_{\mi_d}^d  \rangle \langle F,  \  \phi_{\mi_1}^1 \cdots  \phi_{\mi_d}^d  \rangle
\\  \nonumber 
=&
\sum_{\varphi_0 \in \mathcal B_{0,d} } \langle p^* - \hat p_s , \varphi_0\rangle \langle F , \varphi_0\rangle  + \cdots +    \sum_{\varphi_d\in \mathcal B_{d,d}}\langle p^* -\hat p_s, \varphi_d\rangle \langle F , \varphi_d\rangle 
\\ \label{eq:coefficient tensor multiplication compactible}
=& \sum_{\varphi_0 \in \mathcal B_{0,d} } \langle p^* - \hat p_s , \varphi_0\rangle \langle F , \varphi_0\rangle  + \cdots +    \sum_{\varphi_s\in \mathcal B_{s,d}}\langle p^* -\hat p_s, \varphi_s\rangle \langle F , \varphi_s\rangle .
\end{align}
where in the second equality,  $\hat p_s$ is defined in \eqref{eq:cluster basis function estimator}, the last equality follows the observation that $p^* \in \Span\{ \mathcal B_{j,d}\}_{j=0}^1$ and $\hat p_s\in \Span\{ \mathcal B_{j,d}\}_{j=0}^s$  with $s\geqslant 2$. Therefore,
\begin{itemize}
    \item If  $t >s$, then $\langle F, \varphi_s\rangle =0$ when $\varphi_s \in \mathcal B_{s,d}$. So in this case
$ \eqref{eq:coefficient tensor multiplication compactible} =0$;
\item If $t\leqslant s$, then 
 $ \eqref{eq:coefficient tensor multiplication compactible} =\langle p^*-  \hat p  , F\rangle   $ by \Cref{lemma:coefficient multiplication compactible}.
\end{itemize}
 
\end{proof}

\end{lemma}

\begin{lemma}\label{lemma:coefficient multiplication compactible}

    Let $\Omega$ be a measurable set  in arbitrary dimension. Let $\{ \psi_{k }\}_{k=1}^\infty $ be a collection of orthonormal basis in $\lt(\Omega)$  and let $\mathcal M$ be a $m$-dimensional linear  subspace of $\lt(\Omega)$ such that 
    $$ \mathcal M = \Span\{\psi_k \}_{k=1}^m.  $$ Suppose $ F , G \in \mathcal M  $. Denote $\langle F, G\rangle = \int_{\Omega} F(x) G(x) dx.$ Then    
    $$ \langle F, G\rangle = \sum_{k=1}^m \langle F, \psi_k\rangle \langle G, \psi_k\rangle. $$
 
\end{lemma}
\begin{proof} It suffices to 
    observe that 
    \begin{align*}
        \langle F, G\rangle = \sum_{k=1}^\infty \langle F, \psi_k\rangle \langle G, \psi_k\rangle 
        = \sum_{k=1}^m \langle F, \psi_k\rangle \langle G, \psi_k\rangle,
    \end{align*}
    where the first equality follows from the assumption that $\{ \psi_{k }\}_{k=1}^\infty $ is a complete basis, and the second equality follows from the fact that $F\in \mathcal M$, and so $\langle F, \psi_k\rangle=0 $ if $k >m$.
\end{proof}

\begin{lemma} \label{lemma:density Bernstein}  
  Suppose that   $\{ x^{(i)}\}_{i=1}^N $  are    independently sampled from the density $p^*: \mathbb R^d \to \mathbb R ^+$. 
Denote the empirical distribution $\hat{p}(x)=\frac{1}{N}\sum_{i=1}^N \delta_{x^{(i)}}(x) $. Let $\varphi : \mathbb R^d \to \mathbb R $ be any non-random function. Then 
$$ \p (    \big|  \langle\hat p -p  ^*, \varphi\rangle \big| \geqslant  t     ) \leqslant 2 \exp\bigg(  \frac{ -    N t^2}{  \|p ^*\| _\infty \|\varphi\|_{\lt([0,1]^d)} ^2 +  \| \varphi\|_\infty t  } \bigg) . $$
\end{lemma}
\begin{proof} It suffices to apply the Bernstein's inequality to the random variable $  \langle\hat p -p^*   , \varphi\rangle $. Observe  that
$ \langle\hat p -p ^*  , \varphi\rangle  = \frac{1}{N} \sum_{ i=1}^N \big  [    \varphi(x^{(i)}) - \E [  \varphi(x^{(i)}) ] \big]   $. Since  
\begin{align*}
    \E  ( \big  [   \varphi(x^{(i)}) - \E [  \varphi(x^{(i)})] \big ] ^2  )  
    =  &\var  [    \varphi(x^{(i)})            ]   
    \leqslant  \E (  \varphi^2 (x^{(i)}))  
    = \int \varphi^2(x) p^*(x) \dd x      \leqslant  \|\varphi\|^2_{\lt([0,1]^d)} \|p^*\|_\infty 
\end{align*}
    and that $\| \varphi(x^{(i)}) - \E [  \varphi(x^{(i)}) ] \|_\infty \leqslant 2 \|\varphi\|_\infty$, the desired result follows immediately from the Bernstein's inequality.
\end{proof}

\subsection{Orthonormal matrices}
\begin{lemma} \label{lemma:invariant orthogonal matrix multiplication}
    Suppose $M\in \mathbb R^{n\times m}$ and  $    V\in \mathbb  O^{n\times a} $ in the sense that $  V ^T   V =I_a$. Suppose in addition that  
    $    V^T  M =    U \Sigma   W^T $
    is the SVD of $   V^T M$, where $  U \in \mathbb O^{a\times r}$, $\Sigma\in \mathbb R^{r\times r} $ and $   W \in \mathbb O^{r\times m}$. Then 
    \begin{align} \label{eq:svd preserve}  V   V^T M  = (  V   U) \Sigma    W^T
    \end{align}
    is the SVD of $    V     V^T M$. In particular, 
    $\rank (   V^T M) =\rank(   V    V^T M) $  and 
    $\sigma_j (   V^T M) = \sigma_j(   V   
 V^T M) $ for all $1\leqslant j\leqslant \rank (   V^T M) $.
\end{lemma} 
\begin{proof}
    Note that 
    $$ (   V    U)^T   V    U =    U^T    V^T    V    U =   U^T I_a    U=   U^T    U  = I_r.$$
   Therefore  $   V    U \in \mathbb O^{n\times r } $. Consequently, \eqref{eq:svd preserve} holds by the uniqueness of SVD. Since the singular values both $    
 V^T M $ and $   V    V^T M $ are both determined by $\Sigma$, the desired results follow.
\end{proof}
  
\begin{lemma} \label{lemma:left projection preserves matrix structure}
Suppose $M \in \mathbb R^{n\times m}$ such that $\rank(M) =r $. Suppose in addition that   $ U \in \mathbb O^{n\times r}$  and that 
$  \sigma_r(M) > 2\| UU^T M -M \|_F.$ 
Then  
\begin{itemize}[leftmargin=*]
    \item $\rank (UU^T M)=\rank (U^T M) =r $;
    \item $\row(M) = \row( UU^T M)$;
    \item $\sigma_r( U^T M )=\sigma_r(UU^T M ) \geqslant \frac{\sigma_r(M)}{2}$.
      \item $ UU^T$ is the projection matrix onto $\col(UU^T M)$. In other words,  the column vectors of $U $ form  a basis for  $\col(UU^T M)$. 
\end{itemize}

\end{lemma}
\begin{proof}
    Note that $\row( UU^T M ) \subset \row(M) $. In addition, by \Cref{lemma:Weyl},
    $$ \sigma_r( UU^T M) \geqslant \sigma_r( M) - \| UU^T M -M \|_F >0.$$
    So $\rank(UU^T M) \geqslant r$. Therefore $\dim(\row(UU^T M  ) ) \geqslant r$. Since $ \row( UU^T M ) \subset \row(M)$ and 
    $ \dim (\row(M)) = \rank (M) = r$, it follows that 
    $$ \row( UU^T M ) = \row(M) \quad \text{and} \quad \rank (UU^T M) =r .$$
    By \Cref{lemma:invariant orthogonal matrix multiplication},
 $\rank ( U^T M)  =  \rank (UU^T M) =r $.  
    In addition,
    $$ \sigma_r(U^T M ) = \sigma_r(UU^T M ) \geqslant \sigma_r(  M ) - \| UU^T M - M \|_{F} \geqslant \frac{\sigma_r(  M ) }{2}, $$
    where the   equality follows from  \Cref{lemma:invariant orthogonal matrix multiplication}, and the first inequality follows from by \Cref{lemma:Weyl}.
    \\
    \\
    Denote $U= [u_1,\cdots, u_r]$, where $\{u_k\}_{k=1}^r$ are orthogonal. Then $$\col (UU^T M) \subset  \Span\{u_k\}_{k=1}^r .$$ Since   $\dim ( \Span \{u_k\}_{k=1}^r)  = \dim(\col (UU^T M) ) = r $, it follows that 
    $ \Span (\{u_k\}_{k=1}^r)  = \col (UU^T M) $. Consequently   $\{u_k\}_{k=1}^r$ form  a basis for  $ \col (UU^T M) $.
    
\end{proof}

\begin{lemma}
    \label{lemma:kronecker product orthogonality}
    Let $m, n,r,s \in \mathbb Z^+ $ be positive integers. 
    Suppose $U \in \mathbb O^{m\times r}$ and $W \in \mathbb O^{n \times s}$. Then $ U \otimes W \in \mathbb O^{mn\times rs} .$
 
\end{lemma}
\begin{proof}
    It suffices to observe that 
    $$ (U \otimes W )^T ( U \otimes W)  = (U ^T \otimes W  ^T ) ( U \otimes W)= U ^T U  \otimes W^T W =I_r\otimes I_s = I_{rs} . $$
\end{proof}

\begin{lemma} \label{lemma:composition of orthogonal matrices}
Let $m, r,s \in \mathbb Z^+ $ be positive integers. 
    Suppose $U \in \mathbb O^{m\times r}$ and $V \in \mathbb O^{rn \times s}$. Then $ (U\otimes I_n)V \in \mathbb O^{mn\times s} .$
\end{lemma}
\begin{proof}
    It suffices to observe that 
    \begin{align*}
        [  (U\otimes I_n)V ]^T  (U\otimes I_n)V  =  &  V^T (U^T\otimes I_n)   (U\otimes I_n)V  =  V^T (U^T U \otimes I_n )   V   
         \\
         = & V^T (I_r \times I_n ) V = V^T  I_{rn}    V =I_s.
    \end{align*}
\end{proof}

\begin{lemma} \label{lemma:complement of kronecker product}
    Let $l, m,n, r\in \mathbb Z^+$ be positive integers. Let $U\in \mathbb O^{m\times r}$ and $ U_\perp \in \mathbb O^{m\times (n-r)}$ be the orthogonal complement matrix of $U$. Then 
    $$ (U\otimes I_l) _\perp = U_\perp \otimes I_l. $$
\end{lemma}
\begin{proof}
    It suffices to note that $U_\perp \otimes I_l \in \mathbb O^{ml\times (n-r)l}$ and that $\langle U  \otimes I_l,U_\perp \otimes I_l \rangle =0.   $
\end{proof}

\subsection{Matrix perturbation bound}
\begin{theorem}\label{lemma:Weyl}
    Let $A$, $\hat{A}$ and $E$ be three matrices in $\mathbb R^{n\times m}$ and  that  $\hat{A} = A  + E$.
    Then 
    $$ | \sigma_r(A) - \sigma_r(\hat A)| \leqslant \| E\|_{\op}. $$
\end{theorem}
\begin{proof}
    This is the well-known Weyl's inequality for singular values. 
\end{proof}

\begin{theorem}\label{lemma:Wedin}
    Let $A$, $\hat{A}$ and $E$ be three matrices in $\mathbb R^{n\times m}$ and  that  $\hat{A} = A  + E$. Suppose  the SVD of   $A ,\hat{A}$  satisfies
\begin{equation*}
   A  =  \begin{bmatrix}
   U  & U _\perp
    \end{bmatrix}
    \begin{bmatrix}
    \Sigma _1 & 0 \\ 0 & \Sigma _2 
    \end{bmatrix}
    \begin{bmatrix}
        V^{ T } \\ V^{  T}_\perp
    \end{bmatrix},\quad 
       \hat{A} =  \begin{bmatrix}
   \hat{U} & \hat{U}_\perp
    \end{bmatrix}
    \begin{bmatrix}
    \hat{\Sigma}_1 & 0 \\ 0 & \hat{\Sigma}_2 
    \end{bmatrix}
    \begin{bmatrix}
        \hat{V}^ T  \\ \hat{V}^T_\perp
    \end{bmatrix}.
\end{equation*}
Here $r\leqslant \min\{n,m \}$, and  $\Sigma_1,\hat{\Sigma}_1\in \mathbb{R}^{r\times r}$ contain  the leading $r$ singular values of $A$ and $\hat A$ respectively.
If $ \sigma_r(A) -\sigma_{r+1}(A) -\|E\|_{\op}  >0$, then 
\begin{equation*}
    \|U U^{ T } - \hat{U}\hat{U}^T\|_{\op} \leqslant  \frac { 2   \| E\|_{\op}}{\sigma_r(A) -\sigma_{r+1}(A) -\|E\|_{\op} } .
\end{equation*}
\end{theorem}  
\begin{proof}
    This is the well-known Wedin's theorem.
\end{proof}

\begin{lemma} \label{lemma:Davis Khan} Let $A$, $\hat{A}$ and $E$ be three matrices in $\mathbb R^{n\times m}$ and  that  $\hat{A} = A  + E$. Suppose  the SVD of   $A ,\hat{A}$  satisfies
\begin{equation*}
   A  =  \begin{bmatrix}
   U  & U _\perp
    \end{bmatrix}
    \begin{bmatrix}
    \Sigma _1 & 0 \\ 0 & \Sigma _2 
    \end{bmatrix}
    \begin{bmatrix}
        V^{ T } \\ V^{  T}_\perp
    \end{bmatrix},\quad 
       \hat{A} =  \begin{bmatrix}
   \hat{U} & \hat{U}_\perp
    \end{bmatrix}
    \begin{bmatrix}
    \hat{\Sigma}_1 & 0 \\ 0 & \hat{\Sigma}_2 
    \end{bmatrix}
    \begin{bmatrix}
        \hat{V}^ T  \\ \hat{V}^T_\perp
    \end{bmatrix}.
\end{equation*}
Here $r\leqslant \min\{n,m \}$, and  $\Sigma_1,\hat{\Sigma}_1\in \mathbb{R}^{r\times r}$ contain  the leading $r$ singular values of $A$ and $\hat A$ respectively.
  Denote
$$ 
\alpha = \sigma_{\min}(U^{  T }\hat{A}V  ), \  \beta  = \| U  _\perp ^ T  \hat{A}V  _\perp\|_{\op}, \   z_{21} = \|U_\perp^{ T}EV \|_{\op}, \ \text{and} \ z_{12} = \|U^{  T} E V_\perp \|_{\op}, 
$$ 
If $\alpha ^2\geqslant   \beta ^2  + \min \{ z_{21}^2,z_{12}^2 \}   $, then
\begin{equation*}
    \|U U^{ T } - \hat{U}\hat{U}^T\|_{\op} \leqslant  \frac{\alpha z_{21}+  \beta    z_{12}}{ \alpha ^2 -  \beta  ^2 - \min \{z_{21}^2,z_{12}^2 \} } .
\end{equation*}

\end{lemma}
\begin{proof}
   See  Theorem 1 of  \cite{cai2018rate}.
 
\end{proof} 

 \begin{corollary} \label{corollary:Davis Khan 2} Let $A$, $\hat{A}$ and $E$ be three matrices in $\mathbb R^{n\times m}$ and that  $\hat{A} = A  + E$. Suppose $\rank (A) = r$ and $ \sigma_r(A)  
\geqslant 4 \|      E         \|_{\op }  $. Write the SVD of   $A ,\hat{A}$  as
\begin{equation}
   A  =  \begin{bmatrix}
   U  &  U_\perp
    \end{bmatrix}
    \begin{bmatrix}
    \Sigma _1 & 0 \\ 0 & 0
    \end{bmatrix}
    \begin{bmatrix}
        V^{ T } \\ V_\perp
    \end{bmatrix},\quad 
       \hat{A} =  \begin{bmatrix}
   \hat{U} & \hat{U}_\perp
    \end{bmatrix}
    \begin{bmatrix}
    \hat{\Sigma}_1 & 0 \\ 0 & \hat{\Sigma}_2 
    \end{bmatrix}
    \begin{bmatrix}
        \hat{V}^ T  \\ \hat{V}^T_\perp
    \end{bmatrix}
\end{equation}
where $\Sigma_1,\hat{\Sigma}_1\in \mathbb{R}^{r\times r}$ contain  the leading $r$ singular values of $A$ and $\hat A$ respectively. The for any 
  $ \tilde U \in \mathbb O^{n\times r} $, $ \tilde V \in \mathbb O^{m\times r} $ such that 
  $$  \tilde U  \tilde U ^T =  UU^T \quad \text{and} \quad \tilde V  \tilde V ^T =  VV^T ,$$
there exists a universal constant $C$ such that 
\begin{align} \label{eq:Davis Khan}
    \|\tilde U \tilde U^{ T } - \hat{U}\hat{U}^T \|_{\op } \leqslant  C  \bigg(  \frac{ \| E \tilde V \|_{\op }}{ \sigma_r (A)} + \frac{\| E\|_{\op}^2}{\sigma_r^2(A )} \bigg) .
\end{align}
\end{corollary}
\begin{proof}
 Let $   \alpha, \alpha,  z_{21}, z_{12}     $ be defined as in \Cref{lemma:Davis Khan}. 
 Observe that  
 $$ \alpha  =  \sigma_{\min}(U^{  T }\hat{A}V  ) =  \sigma_{r}(U^{  T }\hat{A}V  ) = \sigma_r(UU^{  T }\hat{A}V V^T  )
= \sigma_r(UU^{  T } (A+ E)  V V^T  )
 $$ 
 where the second equality follows from the fact that $ U^T \hat A V \in \mathbb R^{r\times r}$, and the second equality follows from   \Cref{lemma:invariant orthogonal matrix multiplication}.
 So by \Cref{lemma:Weyl}, 
 \begin{align}\label{eq:davis khan variant term 1}
     | \alpha  - \sigma_r(A)| \leqslant \| UU^{  T }   E   V V^T \| _{\op} = \|    E   \tilde V \tilde V^T \| _{\op} =
      \|    E   \tilde V   \| _{\op}
 \end{align}
  In addition, since
  $ \rank (A) =r$, it holds that  $ U^T _\perp A V_\perp = 0$. Therefore,
$$ \beta=  \|U^{  T}_\perp  (A+E ) V_\perp \|_{\op } = \|U_\perp^{  T}  E   V_\perp \|_{\op }\leqslant \|E\|_{\op }.  $$
In addition, we have   that $z_{12} \leqslant \|E\|_{\op}$ and that
$$ z_{21} \leqslant  \| E   V\|_{\op}  = \| E   V V^T \|_{\op}  = \| E   \tilde V  \tilde V^T \|_{\op}  =\| E   \tilde V   \|_{\op} , $$
where the first and the last  equality both  follow  from \Cref{lemma:invariant orthogonal matrix multiplication}.
Consequently
\begin{align*}
      \alpha ^2 - \beta  ^2  -\min\{ z_{21}^2,z_{12}^2  \} \geqslant & \frac{1}{2}\sigma_r^2(A)   - 2 \| E\tilde V \|_{\op }^2  - \|E\|_{\op}^2 - \|E\|_{\op } ^2 
    \\
    \geqslant & \frac{1}{2}\sigma_r^2(A) -4   \| E\|^2 _{\op}   \geqslant \frac{1}{4} \sigma_r^2(A),
\end{align*}    
where the first inequality follows from \eqref{eq:davis khan variant term 1},   the second inequality follows from $   \| E\tilde V \|_{\op } \leqslant \|E\|_{\op}$, and the last inequality follows from  
$ \sigma_r(A) \geqslant 4 \| E\|_{\op}$.
  Therefore by \Cref{lemma:Davis Khan},
   \begin{align*}  &
    \|\tilde U \tilde U^{ T } - \hat{U}\hat{U}^T \|_{\op } =
    \|U U^{ T } - \hat{U}\hat{U}^T\|_{\op} \leqslant  \frac{\alpha z_{21}+  \beta    z_{12}}{ \alpha ^2 -  \beta  ^2 - \min \{z_{21}^2,z_{12}^2 \} }  
    \\
    \leqslant  &  \frac{   ( \sigma_r(A)+   \| E\tilde V\|_{\op} )\| E\tilde V\|_{\op}  +  \|E\|_{\op}^2}{  \frac{1}{4}   \sigma_r^2(A) }  
    \leqslant       \frac{4  \| E \tilde V \|_{\op}}{ \sigma_r (A)} +  \frac{ 8\| E\|_{\op}^2}{\sigma_r^2(A )}  ,
\end{align*}
where $ \| E\tilde V \|_{\op } \leqslant \|E\|_{\op}$ is used in the last inequality. 
\end{proof}

\subsection{Polynomial Approximation Theory}
Let $    C^{  \alphaspace}([0,1])$ denote  the class of functions with continuous derivatives up to order $  \alphaspace \geqslant 1$. 

\begin{theorem} 
    For any $f\in  C^{  \alphaspace}([0,1]) $, there exists a polynomial $p_f $ of degree at most $n$ such that
    \begin{align}
        \label{eq:Jackson theorem}\|f- p _f \|_\infty \leqslant C \| f^{(  \alphaspace)}\|_\infty n^{-  \alphaspace}, 
    \end{align}
    where $C $ is an absolute constant which does not depend on $f$.
\end{theorem}
\begin{proof}
    This is the well-known Jackson’s Theorem for polynomial approximation. See Section 7.2 of \cite{devore1993constructive} for a proof.
\end{proof}
\begin{corollary} \label{corollary:legendre}
     Let  $f\in  C^{  \alphaspace}([0,1]) $.  Denote $\mathcal S_n$ the class of polynomial up degree $n$ on $[0,1]$, and  
  $\mathcal P_n$   the projection operator from $\lt([0,1])$ to $\mathcal S_n$  respect to $\lt$ norm.   Then there exists a   constant $C$ which does not depend on $f$  such that 
    $$\|  f-  \mathcal P_n(f)  \|_{\lt([0,1])} \leqslant  C\| f^{(  \alphaspace)}\|_\infty  n^{-  \alphaspace } .$$
\end{corollary}
\begin{proof}
Let $p_f$ be the polynomial satisfying \eqref{eq:Jackson theorem}. 
    Then
    \begin{align*}
        \|  f-  \mathcal P_n(f)  \|_{\lt([0,1])}  \leqslant   \|  f-  p_f  \|_{\lt([0,1])}  \leqslant  \|f- p _f \|_\infty \leqslant C \| f^{(  \alphaspace)}\|_\infty n^{-  \alphaspace}.
    \end{align*}
\end{proof}

 \begin{theorem} \label{theorem:infty norm approximation of polynomial}
    Let  $f\in  C^{  \alphaspace}([0,1]) $.  Denote $\mathcal S_n$ the class of polynomial up degree $n$ on $[0,1]$, and  
  $\mathcal P_n$   the projection operator from $\lt([0,1])$ to $\mathcal S_n$  respect to $\lt$ norm.   Then there exists an constant $C$  only depending on $f$ such that 
    $$\|  f-  \mathcal P_n(f)  \|_\infty \leqslant  C n^{-  \alphaspace } .$$
\end{theorem}
\begin{proof}
  See Theorem 4.8 of \cite{wang2021much}.
\end{proof}

\end{document}